\documentclass[reqno,11pt]{amsart}

\usepackage[top=2.0cm,bottom=2.0cm,left=3cm,right=3cm]{geometry}
\usepackage{amsthm,amsmath,amssymb,dsfont}
\usepackage{mathrsfs,amsfonts,functan,extarrows,mathtools}
\usepackage[colorlinks]{hyperref}

\usepackage{marginnote}
\usepackage{xcolor}
\usepackage{stmaryrd}
\usepackage{esint}
\usepackage{graphicx}
\usepackage{bbm}
\usepackage{wasysym}
\usepackage{textcomp}
\usepackage{float}
\usepackage{wasysym}
\usepackage[title]{appendix}

\usepackage{tikz}
\usetikzlibrary{calc,decorations.pathreplacing}

\newtheorem{theorem}{Theorem}[section]
\newtheorem{proposition}{Proposition}[section]
\newtheorem{definition}{Definition}[section]

\newtheorem{remark}{Remark}[section]
\newtheorem{lemma}{Lemma}[section]

\numberwithin{equation}{section}
\allowdisplaybreaks

\makeatletter
\newsavebox{\@brx}
\newcommand{\llangle}[1][]{\savebox{\@brx}{\(\m@th{#1\langle}\)}%
 \mathopen{\copy\@brx\kern-0.5\wd\@brx\usebox{\@brx}}}
\newcommand{\rrangle}[1][]{\savebox{\@brx}{\(\m@th{#1\rangle}\)}%
 \mathclose{\copy\@brx\kern-0.5\wd\@brx\usebox{\@brx}}}
\makeatother

\def\d{\mathrm{d}}
\def\no{\nonumber}
\def\R{\mathbb{R}}

\def\exp{\mathrm{exp}}

\def\C{\mathcal{C}}

\def\L{\mathcal{L}}

\def\P{\mathcal{P}}

\newcounter{wronumber}

\begin{document}
\title[SOHB model and limit from SOKB model]
			{Hydrodynamic limits from the self-organized kinetic system for body attitude coordination}

\author[Naping Guo]{Naping Guo}
\address[Naping Guo]
{\newline School of Mathematics, South China University of Technology, Guangzhou, 510641, P. R. China}
\email{Guonping@163.com}

\author[Yi-Long Luo]{Yi-Long Luo${}^*$}
\address[Yi-Long Luo]
{\newline School of Mathematics, Hunan University, Changsha, 410082, P. R. China}
\email{luoylmath@hnu.edu.cn}
\thanks{${}^*$ Corresponding author \quad \today}

\maketitle

\begin{abstract}
   The self-organized kinetic system for body attitude coordination (SOKB) was recently derived by Degond et al. (Math. Models Methods Appl. Sci. 27(6), 1005-1049, 2017). This system describe a new collective motion for multi-agents dynamics, where each agent is described by its position and body attitude: agents travel at a constant speed in a given direction and their body can rotate round it adopting different configurations (representing by rotation matrix in $\mathrm{SO(3)}$). In this paper, we study the hydrodynamic limit of the scaled SOKB system with the constant intensity of coordination by employing the Generalized Collision Invariants (GCI)-based Hilbert expansion approach. The limit is the self-organized hydrodynamic model for body attitude coordination (SOHB). In spherical coordinates, the SOHB system is singular. To avoid this coordinate singularity, we transfer SOHB system into a non-singular form by stereographic projection. This work provides the first analytically rigorous justification of the modeling and asymptotic analysis in Degond et al. (Math. Models Methods Appl. Sci. 27(6), 1005-1049, 2017).\\

   \noindent\textsc{Keywords.} Hydrodynamic limit; SOKB model; SOHB system; Hilbert expansion approach  \\

   \noindent\textsc{Mathematics Subject Classification 2020.} 35B40, 35C20, 35Q70, 35A01, 35L02
\end{abstract}


%


\section{Introduction}

Collective motions of self-propelled particles are ubiquitous in many disciplines, such as in physical, biological or chemical systems. Well-known examples contain birds flocking, fish schooling, animal collective behaviors, human crowds and social dynamics and so on. Such collective behaviors generates complexity interactions among individuals, presenting many nonlinear and nonlocal phenomena, and hence are viewed as complex system consisting of active (living) particles. Describing the complex emergence collective behaviors exhibited in multiscale levels brings new challenges from the viewpoints of phenomenological interpretation, modeling and numerical simulations, and rigorous analysis, and has already attracted more and more attentions from different areas for last two decades. For more introductions, we may refer the readers to, for example, \cite{BBGO-17b,BDT-17b} or review papers \cite{ABF8-19m3as,BF-17m3as,CHL-17chap,KCBFL-13physD,MT-11jsp,Tad-21nAMS}.
    
\subsection{The Self-Organized Kinetic system for body attitude coordination}

In \cite{DFM-2017-MMMAS}, Degond, Frouvelle and Merino-Aceituno derived a new model for multi-agent dynamics where each agent is described by its position and body attitude: agents travel at a constant speed in a given direction and their body can rotate around it adopting different configurations. The corresponding model in mesoscopic level is the so-called Self-Organized Kinetic system for body attitude coordination (in brief, SOKB system), which describes the evolution of the one-agent distribution function $f(t\,,x\,,A)$ at position $x\in\mathbb{R}^{3}$, with the matrix $A\in\mathrm{SO(3)}$ at time $t\geq 0$. More precisely, the dimensionless SOKB model reads
\begin{equation}\label{MFM}
	\begin{aligned}
		\partial_{t}f+ Ae_{1}\cdot\nabla_{x}f=-\nabla_{A}\cdot(fF[f])+d\Delta_{A}f\,,
	\end{aligned}
\end{equation}
\begin{equation}
	\begin{aligned}
		F[f]=\nu P_{T_{A}}(\overline{\mathbb{M}}[f])\,,
	\end{aligned}
\end{equation}
\begin{equation}
	\begin{aligned}
		\overline{\mathbb{M}}[f]=PD(\mathbb{M}[f])\,,\mathbb{M}[f](x,t)=\int_{\mathbb{{R}}^{3} \times\mathrm{SO(3)}}K(x-x')f(x'\,,A'\,,t)A'dA'dx'\,,
	\end{aligned}
\end{equation}
where the constant $d > 0$ is the diffusion coefficient, $\{ e_1, e_2, e_3 \}$ stands for the canonical basis of $\R^3$. The positive coefficient $\nu$, depending on $\overline{\mathbb{M}}[f] \cdot A$, gives the intensity of coordination. The kernel of influence $K$ is Lipschitz, bounded, with the following properties:
\begin{equation*}
	\begin{aligned}
		K = K(|x|) \geq 0 \,, \int_{\R^3} K (|x|) \d x = 1 \,, \int_{\R^3} |x|^2 K (|x|) \d x < \infty \,.
	\end{aligned}
\end{equation*}
We denote by $PD(\mathbb{M}[f]) \in \mathrm{SO(3)}$ the corresponding orthogonal matrix coming from the Polar Decomposition of $ \mathbb{M}[f] $. 

The equation \eqref{MFM} is a Fokker-Planck equation. The left-hand side expresses the rate of change of $f$ due to the spatial transport of the agent with $ Ae_{1}$ while the first term at the right-hand side denotes the transport in velocity space due to the interaction force $F[f]$. The last term at the right-hand side is a velocity diffusion term which arises as a consequence of the Brownian noise. Note that $\Delta_{A}$ is the Laplace-Beltrami operator on $\mathrm{SO(3)}$. 

From the view of mean-field limit, the SOKB model \eqref{MFM} can be regarded as a mean-field limit of time-continuous individual-based model--the Vicsek model. Let $\mathbf{X_{k}} (t) \in \mathbb{R}^{3}$ and $A_{k}(t)\in\mathrm{SO(3)}$ be the position and body attitude of the $k$-th agent at the $t$, where the large number $N$ denotes the total number of agents contained in the system. Then the time-continuous version of the Vicsek model is
\begin{equation}\label{Xk}
    \begin{aligned}
        d \mathbf{X_{k}}(t)=v_{0}A_{k}(t)e_{1}dt\,,
    \end{aligned}
\end{equation}
\begin{equation}
    \begin{aligned}\label{SDE}
        dA_{k}(t)=P_{T_{A_{k}}}\circ [\nu(PD(\mathbb{M}_{k})\cdot A_{k})PD(\mathbb{M}_{k})dt+2\sqrt{D}dW_{t}^{k}]\,,
    \end{aligned}
\end{equation}
\begin{equation}\label{MK}
    \begin{aligned}
        \mathbb{M}_{k}(t)=\frac{1}{N}\sum_{i=1}^{N}K(|\mathbf{X}_{i}-x|)A_{i}(t)\,,
    \end{aligned}
\end{equation}
where $P_{T_{A_{k}}}$ is the projection on the tangent space of $\mathrm{SO(3)}$, which ensures that the resulting solution $A_{k}(t)$ stays in $\mathrm{SO(3)}$, provided that the SDE is taken in the Stratonovich sense (which is indicated by the symbol $\circ$). The equation \eqref{Xk} implies that agents move in the direction of the first axis with constant speed $v_{0}$. The equation \eqref{SDE} takes the form of a stochastic differential equation (SDE). The first term inside the bracket is the interaction. The second term is a noise consisting of independent Brownian motion $W_{t}^{k}$. Then, formally by letting $N\rightarrow \infty$, the \eqref{MFM} system can be obtained. 

\subsection{Hydrodynamic limits from SOKB system \eqref{MFM}}

In the SOKB model \eqref{MFM}, to carry out the macroscopic limit, one rescale the space and time variables by setting $\tilde{x} = \epsilon x$, $\tilde{t} = \epsilon t$, which mean the large-space and long-time scaling, respectively. The parameter $\epsilon$ stands for the relaxation time scale, or equivalently, the asymptotic cruise speed (the balance between friction and self-propulsion), see \cite{MJR-2017,PMJP-MBE}. Then, as shown in Section 4 of \cite{DFM-2017-MMMAS}, after neglecting  the higher order terms $\mathcal{O} (\epsilon)$ and skipping the tildes, the SOKB system \eqref{MFM} indicates the rescaled SOKB formulation
\begin{equation}\label{KM}
    \begin{aligned}
        & \partial_{t}f^{\epsilon}+Ae_{1} \cdot \nabla_{x}f^{\epsilon}=\frac{1}{\epsilon}Q(f^{\epsilon})\,, \\
        & Q(f):=d\Delta_{A}f-\nabla_{A}\cdot(fF_{0}[f]) \,, \\
        & F_{0}[f]=\nu(\Lambda [f]\cdot A)P_{T_{A}}(\Lambda{[f]})\,, \\
        & \Lambda [f]:=PD(\lambda [f])=\lambda [f]\left(\sqrt{(\lambda [f])^\top (\lambda [f])}\right)^{-1}\,, \\
        & \lambda [f](x,t):=\int_{\mathrm{SO(3)}}f(x,A',t)A'dA' \,,
    \end{aligned}
\end{equation}
where $f^{\epsilon}=f^{\epsilon}(t,x,A)$ denotes the one-particle density distribution function in the space $(x,A)\in \mathbb{R}^{3}\times \mathrm{SO(3)}$ at time $t\geq 0$, and $\Lambda[f],Q(f),F_{0}[f]$ are nonlinear operators of $f$, which only act on the attitude variable $A$. 

The main purpose of this paper is to investigate the hydrodynamic limit of \eqref{KM} as $\epsilon \to 0$. In order to illustrate this formal coarse-graining process stated in \cite{DFM-2017-MMMAS}, we first introduce the equilibrium, which are expressed by the von Mises-Fisher (VMF) distribution with respect to the local mean body attitude $\Lambda\in \mathrm{SO(3)}$, namely, a 4-dimension manifold $\mathcal{E}$ given by
\begin{equation}\label{Equilibrium}
	\begin{aligned}
		\mathcal{E}=\left\{\rho M_{\Lambda}(A)|  \rho>0\,,\Lambda\in\mathrm{SO(3)}\right\}\,,
	\end{aligned}
\end{equation}
where $\rho$ is the total mass and the VMF distrubution is defined as
\begin{equation}\label{M-Lambda}
	\begin{aligned}
		M_{\Lambda}(A)=\frac{1}{Z}\exp\left(\frac{\sigma(A\cdot \Lambda)}{d}\right)
	\end{aligned}
\end{equation}
with a normalizing constant $Z=Z(\nu,d)=\int_{\mathrm{SO(3)}}\exp\left(d^{-1}\sigma(A\cdot \mathrm{Id})\right)dA$ and $\sigma =\sigma (\mu )$ is such that 
\begin{equation}\label{sigma-nu}
	\begin{aligned}
		\tfrac{\d}{\d \mu} \sigma (\mu ) = \nu (\mu) > 0 \,.
	\end{aligned}
\end{equation}
Here the smooth function $\nu (\cdot) > 0$ appears in \eqref{KM}. Based on the VMF distribution $M_\Lambda (A)$, we introduce the linear Fokker-Planck type operator
\begin{equation}\label{FP-Oper}
  \begin{aligned}
    \L_{M_\Lambda} f = d \nabla_A \cdot \left[ M_\Lambda \nabla_A \left( \frac{f}{M_\Lambda}  \right) \right] \,.
  \end{aligned}
\end{equation}
The VMF distribution enjoys the following properties:
\begin{enumerate}
	\item[(i)] $M_{\Lambda}(A)$ is a probability density, i.e., $\int_{\mathrm{SO(3)}} M_{\Lambda}(A) \d A=1$;
	\item[(ii)] For any fixed $\Lambda \in \mathrm{SO(3)}$, $\Lambda=\Lambda[\rho M_{\Lambda}]$.
\end{enumerate}
The relation \eqref{sigma-nu} ensures that the collision operator $Q$ can be rewritten as
\begin{equation}\label{Qf}
	\begin{aligned}
		Q(f)=d\nabla_{A}\cdot \left[M_{\Lambda[f]} \nabla_{A}\left(\frac{f}{M_{\Lambda[f]}}\right)\right] = \L_{M_{\Lambda[f]}} f \,,
	\end{aligned}
\end{equation}
which results in a dissipation relation
\begin{equation}
	\begin{aligned}
		H (f) : = \int_{\mathrm{SO(3)}} Q (f) \frac{f}{M_{\Lambda[f]}} \d A = - d \int_{\mathrm{SO(3)}} M_{\Lambda[f]} \left| \nabla_A \left( \frac{f}{M_{\Lambda[f]}}  \right) \right|^2 \d A \leq 0 \,.
	\end{aligned}
\end{equation}
This implies that $Q (f) = 0$ iff $f \in \mathcal{E}$, iff $H (f) = 0$.

Formally, by letting $f^\epsilon (t,x,A) \to f_0 (t,x,A)$ as $\epsilon \to 0$, the scaled equation \eqref{KM} reads
\begin{equation*}
	\begin{aligned}
		f_0 (t,x,A) = \rho (t,x) M_{\Lambda (t,x)} (A) \in \mathcal{E} \,.
	\end{aligned}
\end{equation*}
The aim of hydrodynamic limit for the SOKB system \eqref{KM} is to explore what equations the macroscopic unknowns $(\rho, \Lambda) (t,x)$ obeys. Compared to the hydrodynamic limits of Boltzmann equation or related models, it is completely different. The macroscopic limit equations of Boltzmann equation or related models can be derived from the mass, momentum, energy conservation laws of the kinetic models by using the Collision Invariants (CI) of the corresponding collision operator, see \cite{BGL-CPAM-1993,C-1990-Book,GS-2005-BOOK} for instance. For these models, the all CIs are exactly spanned the equilibrium manifold $\mathcal{E}$. However, for the SOKB model \eqref{MK}, the main challenge of deriving the macroscopic equations is the lack of conservation laws. Actually, there is only one mass conservation law. To overcome this difficulty, the {\em Generalized Collision Invariants} (GCI) are employed in \cite{DFM-2017-MMMAS}. Remark that the concept of GCI was first introduced in Degond-Motsch's work \cite{PS-2008-MMMAS}, which means that we will consider collision invariants not for all $f$ but for those satisfying $P_{T_{\Lambda_0}} (\lambda [f]) = 0$ when one arbitrary $\Lambda_0 \in \mathrm{SO(3)}$ is given. 

For any $\rho_0 > 0$, the linearized operator of the nonlinear collision operator $Q (f)$ around $f_0 = \rho_0 M_{\Lambda_0}$ (denoted by $\L_{M_{\Lambda_0}}^{\mathrm{SB}} f_1$ with $f = f_1 + f_0$) reads
\begin{equation}\label{L-SB}
	\begin{aligned}
		\L_{M_{\Lambda_0}}^{\mathrm{SB}} f_1 = & \L_{M_{\Lambda_0}} \left( f_1 - (c_1 d)^{-1} \nu (A \cdot \Lambda_0) A \cdot P_{T_{\Lambda_0}} ( \lambda [f_1] ) M_{\Lambda_0} \right) \\
= & d \nabla_A \cdot \left\lbrace M_{\Lambda_0} \nabla_A \left( \frac{ f_1 - (c_1 d)^{-1} \nu (A \cdot \Lambda_0) A \cdot P_{T_{\Lambda_0}} ( \lambda [f_1] ) M_{\Lambda_0} }{ M_{ \Lambda_0 } } \right) \right\rbrace \,,
	\end{aligned}
\end{equation}
as in Proposition \ref{Pp2.1} below. Here the Fokker-Planck type operator $\L_{M_{\Lambda_0}}$ is given in \eqref{FP-Oper}. In \cite{DFM-2017-MMMAS}, Degond et al. introduced the definition and existence of GCI related to the linear operator
\begin{equation*}
	\begin{aligned}
		\mathcal{Q} (f, \Lambda_0) = \L_{M_{\Lambda_0}} \left( f - (c_1 d)^{-1} \nu (A \cdot \Lambda_0) A \cdot P_{T_{\Lambda_0}} ( \lambda [f] ) M_{\Lambda_0} \right) \,.
	\end{aligned}
\end{equation*}

For technical simplicity, we consider the constant intensity of coordination in this paper, hence,
\begin{equation}\label{Const-int}
	\begin{aligned}
		\nu (\cdot) = \nu_0 > 0 \,.
	\end{aligned}
\end{equation}
Consequently, the local Maxwellian $M_{\Lambda_0} (A)$ is exactly expressed by
\begin{equation}\label{M-Lambda-0}
	\begin{aligned}
		M_{\Lambda_0} (A) = \tfrac{1}{Z} \exp \left( \tfrac{\nu_0}{d} A \cdot \Lambda_0 \right) \,.
	\end{aligned}
\end{equation}
We remark that the general nonconstant intensity of coordination can also be treated under some proper assumptions such as the positive lower bound and sufficient regularity.

\begin{definition}[Generalized Collision Invariant, \cite{DFM-2017-MMMAS}]
	For a given $\Lambda_0 \in \mathrm{SO(3)}$, the Fokker-Planck type operator $\L_{M_{\Lambda_0}}$ is defined in \eqref{FP-Oper},  i.e., 
	\begin{equation*}
		\begin{aligned}
			\L_{M_{\Lambda_0}} f = \nabla_A \cdot \left\lbrace d M_{\Lambda_0} \nabla_A \left( \frac{f_1}{ M_{ \Lambda_0 } } \right) \right\rbrace \,.
		\end{aligned}
	\end{equation*}
	We say that a real-valued function $\psi : \mathrm{SO(3)} \to \R$ is a Generalized Collision Invariant (in brief, GCI) associated to $\Lambda_0$, or for short $\psi \in GCI (\Lambda_0)$, if
	\begin{equation*}
		\begin{aligned}
			\int_{SO(3)} \L_{M_{\Lambda_0}} f \psi \d A = 0 \ \textrm{ for all } f \textrm{ such that } P_{T_{\Lambda_0}} (\lambda [f]) = 0 \,.
		\end{aligned}
	\end{equation*}
\end{definition}
Actually, there hold
\begin{equation*}
	\begin{aligned}
		GCI (\Lambda_0) = & \{ \psi ; \textrm{ there exists } B \in T_{\Lambda_0} \textrm{ such that } (\L_{M_{\Lambda_0}}^* \psi ) (A) = B \cdot A \} \\
		= & \{ P \cdot (\Lambda_0^\top A) \bar{\psi}_0 (\Lambda_0 \cdot A) + C ; C \in \R ; P \in \mathcal{A} \} \,,
	\end{aligned}
\end{equation*}
where $\mathcal{A}$ stands for the set of antisymmetric matrices, $\L_{M_{\Lambda_0}}^*$ is the adjoint operator of $\L_{M_{\Lambda_0}}$ with the form
\begin{equation*}
	\begin{aligned}
		\L_{M_{\Lambda_0}}^* \psi = M_{\Lambda_0}^{-1} \nabla_A \cdot (M_{\Lambda_0} \nabla_A \psi) \,.
	\end{aligned}
\end{equation*}
Here the function $\bar{\psi}_0 (\cdot)$ is defined by $\tilde{\psi}_0 (\theta) = \bar{\psi}_0 (\frac{1}{2} + \cos \theta)$, where the $2 \pi$-periodic, even and negative function $\tilde{\psi}_0 : \R \to \R^+$ is the unique solution to 
\begin{equation*}
	\begin{aligned}
		\frac{1}{\sin^2 (\theta / 2)} \partial_\theta \big( \sin^2 (\theta / 2) m (\theta) \partial_\theta ( \sin \theta \tilde{\psi}_0 ) \big) - \frac{m (\theta) \sin \theta}{2 \sin^2 (\theta / 2)} \tilde{\psi}_0 = \sin \theta m (\theta) \,, \\
		m (\theta) = Z^{-1} \exp \big( d^{-1} \sigma (\frac{1}{2} + \cos \theta) \big) \,.
	\end{aligned}
\end{equation*}
It is easy to see that $GCI (\Lambda_0)$ is a 4-dimension manifold.

Based on GCI, as in \cite{DFM-2017-MMMAS}, one can derive that $(\rho, \Lambda) (t,x)$ satisfies the following Self-Organized Hydrodynamics for body attitude coordinate (in short, SOHB) model:
\begin{equation}\label{SOHB}
	\left\{
	  \begin{aligned}
	  	& \partial_{t}\rho+c_{1}\nabla_{x}\cdot (\rho\Lambda e_{1})=0 \,, \\
	  	& \rho(\partial_{t}\Lambda+c_{2}((\Lambda e_{1})\cdot \nabla_{x})\Lambda)+[(\Lambda e_{1}) \times(c_{3}\nabla_{x}\rho+c_{4}\rho r_{x}(\Lambda))+c_{4}\rho\delta_{x}(\Lambda)\Lambda e_{1}]_{\times}\Lambda=0 \,, \\
	  	& \Lambda \in \mathrm{SO(3)} \,,
	  \end{aligned}
	\right.
\end{equation}
which governs the dynamics of density $\rho=\rho(t,x):\mathbb{R}^{+}\times\mathbb{R}^{3}\rightarrow\mathbb{R}$ and the matrix of the mean body attitude $\Lambda=\Lambda(t,x)\in \mathrm{SO(3)}$. For a given vector $u$, we introduce the antisymmetric matrix $[u]_\times$, where $[\cdot]_{\times}$ is the linear operator form $\R^3$ to $\mathcal{A}$ given by
\begin{equation}\label{u-times}
	\begin{aligned}
		[u]_\times : = 
		\left[
		\begin{array}{ccc}
			0 & - u_3 & u_2 \\
			u_3 & 0 & - u_1 \\
			- u_2 & u_1 & 0 \\
		\end{array}
		\right] \,,
	\end{aligned}
\end{equation}
so that for any vectors $u, v \in \R^3$, we have $[u]_{\times} v = u \times v$. 

The scalar $\delta_{x}(\Lambda)$ and the vector $r_{x}(\Lambda)$ are the first-order differential operators intrinsic to the dynamics: if $\Lambda(x)=\exp([\mathbf{b}(x)]_{\times})\Lambda(x_{0})$ with $b$ smooth around $x_{0}$ and $\mathbf{b}(x_{0})=0,$ then 
\begin{equation*}
	\begin{aligned}
		\delta_{x}(\Lambda)(x_{0})=\nabla_{x}\cdot \mathbf{b}(x)|_{x=x_{0}}\,,r_{x}(\Lambda)(x_{0})=\nabla_{x}\times\mathbf{b}(x)|_{x=x_{0}}\,,
	\end{aligned}
\end{equation*}
where $\nabla_{x}\times$ is the curl operator. These operators are well defined as long as $\Lambda$ is smooth: as we will see in the next section, we can always express a rotation matrix as $\exp([\mathbf{b}]_{\times})$ for some vector $\mathbf{b}\in\mathbb{R}^{3},$ and this function $\mathbf{b}\rightarrow\exp([\mathbf{b}]_{\times})$ is a local diffeomorphism between a neighborhood of $0\in \mathbb{R}^{3}$,and the identity of $\mathrm{SO(3)}$. This gives a unique smooth representation of $\mathbf{b}$ in the neighborhood of 0 when $x$ is in the neighborhood of $x_{0}$ since then $\Lambda(x)\Lambda(x_{0})^{-1}$ is in the neighborhood of Id.

For completeness, we mention here that the coefficients in the system satisfy
\begin{equation*}
    \begin{aligned}
        c_1 = \tfrac{2}{3} \big< \tfrac{1}{2} + \cos \theta \big>_{m(\theta) \sin^2 ( \sin \theta / 2)} \in (0,1) \,, \ c_2 = \tfrac{1}{5} \langle 2 + 3 \cos \theta \rangle_{ \tilde{m} ( \theta ) \sin^2 ( \theta / 2)} \,, \\
        c_3 = d \big< \nu ( \tfrac{1}{2} + \cos \theta )^{-1} \big>_{ \tilde{m} ( \theta ) \sin^2 ( \theta / 2 ) } \,, \ c_4 = \tfrac{1}{5} \langle 1 - \cos \theta \rangle_{ \tilde{m} ( \theta ) \sin^2 ( \theta / 2 ) } \,,
    \end{aligned}
\end{equation*}
where the notation
\begin{equation*}
	\begin{aligned}
		\big< g ( \theta ) \big>_{ m ( \theta ) \sin^2 ( \sin \theta / 2 ) } = \int_0^{\pi} g ( \theta ) \tfrac{ m ( \theta ) \sin^2 ( \sin \theta / 2 ) }{ \int_0^{\pi} m ( \theta' ) \sin^2 ( \sin \theta' / 2 ) \d \theta' } \d \theta \,, \ \tilde{m} ( \theta ) = \nu ( \tfrac{1}{2} + \cos \theta ) \sin^2 \theta m ( \theta ) \tilde{\psi}_0 ( \theta ) \,.
	\end{aligned}
\end{equation*}

The first equation of the \eqref{SOHB} is the continuity equation for $\rho$, which ensures the mass conservation. The convection velocity is given by $c_{1}\Lambda e_{1} $ and $\Lambda e_{1}$ gives the direction of motion. The second one describes the evolution of $\Lambda$. We remark that every term in the second equation of the \eqref{SOHB} belongs to the tangent space at $\Lambda$ in $\mathrm{SO(3)}$.

The goal of this paper is to rigorously justify the limit from the SOKB system \eqref{MK} to the SOHB equations \eqref{SOHB} under the case for constant intensity of coordination in \eqref{Const-int} as $\epsilon \to 0$. 

\subsection{Historical remarks on self-organized motions}

In this subsection, we will give a brief review of research of self-organized motions, both modeling analysis and rigorous mathematical theory. 

There are many kinds of models related to the self-organized motions. In \cite{VCBCS-95prl}, the so-called Vicsek model describing the large number micro particles' collective motions was first established by Vicsek, Czir\'ok, Ben-Jacob, Cohen and Shochet. Then, in \cite{CS-07ieee,CS-07jjm} Cucker and Smale derived a particles alignment dynamics model, in which the agents tend to align with their neighbors. This phenomenon is also described by the term {\em flocking} \cite{ABF8-19m3as}. In \cite{PS-2008-MMMAS}, Degond and Motsch proposed the self-organized kinetic model and formally derived a self-organized  hydrodynamic system, in which the concept of GCI was first introduced. Moreover, a kinetic version of Cucker-Smale model was proposed by Ha-Tadmor \cite{HT-08krm}. Recently, Degond, Frouvelle and Merino-Aceituno \cite{DFM-2017-MMMAS} proposed a new model for multi-agents dynamics where each agent was described by its position and body attitude: agents travel at a constant speed in a given direction and their body can rotate round it adopting different configurations (representing by rotation matrix in $\mathrm{SO(3)}$). They first gave the Individual Based Model for this dynamics and formally derived its corresponding kinetic and macroscopic equations, which call the SOKB model and SOHB model, respectively. There are also variants of interesting self-propelled (self-organized) models having been developed to characterize more complex interactions such as repulsion, attractions, nematic alignment, suspensions, et al., see \cite{DFL-15arma,DFLMN-14swz,DFAT-18mms,DA-20m3as,DMVY-19jmfm} and the references therein.

The rigorous mathematical results on the self-propelled motions are relatively fewer than the formal results. The first rigorous result can be revisited from Jiang, Xiong and Zhang's work \cite{JXZ-16sima}, in which the GCI-based Hilbert expansion approach was employed to rigorously justify the hydrodynamic limit from self-organized kinetic model of Vicsek type to self-organized hydrodynamic model. Then, in \cite{JLZ-ARMA-2020}, Jiang, Zhang and the second author of this paper has proved the hydrodynamic limit from a kinetic -fluid model coupling of Vicsek-Navier-Stokes model towards the self-organized hydrodynamics and Navier-Stokes equations by GCI-based Hilbert expansion method, in which the local well-posedness of the limit equations was proved by introducing the stereographic projection transform (avoiding the coordinate singularity from the spherical coordinates transform, see \cite{DLMP-2013-MAA,ZJ-NARWA-2017} for instance). We emphasize that, in these works, the notion of GCI employed in \cite{PS-2008-MMMAS} plays a key role in derivation and analysis of the limit regime. The GCI hypothesis, forcing the class of solutions to be constructed always present a first-moment along with the known limit orientation, restricts our the expansion ansatz to a linear case with respect to the limit orientation. Note this remains consistent with the limit regime, due to the basic fact the limit equilibrium lies in the kernel of Fokker-Planck operator, with the exactly same orientation. From a viewpoint of analysis, under the GCI hypothesis, the linear collision operator is reduced to the Fokker-Planck(-type) operator. Very recently, in \cite{JLZ-arXiv-2023} Jiang, Zhang and the second author of this paper proved the hydrodynamic limit from a kinetic Cucker-Smale type model to the self-organized hydrodynamic model by the Hilbert expansion method. In this work, the authors illustrated the concept of GCI in an other way. More precisely, the GCI is not in the kernel of the linearized operator $\L$, but in that of dual operator $\L^*$ ($ \neq \L$). There are also some rigorous results on related models. For example, Figalli and Kang \cite{FK-19apde} proved the limit from the kinetic Cucker-Smale model to the pressureless Euler system with nonlocal alignment.

\subsection{Notations and main results} 

\subsubsection{Notations}

Before we state the main results, we initially introduce some notations. Let $ \ss = ( \ss_1 , \ss_2 , \ss_3 ) \in \mathbb{N}^{3}$ be a multi-index with its length defined as $ | \ss | = \ss_1 + \ss_2 + \ss_3 $. The symbol $\ss' \leq \ss$ means $\ss'_i \leq \ss_i$ for $i = 1,2,3$. Moreover, $\ss' < \ss$ stands for $\ss' \leq \ss$ and $|\ss'| < | \ss |$. We then define here the multi-derivative operator 
$$ \partial_x^{\ss} = \frac{ \partial^{ | \ss | } }{ \partial x_1^{ \ss_1 } \partial x_2^{ \ss_2 } \partial x_3^{ \ss_3 } } \,. $$ 
In addition, the notation $ A \lesssim B $ means that there exists some harmless positive constant $ C > 0 $ such that $ A \leq C B $. For two given matrices $M, N \in \R^{3 \times 3}$, the notation $M \cdot N = \mathrm{tr} (M N^\top) = \sum_{i,j=1}^3 M_{ij} N_{ij}$ stands for the Frobenius inner product of $M$ and $N$. Furthermore, $|M|^2 = M \cdot M$ for $M \in \R^{3 \times 3}$.

We will work in Sobolev spaces with respect to $x \in \R^3$, and in weighted Sobolev spaces with respect to the microscopic velocity variables $A \in \mathrm{SO(3)}$. The spaces $L^p_x := L^p (\R^3) $ for $1 \leq p \leq \infty$ endows with the norms
\begin{equation*}
	\begin{aligned}
		\| f \|_{L^p_x} = \big( \int_{\R^3} |f|^p \d x \big)^\frac{1}{p} \, ( 1 \leq p < \infty ) \,, \quad \| f \|_{L^\infty_x} = \sup_{x \in \R^3} |f (x)| \,.
	\end{aligned}
\end{equation*}
The Sobolev space $H^s_x : = H^s (\R^3)$ is defined by the norm
\begin{equation*}
	\begin{aligned}
		\| f \|_{H^s_x} = \big( \sum_{|\ss| \leq s} \| \partial_x^{\ss} f \|^2_{L^2_x} \big)^\frac{1}{2}
	\end{aligned}
\end{equation*}
for integer $s \geq 0$. In particular, $H^0_x = L^2_x$. The weighted $L^2$-space $L^2_A ( M_{ \Lambda_0 } )$ with respect to the variable $A \in \mathrm{SO(3)}$ is defined by
\begin{equation*}
	\begin{aligned}
		\| f \|_{ L^2_A ( M_{ \Lambda_0 } ) } = \big( \int_{\mathrm{SO}(3)} |f|^2 M_{ \Lambda_0 } \d A \big)^\frac{1}{2} < \infty \,.
	\end{aligned}
\end{equation*} 
Moreover, we will also introduce the mixed weighted space $L^p_x L^2_A ( M_{ \Lambda_0 } )$ by
\begin{equation*}
	\begin{aligned}
		\| f \|_{ L^p_x L^2_A ( M_{ \Lambda_0 } ) } = \big( \int_{\R^3} \big( \int_{\mathrm{SO}(3)} | f |^2 M_{ \Lambda_0 } \d A \big)^\frac{p}{2} \d x \big)^\frac{1}{p} \ (1 \leq p < \infty ) \,, 
	\end{aligned}
\end{equation*}
and
\begin{equation*}
	\begin{aligned}
		\| f \|_{ L^\infty_x L^2_A ( M_{ \Lambda_0 } ) } = \sup_{x \in \R^3} \big( \int_{\mathrm{SO}(3)} | f |^2 M_{ \Lambda_0 } \d A \big)^\frac{1}{2} < \infty \,.
	\end{aligned}
\end{equation*}
Particularly, we denote by $L^2_{x,A} ( M_{ \Lambda_0 } ) = L^2_x L^2_A ( M_{ \Lambda_0 } )$, which endows the norm
\begin{equation*}
	\begin{aligned}
		\| f \|_{ L^2_{x,A} ( M_{ \Lambda_0 } ) } = \big( \iint_{\R^3 \times \mathrm{SO(3)}} |f|^2 M_{ \Lambda_0 } \d A \d x \big)^\frac{1}{2} < \infty \,.
	\end{aligned}
\end{equation*}
At the end, we define the mixed Sobolev space $H^s_x L^2_A ( M_{ \Lambda_0 } )$ by
\begin{equation*}
	\begin{aligned}
		\| f \|_{ H^s_x L^2_A ( M_{ \Lambda_0 } ) } = \big( \sum_{|\ss| \leq s} \| \partial_x^{\ss} f \|^2_{ L^2_{x,A} ( M_{ \Lambda_0 } ) } \big)^\frac{1}{2} < \infty 
	\end{aligned}
\end{equation*}
for $s \geq 0$. Observe that $H^0_x L^2_A ( M_{ \Lambda_0 } ) = L^2_{x,A} ( M_{ \Lambda_0 } ) $.

\subsubsection{Main results}

We first sate our main result of the local well-posedness of the system \eqref{SOHB}.

\begin{theorem}[Well-posedness of \eqref{SOHB} system]\label{LWP}
	Let integer $m \geq 3$ and constant $\rho_* > 0$. Let the function $\rho^{in} (x) > 0$ satisfy 
	\begin{equation}\label{IC-2}
		\begin{aligned}
			\inf_{x \in \R^3} \rho^{in} (x) > 0 \,, \ \rho^{in} - \rho_* \in H^m_x \,.
		\end{aligned}
	\end{equation}
	Given the functions $\phi_{i}^{in}(x) \,, \ \theta_{i}^{in}(x) \in H^m_x$ $(i = 1,2,3)$ such that the unit column vector fields
	\begin{equation}\label{IC-3}
		\begin{aligned}
			\Omega^{in} = ( \tfrac{2\phi_{1}^{in}}{W_{1}^{in}}, \tfrac{2\theta_{1}^{in}}{W_{1}^{in}}, 1 - \tfrac{2}{W_{1}^{in}})^\top \,, \ \mathbf{u}^{in} = ( \tfrac{2\phi_{2}^{in}}{W_{2}^{in}}, \tfrac{2 \theta_{2}^{in}}{W_{2}^{in}}, 1 - \tfrac{2}{W_{2}^{in}})^\top \,, \\ 
			\mathbf{v}^{in} = ( \tfrac{2\phi_{3}^{in}}{W_{3}^{in}}, \tfrac{2\theta_{3}^{in}}{W_{3}^{in}}, 1 - \tfrac{2}{W_{3}^{in}} )^\top \,, \ W_{i}^{in} = 1 + ( \phi_{i}^{in} )^{2} + ( \theta_{i}^{in} )^{2} \,, 
		\end{aligned}
	\end{equation}
	obeying $\Omega^{in} \cdot \mathbf{u}^{in} = \mathbf{u}^{in} \cdot \mathbf{v}^{in} = \mathbf{v}^{in} \cdot \Omega^{in} = 0$, namely, $\Lambda^{in} = (\Omega^{in}, \mathbf{u}^{in}, \mathbf{v}^{in}) \in \mathrm{SO} (3)$, the initial data of the \eqref{SOHB} system is imposed on
	\begin{equation}\label{IC-1}
		\begin{aligned}
			(\rho, \Lambda) (0, x) = (\rho^{in}, \Lambda^{in}) (x) \in \R \times \mathrm{SO(3)}  \,,
		\end{aligned}
	\end{equation}
	where $\Lambda^{in}$ actually enjoys the regularity $\nabla_x \Lambda^{in} \in H^{m-1}_x$. Then there is a time $T > 0$ such that the \eqref{SOHB} system with initial conditions \eqref{IC-1} admits a unique solution $(\rho, \Lambda) (t,x) \in \R \times \mathrm{SO} (3)$ over the time interval $[0, T]$ with the form
	\begin{equation}
		\begin{aligned}
			\Lambda e_1 & = ( \tfrac{2 \phi_1}{W_1}, \tfrac{2 \theta_1}{W_1} , 1 - \tfrac{2}{W_1} )^\top \,, \ \Lambda e_2 = ( \tfrac{2 \phi_2}{W_2}, \tfrac{2 \theta_2}{W_2} , 1 - \tfrac{2}{W_2} )^\top \,, \\ 
			\Lambda e_3 & = ( \tfrac{2 \phi_3}{W_3}, \tfrac{2 \theta_3}{W_3} , 1 - \tfrac{2}{W_3} )^\top \,, \ W_i = 1 + (\theta_i)^2 + (\phi_i)^2 \, (i = 1,2,3) \,.
		\end{aligned}
	\end{equation}
	Moreover, $(\rho, \phi_1, \theta_1, \phi_2, \theta_2, \phi_3, \theta_3)$ subjects to
	\begin{equation}
		\begin{aligned}
			& \inf_{(t,x)\in [0, T] \times \R^3} \rho (t,x) > 0 \,, \\ 
			& \rho (t,x) - \rho_* \,, \phi_i (t,x), \theta_i (t,x) \in L^\infty ( [0,T], H^m_x ) \cap H^1 ( [0,T] ,H^{m-1}_x)
		\end{aligned}
	\end{equation}
	for $i = 1,2,3$, which further follows that
	\begin{equation}
		\begin{aligned}
			\nabla_x \Lambda (t,x) \in L^\infty ([0, T], H^{m-1}_x) \cap H^1 ([0, T], H^{m-2}_x) \,.
		\end{aligned}
	\end{equation}
\end{theorem}

We then display the results of fluid limit from the scaling SOKB equation \eqref{KM} to the macroscopic \eqref{SOHB} system. In this paper, we will use the Hilbert expansion approach to achieve our goal. More precisely, we will seek a solution to \eqref{KM} with special form
\begin{equation}\label{Exp-form}
	\begin{aligned}
		f^\epsilon (t,x,A) = f_0 (t,x,A) + \epsilon f_1 (t,x,A) + \epsilon f_R^\epsilon (t,x,A) \,,
	\end{aligned}
\end{equation}
where $ f_0 (t,x,A) = \rho_0 (t,x) M_{ \Lambda_0 (t, x) } (A) $ with $(\rho_0, \Lambda_0)$ subjecting the SOHB system \eqref{SOHB}, and $f_1 (t,x,A)$ is determined by
\begin{equation}\label{f1}
	\begin{aligned}
		\L_{M_{\Lambda_0}} f_1 = \P_{\L}^\perp \big( \partial_t f_0 + A e_1 \cdot \nabla_x f_0 \big) \,, \quad f_1 = \P_{\L}^\perp f_1 \,.
	\end{aligned}
\end{equation} 
with GCI constraint $P_{T_{\Lambda_0}} (\lambda [f_1]) = 0$. Here $\P_{\L}^\perp = I - \P_{\L}$, where $\P_{\L}$ is the projection from $L^2_A (M_{\Lambda_0})$ to the kernel $\textrm{Ker} (\L_{M_{\Lambda_0}})$. Moreover, the remainder $f_R^\epsilon (t,x,A)$ obeys the equation
\begin{equation}\label{RE}
	\begin{aligned}
		\partial_{t} f_{R}^{\epsilon} + A e_{1} \cdot \nabla_{x} f_{R}^{\epsilon} - \frac{1}{\epsilon} \mathcal{L}_{M_{\Lambda_{0}}} f_{R}^{\epsilon} + \frac{1}{\epsilon} L_R f_R^\epsilon = R (f_1) + \widetilde{Q}(f_{R}^{\epsilon}) \,,
	\end{aligned}
\end{equation}
where the terms $L_R f_R^\epsilon$, $ \widetilde{Q}(f_{R}^{\epsilon}) $ and $ R (f_1) $ are later defined in \eqref{LR-opt}, \eqref{RWQ} and \eqref{RRequation}, respectively. The above formal analysis will be given in Section \ref{Sec:Formal} later.

We now impose the following well-prepared initial data on the SOKB equation \eqref{KM}:
\begin{equation}\label{IC-wp}
	\begin{aligned}
		f^\epsilon ( 0, x, A ) = f^{ \epsilon, in} ( x, A ) : = f_0^{in} ( x , A ) + \epsilon f_1^{in} ( x , A ) + \epsilon f_{R}^{\epsilon, in} (x,A) \,,
	\end{aligned}
\end{equation}
where 
\begin{equation*}
	\begin{aligned}
		f_0^{in} ( x , A ) = \rho_0^{in} (x) M_{ \Lambda_0^{in} (x) } (A) \,.
	\end{aligned}
\end{equation*}
Here $( \rho_0^{in}, \Lambda_0^{in} ) (x)$ is the initial data of $(\rho_0, \Lambda_0)$. $f_1^{in} ( x , A )$ is determined by the same way of \eqref{f1}, just replacing $\Lambda_0$ and $f_0$ by $\Lambda_0^{in}$ and $f_0^{in}$, respectively. Moreover, $ f_{R}^{\epsilon, in} (x,A) $ is the initial data of the remainder equation \eqref{RE}, i.e.,
\begin{equation}\label{IC-RE}
	\begin{aligned}
		f_R^\epsilon (0,x,A) = f_{R}^{\epsilon, in} (x,A) \,.
	\end{aligned}
\end{equation}

Now we precisely state our results.
\begin{theorem}[Hydrodynamic limit from SOKB to SOHB]\label{Hydrodynamic}
	Let integer $s\geq 2$. Assume that $ ( \rho_{0}^{in} \,, \Lambda_{0}^{in} ) (x) \in \R \times \mathrm{SO(3)} $ satisfies the assumptions in Theorem \ref{LWP} with $m = s + 4$, such that
	\begin{itemize}
		\item $(\rho_0, \Lambda_0) (t,x)$ is solved over the time interval $[0, T]$ by Theorem \ref{LWP};
		\item The expanded term $f_1$ can be dominated in terms of $ ( \rho_{0}^{in} \,, \Lambda_{0}^{in} ) $ (see Lemma \ref{Lmm-f1} later).
	\end{itemize}
	We further assume that $ d > \frac{25 \sqrt[4]{3} \nu_0 }{ c_1 \lambda_0 } $, where $\lambda_0 > 0$ is the Poincar\'e constant given in Lemma \ref{CE} later, and
	\begin{equation}\label{IC-fR-bnd}
		\begin{aligned}
			\sup_{\epsilon \in (0,1)} \big( \| \rho_R^{\epsilon , in} \|_{H^s_x} + \| \tfrac{ f_R^{ \epsilon, in } }{ M_{ \Lambda_0^{in} } } - \rho_R^{\epsilon , in} \|_{ H^s_x L^2_A ( M_{ \Lambda_0^{in} } ) } \big) < \infty \,,
		\end{aligned}
	\end{equation}
    where $ \rho_R^{\epsilon , in} = \int_{\mathrm{SO}(3)} f_R^{ \epsilon, in } \d A $.
	
	Then there exists an $\epsilon_{0}>0$ such that, for all $ \epsilon \in ( 0 , \epsilon_0 ) $, the Cauchy problem of \eqref{KM}-\eqref{IC-wp} admits a unique solution $ f^{\epsilon} (t,x,A) \in L^\infty ( [0, T]; H^s_x L^2_A ( M_{ \Lambda_0 } ) )$ with the form \eqref{Exp-form}, where the remainder $ f_R^\epsilon ( t, x, A ) $ enjoys the uniform bound
	\begin{equation}\label{Unf-Bnd}
		\begin{aligned}
			\sup_{t \in [0, T]} \Big( \| \rho_R^\epsilon \|^2_{H^s_x} + \| \tfrac{ f_R^\epsilon }{ M_{ \Lambda_0 } } - \rho_R^\epsilon \|^2_{ H^s_x L^2_A ( M_{ \Lambda_0 } ) } \Big) + \frac{1}{\epsilon} \int_0^T \| \nabla_A ( \tfrac{ f_R^\epsilon }{ M_{ \Lambda_0 } } ) \|^2_{ H^s_x L^2_A ( M_{ \Lambda_0 } ) } \d t \leq C \,.
		\end{aligned}
	\end{equation}
	Here $\rho_R^\epsilon = \int_{\mathrm{SO}(3)} f_R^\epsilon \d A$, and the constant $C > 0$ is independent of $\epsilon$. 
\end{theorem}

\begin{remark}
	Based on Lemma \ref{Lmm-f1} and the uniform bound \eqref{Unf-Bnd}, one easily has
	\begin{equation*}
		\begin{aligned}
			\sup_{t \in [0, T]} \| f^\epsilon - f_0 \|^2_{H^s_x L^2_A ( M_{ \Lambda_0 } )} \leq C \epsilon \to 0 
		\end{aligned}
	\end{equation*}
    as $\epsilon \to 0$. This shows the limit from SOKB equation \eqref{KM} to the SOHB system \eqref{SOHB} with the convergence rate $ \epsilon $.
\end{remark}

\begin{remark}
	The assumption $ d > \frac{25 \sqrt[4]{3} \nu_0 }{ c_1 \lambda_0 } $ means that micro diffusion effect is stronger than the effect of intensity of coordination. The structure of the micro diffusion effect corresponds to the linearized operator $\L_{ M_{ \Lambda_0 } }$, which can offer the diffusive mechanism. However, the effect of intensity of coordination corresponds to the error linear operator $ L_R $ defined in \eqref{LR-opt} later, which will weaken the micro diffusion effect. This coefficient assumption is such that the micro diffusion effect plays a dominant role, compared to the effect of intensity of coordination.
\end{remark}

\subsection{Sketch of proofs and novelties}

In this paper, the main goal is to prove the hydrodynamic limit of the SOKB equation \eqref{KM}, which is the conections to the corresponding fluid equations--SOHB system \eqref{SOHB}. Justifying the hydrodynamic limit rigorously is a huge issue, which has been an active reserch field from late 70's, such as the popular contributions for the Navier-Stokes and Euler limits from the Boltzmann equation (see \cite{ASR-19b,BGL-CPAM-1993,BGL-91jsp,Briant-JDE-2016,Caf-80cpam,GSR-04inv,Guo-2006-CPAM,JL-APDE-2022,JLZ-ARMA-2023,JM-CPAM-2017,JXZ-IUMJ-2018} and the references therein). As stated in \cite{JL-APDE-2022} for instance, there are two types of results in this field:
\begin{enumerate}
	\item First obtaining the solutions fo the scaled kinetic equation with {\em uniform} bounds in the small parameter $\epsilon$, then extracting a subsequence converging to the solutions of the fluid equations as $\epsilon \to 0$.
	
	\item First obtaining the solutions for the limit fluid equations, then constructing a sequence fo special solutions (near the Maxwellians) of the scaled kinetic equation for small parameter $\epsilon > 0$.
\end{enumerate}

Usually, the results of type (1) is harder to be obtained than that of type (2), due to the mixture of small parameter singularity and nonlinearity. In this paper, we thereby try to obtain the result of type (2), in which we will employ the so-called Hilbert expansion approach. The advantage of this method is to separate the small parameter singularity and nonlinearity in the remainder equation. Then we will finish our proof by two steps: 1. well-posedness of the SOHB system \eqref{SOHB}; 2. uniform-in-$\epsilon$ bounds for remainder equation \eqref{RE}.

\subsubsection{Well-posedness of SOHB system}

We first investigate the geometric constraint $\Lambda \in \mathrm{SO(3)}$. By employing the Gr\"onwall inequality, we can prove that the constraint holds provided that the initial data $\Lambda_0^{in} \in \mathrm{SO(3)}$, see Lemma \ref{Initial}. Then we rewrite the form of \eqref{SOHB} as the form \eqref{SOHB-1}. Therefore, by adopting the {\em stereographic projection transform} inspired by \cite{JLZ-ARMA-2020}, we can represent the system \eqref{SOHB-1} as the coordinates form \eqref{SHOB-spt}, see Lemma \ref{Lmm-SOHB-spt}. Compared to the spherical coordinates transform (see \cite{DLMP-2013-MAA,ZJ-NARWA-2017} for instance), the stereographic projection transform can avoid the coordinate singularity. As shown in subsection \ref{subsec:Sym} below, the system \eqref{SHOB-spt} can be expressed by the symmetric form \eqref{U-equ}. Therefore, the Proposition 2.1 in Page 425 of \cite{MET-2011-AMSS} can conclude the results of local well-posedness of the SOHB system \eqref{SOHB} given in Theorem \ref{LWP}.

\subsubsection{Uniform-in-$\epsilon$ bounds for remainder equation \eqref{RE}}

The proof of Theorem \ref{Hydrodynamic} relies on the estimate uniform in small $0 < \epsilon < \epsilon_0$, i.e., the estimate \eqref{UEE-ineq}. This depends on the careful design of the energy and energy dissipation functionals, the control of the singular terms (the terms with $\epsilon$ in their denominators). The general principle is to find sufficient dissipative and decay structures. The key is: these ``good" structures should come from the (micro) kinetic part of the system.

From Lemma \ref{CE} later, the dissipative structure comes from the operator $- \L_{M_{ \Lambda_0 }} f_R^\epsilon$ by multiplying the unknowns $ \tfrac{ f_{R}^{\epsilon} }{ M_{ \Lambda_{0} } } - \rho_R^\epsilon$. Then we can obtain the dissipative structure
\begin{equation*}
	\begin{aligned}
		D_0 (t) = \frac{d}{\epsilon} \| \nabla_A ( \tfrac{f_R^\epsilon}{ M_{ \Lambda_0 } } ) \|^2_{ L^2_{x,A} ( M_{ \Lambda_0 } ) } \,.
	\end{aligned}
\end{equation*}
Moreover, under the key cancellation
\begin{equation*}
	\begin{aligned}
		\int_{\mathrm{SO}(3)} ( \tfrac{ f_R^\epsilon }{ M_{ \Lambda_0 } } - \rho_R^\epsilon ) M_{ \Lambda_0 } \d A = 0 \,,
	\end{aligned}
\end{equation*}
the Poincar\'e inequality in Lemma \ref{CE} tells us that some terms in the right-hand side of energy estimates can be controlled in terms of $\tfrac{1}{\epsilon} \| \tfrac{ f_R^\epsilon }{ M_{ \Lambda_0 } } - \rho_R^\epsilon \|^2_{ L^2_{x,A} ( M_{ \Lambda_0 } ) } $.

Furthermore, $\rho_R^\epsilon$ is exactly the coefficient of the projection $\P_{\L} : L^2_A (M_{ \Lambda_0 }) \to \mathrm{Ker} (\L_{M_{ \Lambda_0 }})$ acting on the remainder $f_R^\epsilon / M_{ \Lambda_0 }$. In this sense, the dissipative structure is from the modulo $L^2_A (M_{ \Lambda_0 }) / \mathrm{Ker} (\L_{M_{ \Lambda_0 }})$-part of $f_R^\epsilon$. As a result, the kernel $\mathrm{Ker} (\L_{M_{ \Lambda_0 }})$-part of $f_R^\epsilon$ should be dominated separately. Inspired by the so-called {\em micro-macro decomposition approach} for Boltzmann equation (see Guo's work \cite{Guo-2006-CPAM}, for instance), we will employ the micro-macro decomposition approach to the SOKB system. In the remainder equation \eqref{RE}, the operators $\L_{ M_{ \Lambda_0 } } f_R^\epsilon $, $ L_R f_R^\epsilon $ and $\widetilde{Q} (f_R^\epsilon)$ are all divergence form with respect to the variable $A \in \mathrm{SO(3)}$, which means that
\begin{equation*}
	\begin{aligned}
		\int_{\mathrm{SO}(3)} \big( \tfrac{1}{\epsilon} \L_{ M_{ \Lambda_0 } } f_R^\epsilon - \tfrac{1}{\epsilon} L_R f_R^\epsilon + \widetilde{Q} (f_R^\epsilon) \big) \d A = 0 \,.
	\end{aligned}
\end{equation*}
In other words, $\tfrac{1}{\epsilon} \L_{ M_{ \Lambda_0 } } f_R^\epsilon - \tfrac{1}{\epsilon} L_R f_R^\epsilon + \widetilde{Q} (f_R^\epsilon) \in L^2_A (M_{ \Lambda_0 }) / \mathrm{Ker} (\L_{M_{ \Lambda_0 }})$. Then, from from projecting the remainder equation \eqref{RE} into $\mathrm{Ker} (\L_{M_{ \Lambda_0 }})$, the macro-equation \eqref{macro-eq} holds, hence,
\begin{equation*}
	\begin{aligned}
		\partial_{t} \rho_{R}^{\epsilon} + \int_{\mathrm{SO(3)}} A e_{1} \cdot \nabla_{x} f_{R}^{\epsilon} \d A = \int_{\mathrm{SO}(3)} R (f_1) \d A \,.
	\end{aligned}
\end{equation*}
Based on the above micro-macro decomposition arguments, we can design a $L^2$ energy
\begin{equation*}
	\begin{aligned}
		E_0 (t) = \| \| \tfrac{ f_R^\epsilon }{ M_{ \Lambda_0 } } - \rho_R^\epsilon \|^2_{ L^2_{x,A} ( M_{ \Lambda_0 } ) } + \| \rho_R^\epsilon \|^2_{ L^2_x } \,.
	\end{aligned}
\end{equation*}

We also have to deal with the error linear operator $\tfrac{1}{\epsilon} L_R f_R^\epsilon$ given in \eqref{LR-opt}, which will weaken the effect of micro diffusion. Thanks to Lemma \ref{CRFM}, one has
\begin{equation*}
	\begin{aligned}
		P_{T_{\Lambda_0}} ( \lambda [f_R^\epsilon - \rho_R^\epsilon M_{\Lambda_0}] ) = P_{T_{\Lambda_0}} ( \lambda [f_R^\epsilon] ) \,,
	\end{aligned}
\end{equation*}
which gives
\begin{equation*}
	\begin{aligned}
		L_R f_R^\epsilon = L_R (f_R^\epsilon - \rho_R^\epsilon M_{\Lambda_0}) \,.
	\end{aligned}
\end{equation*}
Then, under the coefficient assumption $ d > \frac{25 \sqrt[4]{3} \nu_0 }{ c_1 \lambda_0 } $, the Poincar\'e inequality 
$$ \| \nabla_{A} ( \tfrac{ f_{R}^{\epsilon} }{ M_{ \Lambda_{0} } } ) \|^2_{L^2_A (M_{\Lambda_0})} \geq \lambda_{0} \| \tfrac{ f^{\epsilon}_{R} }{ M_{ \Lambda_{0} } } - \rho_R^\epsilon \|^{2}_{ L^2_A ( M_{ \Lambda_{0} } )}$$
guarantees that $\tfrac{1}{\epsilon} L_R f_R^\epsilon$ can be absorbed by $- \tfrac{1}{\epsilon} \L_{ M_{ \Lambda_0 } } f_R^\epsilon $. As a result, we obtain the $L^2$ estimate \eqref{L2-bnd}, i.e.,
\begin{equation*}
	\begin{aligned}
		\tfrac{\d}{\d t} E_0 (t) + d_\star D_0 (t) \lesssim \textrm{ some controllable terms by higher order derivatives} \,.
	\end{aligned}
\end{equation*}

However, the above $L^2$ estimate is not closed. We need to derive the higher order spatial derivatives estimates. For $\ss \in \mathbb{N}^3$ with $1 \leq |\ss| = k \leq s$ ($ s \geq 2 $), multiplying \eqref{Hmicro-eq} by $\partial_x^{\ss} ( \tfrac{ f_R^\epsilon }{ M_{ \Lambda_0 } } - \rho_R^\epsilon ) M_{ \Lambda_0 } $ and combining with the macro-equation \eqref{Hmacro-eq}, i.e.,
\begin{equation*}
	\begin{aligned}
		\partial_{t} \partial_x^{\ss} \rho_{R}^{\epsilon} + \int_{\mathrm{SO(3)}} A e_{1} \cdot \nabla_{x} \partial_x^{\ss} f_{R}^{\epsilon} \d A = \int_{\mathrm{SO}(3)} \partial_x^{\ss} R (f_1) \d A \,,
	\end{aligned}
\end{equation*}
we can obtain the energy functional with $k$-th order spatial derivatives
\begin{equation*}
	\begin{aligned}
		\| \partial_x^{\ss} ( \tfrac{ f_R^\epsilon }{ M_{ \Lambda_0 } } - \rho_R^\epsilon ) \|^2_{L^2_{x,A} ( M_{ \Lambda_0 } ) } + \| \partial_x^{\ss} \rho_R^\epsilon \|^2_{L^2_x} \,,
	\end{aligned}
\end{equation*}
see Step 1 and 2 of Subsection \ref{Subsec:DerEst} later. The dissipative structure comes from the term $ - \frac{1}{\epsilon} \partial_x^{\ss} \big[ \tfrac{1}{ M_{ \Lambda_0 } } \L_{ M_{ \Lambda_0 } } f_R^\epsilon \big] $, which gives the dissipation
\begin{equation*}
	\begin{aligned}
		\frac{d}{\epsilon} \| \nabla_A \partial_x^{\ss} ( \tfrac{ f_R^\epsilon }{ M_{ \Lambda_0 } } ) \|^2_{ L^2_{x,A} ( M_{ \Lambda_0 } ) } \,.
	\end{aligned}
\end{equation*}
But there are two remainders $\C_2$ and $\C_3$ required to be dominated, see \eqref{C1C2C3}. Specially, there is a norm $\| \partial_x^{\ss} ( \tfrac{ f_R^\epsilon }{ M_{ \Lambda_0 } } - \rho_R^\epsilon ) \|^2_{ L^2_{x,A} ( M_{ \Lambda_0 } ) }$ which should be controlled in terms of $\| \nabla_A \partial_x^{\ss} ( \tfrac{ f_R^\epsilon }{ M_{ \Lambda_0 } } ) \|^2_{ L^2_{x,A} ( M_{ \Lambda_0 } ) }$. Unfortunately, it cannot directly apply the Poincar\'e inequality in Lemma \ref{CE}, due to 
$$ \int_{\mathrm{SO}(3)} \partial_x^{\ss} ( \tfrac{ f_R^\epsilon }{ M_{ \Lambda_0 } } - \rho_R^\epsilon ) M_{ \Lambda_0 } \d A \neq 0 \,. $$
Thanks to the cancellation $ \int_{\mathrm{SO}(3)} \partial_x^{\ss} [ ( \tfrac{ f_R^\epsilon }{ M_{ \Lambda_0 } } - \rho_R^\epsilon ) M_{ \Lambda_0 } ] \d A = 0 $ for any $\ss \in \mathbb{N}^3$, we can obtain the Poincar\'e type inequality \eqref{HPoincare} for higher order spatial derivatives, i.e.,
\begin{equation*}
	\begin{aligned}
		\lambda_0 \| \partial_x^{\ss} ( \tfrac{ f_R^\epsilon }{ M_{ \Lambda_0 } } - \rho_R^\epsilon ) \|_{ L^2_{x,A} ( M_{ \Lambda_0 } ) }^2 \leq & \| \nabla_A \partial_x^{\ss} ( \tfrac{ f_R^\epsilon }{ M_{ \Lambda_0 } } ) \|_{ L^2_{x,A} ( M_{ \Lambda_0 } ) }^2 \\
		& + C ( \Lambda_0 ) \sum_{0 \neq \ss' \leq \ss} \| \nabla_A \partial_x^{\ss - \ss'} ( \tfrac{ f_R^\epsilon }{ M_{ \Lambda_0 } } ) \|^2_{ L^2_{x,A} ( M_{ \Lambda_0 } ) } \,,
	\end{aligned}
\end{equation*}
which satisfies our requirements.

Moreover, under the coefficient assumption $ d > \frac{25 \sqrt[4]{3} \nu_0 }{ c_1 \lambda_0 } $ and the Poincar\'e type inequality \eqref{HPoincare}, the error linear operator with higher order derivatives $ \frac{1}{\epsilon} \partial_x^{\ss} \big[ \tfrac{1}{ M_{ \Lambda_0 } } L_R ( f_R^\epsilon - \rho_R^\epsilon M_{ \Lambda_0 } ) \big] $ can be successfully absorbed by $ - \frac{1}{\epsilon} \partial_x^{\ss} \big[ \tfrac{1}{ M_{ \Lambda_0 } } \L_{ M_{ \Lambda_0 } } f_R^\epsilon \big] $. Therefore, we gain the $k$-th order derivative estimate \eqref{Hk-bnd}, i.e.,
\begin{equation*}
	\begin{aligned}
		\tfrac{\d}{\d t} E_k (t) + d_\star D_k (t) \lesssim \sum_{0 \leq \mathfrak{j} \leq k-1} D_{ \mathfrak{j} } (t) + \textrm{ some controllable terms} \,.
	\end{aligned}
\end{equation*}
In order to control the quantity $ \sum_{0 \leq \mathfrak{j} \leq k-1} D_{ \mathfrak{j} } (t) $, based on \eqref{L2-bnd} and \eqref{Hk-bnd}, the induction arguments for $k = 0, 1, \cdots, s$ can be applied to finish our uniform-in-$\epsilon$ bounds. At the end, the continuity arguments can finish the proof of Theorem \ref{Hydrodynamic}.

\subsection{Organization of current paper}

In the next section, the formal analysis of the Hilbert expansion are carried, in which the remainder equation is written down. Section \ref{Sec:LWP} gives the proof of well-posedness to the SOHB system \eqref{SOHB}. Section \ref{Sec:UnifEst} carries out the uniform-in-$\epsilon$ estimates of the remainder equation and finish the proof of Theorem \ref{Hydrodynamic}. Section \ref{Sec:f1}, the bounds for the expanded term $f_1$ are derived, i.e., proof of Lemma \ref{Lmm-f1}.

\section{Formal Analysis from the Hilbert Expansion Method}\label{Sec:Formal}

In this section, we formally study the macroscopic limit of the kinetic model \eqref{KM} by using Hilbert expansion. The key is to compute the leading term, the higher order expanded terms and truncate the expansion in a proper order. Then the equation of remainders should be write down. Before carrying out this, we first do some preparations.

\subsection{Preparations}

\begin{lemma}[Consistency relation for the ``flux", Lemma 4.8 of \cite{DFM-2017-MMMAS}]\label{CRFM}
\begin{equation}\label{CRF}
    \begin{aligned}
        \lambda [M_{\Lambda_{0}}]=c_{1}\Lambda_{0}\,,
    \end{aligned}
\end{equation}
    where $c_{1}\in (0,1)$ is explicitly expressed as
\begin{equation*}
    \begin{aligned}
        c_1 = \tfrac{2}{3} \big< \tfrac{1}{2} + \cos \theta \big>_{ m ( \theta ) \sin^2 ( \sin \theta / 2 )}
    \end{aligned}
\end{equation*}
for $m(\theta)=\exp(d^{-1}\sigma(\frac{1}{2} + \cos \theta ))$. Moreover, for any fixed $(\rho_0 , \Lambda_0) \in \R_+ \times \mathrm{SO(3)}$,
\begin{equation}
	\begin{aligned}
		\Lambda_0 = \Lambda [\rho_0 M_{\Lambda_0}] \,.
	\end{aligned}
\end{equation}
\end{lemma}

\begin{lemma}[Projection operator on the tangent space, Proposition A.3 of \cite{DFM-2017-MMMAS}]\label{Lmm-Proj-SO3}
	Let $A \in \mathrm{SO(3)}$ and $M \in \mathcal{M}$ (set of square matrices). Let $P_{T_A}$ be the orthogonal projection on $T_A$ (tangent space at $A$), then
	\begin{equation}
		\begin{aligned}
			P_{T_A} (M) = \frac{1}{2} (M - A M^\top A) = \nabla_A (A \cdot M) \,.
		\end{aligned}
	\end{equation}
    Notice that then
    \begin{equation}
    	\begin{aligned}
    		P_{T_A^\perp} (M) = \frac{1}{2} (M + A M^\top A) \,.
    	\end{aligned}
    \end{equation}
\end{lemma}

Then we study the properties of $\Lambda [f]$ defined in \eqref{KM}.

\begin{lemma}\label{Lmm-Lambda}
	Let $f_0 = \rho_0 M_{\Lambda_0} \in \mathcal{E}$, where $\rho_0 > 0$ and $\Lambda_0 \in \mathrm{SO(3)}$ are both any fixed. Then  
	\begin{equation}
		\begin{aligned}
			\tfrac{\d}{\d \epsilon}|_{\epsilon=0} \Lambda[f_0 + \epsilon f_1 ] = & (c_1 \rho_0 )^{-1} P_{T_{\Lambda_0}} ( \lambda[f_1]) \,, \\
			\tfrac{\d^2}{\d \epsilon^2} \big|_{\epsilon = 0} \Lambda [ f_0 + \epsilon f_1 ] = & - \tfrac{1}{2 ( c_1 \rho_0 )^2} \lambda [ f_1 ] ( \lambda [ f_1 ]^\top \Lambda_0 + \Lambda_0^\top \lambda [ f_1 ] ) \\
			& - \tfrac{1}{( c_1 \rho_0 )^3} \lambda [ f_1 ]^\top \lambda [ f_1 ] + \tfrac{3}{4 ( c_1 \rho_0 )^3} \big( \lambda [ f_1 ]^\top \Lambda_0 + \Lambda_0^\top \lambda [ f_1 ] \big)^2 \,,
		\end{aligned}
	\end{equation}
    where the constant $c_1$ is given in Lemma \ref{CRFM}. 
\end{lemma}

\begin{proof}
	Recall that
	\begin{equation}
		\begin{aligned}
			\Lambda[f] = PD (\lambda [f]) = \lambda[f] \left( \sqrt{\lambda[f]^\top \lambda[f]} \right)^{-1} \,, \ \lambda[f] = \int_{\mathrm{SO(3)}} f(x,A',t) A' \d A'\,.
		\end{aligned}
	\end{equation}
    A direct calculation implies
    \begin{equation*}
    	\begin{aligned}
    		\tfrac{ \d }{ \d \epsilon } \Lambda [ f_{0} + \epsilon f_{1} ] = & \lambda [ f_{1} ] V [ f_0 + \epsilon f_1 ] + \lambda [ f_{0} + \epsilon f_{1} ] \tfrac{ \d }{ \d \epsilon } V [ f_0 + \epsilon f_1 ] \,,
    	\end{aligned}
    \end{equation*}
    where
    \begin{equation*}
    	\begin{aligned}
    		V [ f_0 + \epsilon f_1 ] = \left( \sqrt{ \lambda [ f_{0} + \epsilon f_{1} ]^\top \lambda [ f_{0} + \epsilon f_{1} ] } \right)^{-1} \,.
    	\end{aligned}
    \end{equation*}
    Moreover,
	\begin{equation}\label{D1}
		\begin{aligned}
			& \tfrac{ \d }{ \d \epsilon } V [ f_0 + \epsilon f_1 ] = - \frac{1}{2} \left( V [ f_0 + \epsilon f_1 ] \right)^2 \left( \lambda [ f_{1} ]^\top \lambda [ f_{0} + \epsilon f_{1} ] + \lambda [ f_{0} + \epsilon f_{1} ]^\top \lambda [ f_{1} ] \right) V [ f_0 + \epsilon f_1 ] \,.
		\end{aligned}
	\end{equation}
	Then
	\begin{equation}\label{D2}
		\begin{aligned}
			\tfrac{ \d }{ \d \epsilon } \Lambda [ f_{0} + \epsilon f_{1} ] = & \lambda [ f_{1} ] V [ f_0 + \epsilon f_1 ] \\
			& - \tfrac{1}{2} \lambda [ f_{0} + \epsilon f_{1} ] \left( V [ f_0 + \epsilon f_1 ] \right)^2 \\
			& \times \left( \lambda [ f_{1} ]^\top \lambda [ f_{0} + \epsilon f_{1} ] + \lambda [ f_{0} + \epsilon f_{1} ]^\top \lambda [ f_{1} ] \right) V [ f_0 + \epsilon f_1 ] \,.
		\end{aligned}
	\end{equation}
	
	So, together with the fact $\lambda [f_0] = c_1 \rho_0 \Lambda_0$ derived from Lemma \ref{CRFM}, we can get
	\begin{equation}
		\begin{aligned}
			& \tfrac{\d}{\d\epsilon}|_{\epsilon=0}\Lambda [f_{0}+\epsilon f_{1}]\\
			=&\lambda[f_{1}]\left(\sqrt{\lambda[f_{0}]^\top\lambda[f_{0}]}\right)^{-1}-\frac{1}{2}\lambda[f]\left(\sqrt{\lambda[f_{0}]^\top\lambda[f_{0}]}\right)^{-2}\\
			&\times\left(\lambda[f_{1}]^\top\lambda[f_{0}]+\lambda[f_{0}]^\top\lambda[f_{1}]\right)\left(\sqrt{\lambda[f_{0}]^\top\lambda[f_{0}]}\right)^{-1}\\
			=&-\frac{1}{2}c_{1}\rho_{0}\Lambda_{0}\frac{1}{(c_{1}\rho_{0})^{2}}\left(\lambda[f_{1}]^\top c_{1}\rho_{0}\Lambda_{0}+c_{1}\rho_{0}\Lambda_{0}^\top\lambda[f_{1}]\right)(c_{1}\rho_{0})^{-1} +(c_{1}\rho_{0})^{-1}\Lambda[f_{1}]\\
			=&(c_{1}\rho_{0})^{-1}\frac{1}{2}(\lambda[f_{1}]-\Lambda_{0}\lambda[f_{1}]^\top\Lambda_{0})\\
			=&(c_{1}\rho_{0})^{-1}P_{T_{\Lambda_{0}}}(\lambda[f_{1}])\,,
		\end{aligned}
	\end{equation}
    where the last equality is deduced from Lemma \ref{Lmm-Proj-SO3}. 
    
    By the relations \eqref{D1} and \eqref{D2}, it is easy to derive that
    \begin{equation}
    	\begin{aligned}
    		& \tfrac{\d^2}{\d \epsilon^2} \Lambda [ f_0 + \epsilon f_1 ] \\
    		= & - \tfrac{1}{2} \lambda [ f_1 ] ( V [ f_0 + \epsilon f_1 ] )^2 ( \lambda [f_1]^\top \lambda [ f_0 + \epsilon f_1 ] + \lambda [ f_0 + \epsilon f_1 ]^\top \lambda [ f_1 ] ) V [ f_0 + \epsilon f_1 ] \\
    		& + \tfrac{1}{2} ( V [ f_0 + \epsilon f_1 ] )^3 \Big\{ ( \lambda [ f_1 ]^\top \lambda [ f_0 + \epsilon f_1 ] + \lambda [ f_0 + \epsilon f_1 ]^\top \lambda [ f_1 ] ) V [ f_0 + \epsilon f_1 ] \Big\}^2 \\
    		& - ( V [ f_0 + \epsilon f_1 ] )^2 ( \lambda [ f_1 ]^\top \lambda [ f_1 ] ) V [ f_0 + \epsilon f_1 ]  \\
    		& + \tfrac{1}{4} \Big\{ ( V [ f_0 + \epsilon f_1 ] )^2 ( \lambda [ f_1 ]^\top \lambda [ f_0 + \epsilon f_1 ] + \lambda [ f_0 + \epsilon f_1 ]^\top \lambda [ f_1 ] ) \Big\}^2 V [ f_0 + \epsilon f_1 ] \,.
    	\end{aligned}
    \end{equation}
    Then by $\lambda [ f_0 ] = c_1 \rho_0 \Lambda_0$ and $\Lambda_0 \in \mathrm{SO(3)}$, one has
    \begin{equation*}
    	\begin{aligned}
    		V [ f_0 ] = ( \sqrt{ \lambda [ f_0 ]^\top \lambda [ f_0 ] } )^{-1} = ( \sqrt{ c_1 \rho_0 \Lambda_0^\top c_1 \rho_0 \Lambda_0 } )^{-1} = \tfrac{1}{c_1 \rho_0} I_3 \,.
    	\end{aligned}
    \end{equation*}
    It thereby follows that
    	\begin{align*}
    		\tfrac{\d^2}{\d \epsilon^2} \big|_{\epsilon = 0} \Lambda [ f_0 + \epsilon f_1 ] = & - \tfrac{1}{2} \lambda [ f_1 ] \tfrac{1}{( c_1 \rho_0 )^2} \big( \lambda [ f_1 ]^\top c_1 \rho_0 \Lambda_0 + c_1 \rho_0 \Lambda_0^\top \lambda [ f_1 ] \big) \tfrac{1}{c_1 \rho_0} \\
    		& + \tfrac{1}{2} \tfrac{1}{( c_1 \rho_0 )^3} \Big\{ \big( \lambda [ f_1 ]^\top c_1 \rho_0 \Lambda_0 + c_1 \rho_0 \Lambda_0^\top \lambda [ f_1 ] \big) \tfrac{1}{c_1 \rho_0} \Big\}^2 \\
    		& - \tfrac{1}{(c_1 \rho_0)^2} ( \lambda [ f_1 ]^\top \lambda [ f_1 ] ) \tfrac{1}{c_1 \rho_0} \\
    		& + \tfrac{1}{4} \Big\{ \tfrac{1}{(c_1 \rho_0)^2} \big( \lambda [ f_1 ]^\top c_1 \rho_0 \Lambda_0 + c_1 \rho_0 \Lambda_0^\top \lambda [ f_1 ] \big) \Big\}^2 \tfrac{1}{c_1 \rho_0} \\
    		= & - \tfrac{1}{2 ( c_1 \rho_0 )^2} \lambda [ f_1 ] ( \lambda [ f_1 ]^\top \Lambda_0 + \Lambda_0^\top \lambda [ f_1 ] ) \\
    		& - \tfrac{1}{( c_1 \rho_0 )^3} \lambda [ f_1 ]^\top \lambda [ f_1 ] + \tfrac{3}{4 ( c_1 \rho_0 )^3} \big( \lambda [ f_1 ]^\top \Lambda_0 + \Lambda_0^\top \lambda [ f_1 ] \big)^2 \,.
    	\end{align*}
    The proof of Lemma \ref{Lmm-Lambda} is thereby completed.
\end{proof}

Next we investigate the linearized operator of $Q (f)$ around the equilibrium $f_0 = \rho_0 M_{\Lambda_0}$.

\begin{proposition}\label{Pp2.1}
	Assume $\sigma (\mu)$ satisfies \eqref{sigma-nu}, i.e., $\frac{\d}{\d \mu} \sigma (\mu) = \nu (v) > 0$. Then the linearized collision operator $$\mathcal{L}^{\mathrm{SB}}_{M_{\Lambda_{0}}}f_{1}=\frac{d}{d\epsilon}\big|_{\epsilon=0}Q(f_{0}+\epsilon f_{1})$$
	defined by the first variation of $Q$ with respect to $f_{1}$ can be expressed by
	\begin{equation}\label{Lf_{1}}
		\begin{aligned}
			\L_{M_{\Lambda_0}}^{\mathrm{SB}} f_1 = & \L_{M_{\Lambda_0}} \left( f_1 - (c_1 d)^{-1} \nu (A \cdot \Lambda_0) A \cdot P_{T_{\Lambda_0}} ( \lambda [f_1] ) M_{\Lambda_0} \right) \\
= & d \nabla_A \cdot \left\lbrace M_{\Lambda_0} \nabla_A \left( \frac{ f_1 - (c_1 d)^{-1} \nu (A \cdot \Lambda_0) A \cdot P_{T_{\Lambda_0}} ( \lambda [f_1] ) M_{\Lambda_0} }{ M_{ \Lambda_0 } } \right) \right\rbrace \,,
		\end{aligned}
	\end{equation}
    where the constant $c_1 > 0$ is given in Lemma \ref{CRFM}. Namely, $\mathcal{L}^{\mathrm{SB}}_{M_{\Lambda_{0}}}$ is a Fokker-Planck type operator defined on $\mathrm{SO(3)}$.
\end{proposition}

\begin{proof}
	Recalling the expressions of $M_\Lambda$ in \eqref{M-Lambda} and $Q (f)$ in \eqref{Qf}, one has
	\begin{equation}\label{E1}
		\begin{aligned}
			& \tfrac{\d}{\d \epsilon} \big|_{\epsilon=0} Q( f_{0} + \epsilon f_{1})
			\\
			= & d \nabla_{A} \cdot \left( \tfrac{\d}{\d \epsilon} \big|_{\epsilon = 0} \left[ M_{\Lambda [f_{0} + \epsilon f_{1} ] } \nabla_{A} \left( \frac{f_{0} + \epsilon f_{1}}{ M_{ \Lambda [ f_{0} + \epsilon f_{1} ] } } \right) \right] \right) \\
			= & d \nabla_{A} \cdot \left \lbrace
			M_{ \Lambda [f_{0}] } \nabla_{A} \left( \frac{f_{1}}{ M_{ \Lambda [f_{0}] } } \right) \right. \\
			& + M_{ \Lambda [f_{0}] } A \cdot \tfrac{\d}{\d \epsilon} \big|_{ \epsilon = 0 } \Lambda [ f_{0} + \epsilon f_{1} ] d^{-1} \nu ( A \cdot \Lambda [f_{0}] ) \nabla_{A} \left( \frac{ f_{0} }{ M_{ \Lambda [f_{0}] } } \right) \\
			& \left. - M_{ \Lambda [f_{0}] } \nabla_{A} \left( \frac{ f_{0} }{ M_{ \Lambda [f_{0}] } } A \cdot \tfrac{\d}{\d \epsilon} \big|_{ \epsilon = 0 } \Lambda [ f_{0} + \epsilon f_{1} ] d^{-1} \nu ( A \cdot \Lambda [f_{0}] ) \right) \right \rbrace \,.
		\end{aligned}
	\end{equation}
Observe that $f_0/ M_{\Lambda[f_0]} = \rho_0$, which means that
\begin{equation}\label{f0-0}
  \begin{aligned}
    \nabla_{A} \left( \frac{ f_{0} }{ M_{ \Lambda [f_{0}] } } \right) = 0 \,.
  \end{aligned}
\end{equation}
Thanks to Lemma \ref{CRFM} and \ref{Lmm-Lambda}, one has
	\begin{equation}\label{E2}
		\begin{aligned}
			\Lambda [f_{0}]=\Lambda [\rho_{0}M_{\Lambda_{0}}]=\Lambda_{0}
		\end{aligned}
	\end{equation}
	and
	\begin{equation}\label{E3}
		\begin{aligned}
			\tfrac{\d}{\d \epsilon}|_{\epsilon=0} \Lambda[f_0 + \epsilon f_1 ] = (c_1 \rho_0 )^{-1} P_{T_{\Lambda_0}} ( \lambda[f_1]) \,.
		\end{aligned}
	\end{equation}
	Then the relations \eqref{E1}, \eqref{f0-0}, \eqref{E2} and \eqref{E3} conclude the equality \eqref{Lf_{1}}. The proof of Proposition \ref{Pp2.1} is therefore finished.
\end{proof}

\subsection{Hilbert asymptotic expansion for SOKB system \eqref{KM}}

\subsubsection{{Orders analysis}}

We aim here to justify the asymptotic limit by performing a Hilbert expansion into the SOKB equation:
\begin{equation}\label{f}
	\begin{aligned}
		f^{\epsilon} = f_{0} + \epsilon f_{1} + \epsilon^{2} f_{2} + \cdots \,.
	\end{aligned}
\end{equation}
Then by the orders analysis, we can formally derive the limit equation and the required higher order terms. By plugging \eqref{f} into \eqref{KM}, one has
\begin{equation*}
	\begin{aligned}
		( \partial_{t} + A e_{1} \cdot \nabla_{x} ) ( f_{0} + \epsilon f_{1} + \epsilon^{2} f_{2} + \cdots ) = \frac{1}{\epsilon} Q ( f_{0} + \epsilon f_{1} + \epsilon^{2} f_{2} + \cdots ) \,.
	\end{aligned}
\end{equation*}
Recalling that the definition of $Q (f)$ in \eqref{KM}, one has
\begin{equation}\label{Q1}
    \begin{aligned}
         Q ( f^{\epsilon} ) = & d \Delta_{A} f^{\epsilon} - \nabla_{A} \cdot ( f^{\epsilon} F_{0} [ f^{\epsilon} ] ) \\
         = & d \Delta_{A} ( f_{0} + \epsilon f_{1} + \epsilon^{2} f_{2} + \cdots ) \\
         & - \nabla_{A} \cdot ( ( f_{0} + \epsilon f_{1} + \epsilon^{2} f_{2} + \cdots ) F_{0} [ f_{0} + \epsilon f_{1} + \epsilon^{2} f_{2} + \cdots ] ) \,.
    \end{aligned}
\end{equation}
It therefore infers
    \begin{align}\label{KM-1}
        \no & ( \partial_{t} + A e_{1} \cdot \nabla_{x} ) ( f_{0} + \epsilon f_{1} + \epsilon^{2} f_{2} + \cdots ) \\
        = & \frac{1}{\epsilon} \Big\{ d \Delta_{A} ( f_{0} + \epsilon f_{1} + \epsilon^{2} f_{2} + \cdots ) \\
        \no & - \nabla_{A} \cdot ( ( f_{0} + \epsilon f_{1} + \epsilon^{2} f_{2} + \cdots ) F_{0} [ f_{0} + \epsilon f_{1} + \epsilon^{2} f_{2} + \cdots ] ) \Big\} \,.
    \end{align}
Under the assumption \eqref{Const-int},
\begin{equation*}
	\begin{aligned}
		F_0 [f] = \nu_0 \nabla_A (A \cdot \Lambda [f]) \,.
	\end{aligned}
\end{equation*}
Then, by Lemma \ref{Lmm-Lambda} and the Taylor expansion, one has
\begin{equation}\label{F0}
	\begin{aligned}
		F_{0} [ f_{0} + \epsilon f_{1} + \epsilon^{2} f_{2} + \cdots ] = & F_{0} [ f_{0} ] + \epsilon \nu_0 \nabla_A \Big( \tfrac{\d}{\d \epsilon} \big|_{\epsilon = 0} \Lambda [f_0 + \epsilon f_1] \cdot A \Big) + \mathcal{O} (\epsilon^2) \\
		= & F_{0} [ f_{0} ] + \epsilon \nu_0 (c_1 \rho_0 )^{-1} \nabla_A \Big( P_{T_{\Lambda_0}} ( \lambda[f_1]) \cdot A \Big) + \mathcal{O} (\epsilon^2) \,.
	\end{aligned}
\end{equation}

{\bf Order $\mathcal{O}( \epsilon^{-1} )$ in \eqref{KM-1}:} From \eqref{KM-1} and \eqref{F0}, it follows 
\begin{equation}\label{0Qf0}
    \begin{aligned}
        0 = Q ( f_{0} ) = d \Delta_{A} f_{0} - \nabla_{A} \cdot ( f_{0} F_{0} [ f_{0} ] ) \,.
    \end{aligned}
\end{equation}
Then $f_0 \in \mathcal{E}$ defined in \eqref{Equilibrium}, hence, 
\begin{equation}
    \begin{aligned}
        f_{0} ( t , x , A ) = \rho_{0} ( t, x ) M_{\Lambda_{0} (t,x)} (A) 
    \end{aligned}
\end{equation}
for some parameters $(\rho_0 , \Lambda_0) (t,x) \in \R_+ \times \mathrm{SO(3)}$ to be determined in the next order.

{\bf Order $\mathcal{O}( \epsilon^0 )$ in \eqref{KM-1}:} By \eqref{KM-1} and \eqref{F0}, one has
\begin{equation*}
	\begin{aligned}
		\partial_t f_0 + A e_1 \cdot \nabla_x f_0 = & d \Delta_A f_1 - \nabla_A \cdot (f_1 F_0 [f_0]) \\
		& - \nabla_A \cdot (f_0 \nu_0 (c_1 \rho_0)^{-1} \nabla_A ( P_{T_{\Lambda_0}} (\lambda [f_1]) \cdot A ) ) \\
		= & \L_{M_{\Lambda_0}} f_1 - \tfrac{\nu_0}{c_1 d} d \nabla_A \cdot \left[ M_{\Lambda_0} \nabla_A \left( \frac{P_{T_{\Lambda_0}} (\lambda [f_1]) \cdot A M_{\Lambda_0}}{M_{\Lambda_0}} \right) \right] \\
		= & \L_{M_{\Lambda_0}} \big( f_1 - \tfrac{\nu_0}{c_1 d} P_{T_{\Lambda_0}} (\lambda [f_1]) \cdot A M_{\Lambda_0} \big) \,,
	\end{aligned}
\end{equation*}
hence,
\begin{equation}\label{f0-expand}
	\begin{aligned}
		\partial_t f_0 + A e_1 \cdot \nabla_x f_0 = \L_{M_{\Lambda_0}} \big( f_1 - \tfrac{\nu_0}{c_1 d} P_{T_{\Lambda_0}} (\lambda [f_1]) \cdot A M_{\Lambda_0} \big) \,.
	\end{aligned}
\end{equation}
From the arguments in Theorem 4.16 of \cite{DFM-2017-MMMAS}, one can deduce from integrating the equation \eqref{f0-expand} over $\mathrm{SO(3)}$ that
\begin{equation*}
	\begin{aligned}
		\partial_t \rho_0 + c_1 \nabla_x \cdot (\rho_0 \Lambda_0 e_1) = 0 \,.
	\end{aligned}
\end{equation*}
Denote by $\varphi_P^{\Lambda_0} (A) = P \cdot (\Lambda_0^\top A) \bar{\psi}_0 (\Lambda_0 \cdot A)$ with $P \in \mathcal{A}$ (the set of antisymmetric matrices), which is the non-constant GCI under the constraint $P_{T_{\Lambda_0}} (\lambda [f_1]) = 0$. It then infers from multiplying \eqref{f0-expand} by $\varphi_P^{\Lambda_0} (A)$ and integrating the resultant equation over $A \in \mathrm{SO(3)} $ that
\begin{equation*}
	\begin{aligned}
		\rho_0 ( \partial_t \Lambda_0 + c_2 \Lambda_0 e_1 \cdot \nabla_x \Lambda_0 ) + [(\Lambda_0 e_1) \times ( c_3 \nabla_x \rho_0 + c_4 \rho_0 r_x (\Lambda_0) ) + c_4 \rho_0 \delta_x (\Lambda_0) \Lambda_0 e_1]_\times \Lambda_0 = 0 \,.
	\end{aligned}
\end{equation*}
In summary, the functions $(\rho_0, \Lambda_0) (t,x)$ subject to the SOHB system \eqref{SOHB}.

\subsubsection{{Truncation and the remainder equation}}

In the Hilbert expansion, we hope the number of expanded terms is as small as possible, because the more expanded terms there are, the more special form of the solution becomes. Notice that the limit equations \eqref{SOHB} is derived from the equation \eqref{f0-expand}. So we only require to find a proper $f_1$ such that the equation \eqref{f0-expand} holds.

Note that $\L_{M_{\Lambda_0}}$ is a Fokker-Planck operator. Once $f_0 = \rho_0 M_{\Lambda_0}$ is solved by the SOHB system \eqref{SOHB}, it is easy to know that there is a $f_1$ such that the linear equation \eqref{f0-expand} holds with the constraint $P_{T_{\Lambda_0}} (\lambda [f_1]) = 0$. In the previous sense, the expanded term $f_1$ in \eqref{f} is only partially determined. As a result, we take the truncated form
\begin{equation}\label{f-truncated}
	\begin{aligned}
		f^\epsilon = f_0 + \epsilon f_1 + \epsilon f_R^\epsilon \,,
	\end{aligned}
\end{equation}
where $f_0 = \rho_0 M_{\Lambda_0}$ satisfies \eqref{0Qf0} with $(\rho_0, \Lambda_0)$ solving \eqref{SOHB}, and $f_1$ obeys the equation \eqref{f1}, i.e.,
\begin{equation*}
	\begin{aligned}
		\L_{M_{\Lambda_0}} f_1 = \P_{\L}^\perp \big( \partial_t f_0 + A e_1 \cdot \nabla_x f_0 \big) \,, \quad f_1 = \P_{\L}^\perp f_1 \,, \quad P_{T_{\Lambda_0}} (\lambda [f_1]) = 0 \,.
	\end{aligned}
\end{equation*}

We then insert the truncated expansion \eqref{f-truncated} into the SOKB system \eqref{KM}. Together with \eqref{0Qf0} and \eqref{f0-expand}, one has
\begin{equation}
	\begin{aligned}
		& \epsilon (\partial_t f_R^\epsilon + A e_1 \cdot \nabla_x f_R^\epsilon) + \epsilon (\partial_t f_1 + A e_1 \cdot \nabla_x f_1) \\
		= & d \Delta_A f_R^\epsilon - \nu_0 \nabla_A \cdot \big[ f_R^\epsilon \nabla_A (A \cdot \Lambda [f_0 + \epsilon f_1 + \epsilon f_R^\epsilon]) \big] \\
		& - \nu_0 \nabla_A \cdot \big[ f_0 \nabla_A (A \cdot \tfrac{\d}{\d \epsilon}|_{\epsilon = 0} \Lambda [f_0 + \epsilon g]_{g = f_R^\epsilon} ) \big] \\ 
		& - \nu_0 \nabla_A \cdot \big[ f_1 \nabla_A \{ \Lambda [f_0 + \epsilon f_1 + \epsilon f_R^\epsilon] - \Lambda [f_0] \} \big] \\
		& - \tfrac{1}{\epsilon} \nu_0 \nabla_A \cdot \big[ f_0 \nabla_A \big( A \cdot \{ \Lambda [f_0 + \epsilon f_1 + \epsilon f_R^\epsilon] - \Lambda [f_0] - \epsilon \tfrac{\d}{\d \epsilon} |_{\epsilon = 0} \Lambda [f_0 + \epsilon g]_{g = f_1 + f_R^\epsilon} \} \big) \big] \,.
	\end{aligned}
\end{equation}
By Proposition \ref{Pp2.1}, one has
\begin{equation*}
	\begin{aligned}
		\L_{M_{\Lambda_0}} f_R^\epsilon = d \Delta_A f_R^\epsilon - \nu_0 \nabla_A \cdot [ f_R^\epsilon \nabla_A (A \cdot \Lambda [f_0]) ] \,.
	\end{aligned}
\end{equation*}
Moreover, Lemma \ref{Lmm-Lambda} indicates that
\begin{equation}
	\begin{aligned}
		\tfrac{\d}{\d \epsilon} |_{\epsilon = 0} \Lambda [f_0 + \epsilon g] = (c_1 \rho_0)^{-1} P_{T_{\Lambda_0}} (\lambda [g]) \,.
	\end{aligned}
\end{equation}
As a consequence, the remainder $f_R^\epsilon$ satisfies the equation \eqref{RE}, i.e.,
\begin{equation*}
	\begin{aligned}
		\partial_{t} f_{R}^{\epsilon} + A e_{1} \cdot \nabla_{x} f_{R}^{\epsilon} - \frac{1}{\epsilon} \mathcal{L}_{M_{\Lambda_{0}}} f_{R}^{\epsilon} + \frac{1}{\epsilon} L_R f_R^\epsilon = R (f_1) + \widetilde{Q}(f_{R}^{\epsilon})
	\end{aligned}
\end{equation*}
where the error linear operator $L_R f_R^\epsilon$ reads
\begin{equation}\label{LR-opt}
	\begin{aligned}
		L_R f_R^\epsilon = \tfrac{\nu_0}{c_1} \nabla_A \cdot \big[ M_{\Lambda_0} \nabla_A \big( A \cdot P_{T_{\Lambda_0}} ( \lambda [f_R^\epsilon] ) \big) \big] \,,
	\end{aligned}
\end{equation}
and the nonlinear term $\widetilde{Q}(f_{R}^{\epsilon})$ is expressed by
\begin{equation}\label{RWQ}
	\begin{aligned}
		\widetilde{Q} ( f_{R}^{\epsilon} ) = & - \frac{1}{\epsilon} \nu_{0} \nabla_{A} \cdot [ f_{R}^{\epsilon} \nabla_{A} ( A \cdot \Lambda [ f_{0} + \epsilon f_{1} + \epsilon f_{R}^{\epsilon} ] - A \cdot \Lambda [f_0] ) ] \\
		&- \frac{1}{\epsilon} \nu_0 \nabla_A \cdot [ f_1 \nabla_A ( A \cdot \Lambda [f_0 + \epsilon f_1 + \epsilon f_R^\epsilon] - A \cdot \Lambda [f_0] ) ] \\ 
		& - \frac{1}{\epsilon^2} \nu_0 \nabla_A \cdot \big[ f_0 \nabla_A \big( A \cdot \{ \Lambda [f_0 + \epsilon f_1 + \epsilon f_R^\epsilon] - \Lambda [f_0] - \epsilon \tfrac{\d}{\d \epsilon} |_{\epsilon = 0} \Lambda [f_0 + \epsilon g]_{g = f_1 + f_R^\epsilon} \} \big) \big] \,.
	\end{aligned}
\end{equation}
Moreover, the source term $R (f_1)$ in \eqref{RE} is
\begin{equation}\label{RRequation}
	\begin{aligned}
		R (f_1) = - \partial_{t} f_{1} - A e_{1} \cdot \nabla_{x} f_{1} \,.
	\end{aligned}
\end{equation}
We remark that the nonlinear term $\widetilde{Q} ( f_{R}^{\epsilon} )$ actually does not involve the singularity of parameter $\frac{1}{\epsilon}$ in the view of Taylor expansion.

\section{Well-posedness of SOHB system \eqref{SOHB}: Proof of Theorem \ref{LWP}}\label{Sec:LWP}

In this section, we mainly justify the local existence of the SOHB system \eqref{SOHB} stated in Theorem \ref{LWP}.

\subsection{Geometric constraint $\Lambda \in \mathrm{SO(3)}$}

In this subsection, we will show that the solution $\Lambda (t,x)$ to the $\Lambda$-equation in \eqref{SOHB} will be restricted onto the manifold $\mathrm{SO(3)}$ provided that the corresponding initial data is onto the manifold $\mathrm{SO(3)}$. The spirit can be first found in Eells-Sampson's work \cite{Eells-Sampson-AJM-1964}, in which the authors proved that the heat flow was restricted the same target manifold as the initial data. Moreover, in the works \cite{JLL-JDE-2023,JL-SIMA-2019,JL-JFA-2022,JLT-M3AS-2019}, the similar spirit to deal with the geometric constraint $|d| = 1$ for the hyperbolic Ericksen-Leslie liquid crystal has been employed. More precisely, our result is stated as follows.

\begin{lemma}\label{Initial}
	Assume that the sufficiently smooth matrix-valued function $\Lambda (t,x)$ satisfies the $\Lambda$-equation in \eqref{SOHB}
	\begin{equation}\label{Lambda-eq}
		\begin{aligned}
			\rho(\partial_{t}\Lambda+c_{2}((\Lambda e_{1})\cdot \nabla_{x})\Lambda)+[(\Lambda e_{1}) \times(c_{3}\nabla_{x}\rho+c_{4}\rho r_{x}(\Lambda))+c_{4}\rho\delta_{x}(\Lambda)\Lambda e_{1}]_{\times}\Lambda=0
		\end{aligned}
	\end{equation}
	with initial data $\Lambda (0, x) = \Lambda^{in} (x)$. If the initial data $\Lambda^{in} (x) \in \mathrm{SO(3)}$ and $\rho > 0$, then the solution $\Lambda (t,x) \in \mathrm{SO(3)}$ holds for any $t \geq 0$.
\end{lemma}

\begin{proof}
	Let $\{ e_1, e_2, e_3 \}$ be the canonical basis of $\R^3$. Then the equation \eqref{Lambda-eq} reads 
	\begin{equation}\label{A4}
		\begin{aligned}
			\rho ( \partial_{t} + c_{2} \Lambda e_{1} \cdot \nabla_{x} ) \Lambda e_{i} + {\bf w} \times ( \Lambda e_{i} ) = 0 
		\end{aligned}
	\end{equation}
    for $i = 1,2,3$, where we have used the relation $[ {\bf w} ]_{\times} \Lambda e_{i} = {\bf w} \times ( \Lambda e_{i} )$ with ${\bf w} = ( \Lambda e_{1} ) \times ( c_{3} \nabla_{x} \rho + c_{4} \rho r_{x} ( \Lambda ) ) + c_{4} \rho \delta_{x} ( \Lambda ) \Lambda e_{1} \in \R^3$.
	
	From multiplying \eqref{A4} by $\Lambda e_{i}$, it follows that
	\begin{equation}\label{ii-eq}
		\begin{aligned}
			\rho ( \partial_{t} + c_{2} \Lambda e_{1} \cdot \nabla_{x} ) ( \tfrac{1}{2} | \Lambda e_{i} |^{2} )=0.
		\end{aligned}
	\end{equation}
	If one further dot \eqref{A4} with $\Lambda e_{j} \, ( i \neq j)$, then one gets 
	\begin{equation}\label{ij-eq}
		\begin{aligned}
			\rho ( \partial_{t} + c_{2} \Lambda e_{1} \cdot \nabla_{x} ) \Lambda e_{i} \cdot \Lambda e_{j} + ( {\bf w} \times \Lambda_{i} ) \cdot \Lambda e_{j} = 0 \,.
		\end{aligned}
	\end{equation}
	Similarly, for $i \neq j$, one has
	\begin{equation}\label{ji-eq}
		\begin{aligned}
			\rho ( \partial_{t} + c_{2} \Lambda e_{1} \cdot \nabla_{x} ) \Lambda e_{j} \cdot \Lambda e_{i} + ( {\bf w} \times \Lambda_{j} ) \cdot \Lambda e_{i} =0 \,.
		\end{aligned}
	\end{equation}
	Note that $( {\bf w} \times \Lambda_{i} ) \cdot \Lambda e_{j} + ( {\bf w} \times \Lambda_{j} ) \cdot \Lambda e_{i} = 0$ for $i \neq j$. It therefore follows from \eqref{ij-eq} and \eqref{ji-eq} that
	\begin{equation}\label{ij-equ}
		\begin{aligned}
			\rho ( \partial_{t} + c_{2} \Lambda e_{1} \cdot \nabla_{x} ) ( \Lambda e_{i} \cdot \Lambda e_{j} ) = 0 \,, \ \forall \, 1 \leq i \neq j \leq 3 \,.
		\end{aligned}
	\end{equation}
	Let 
	$$ f_{ij} = \Lambda e_{i} \cdot \Lambda e_{j} - \delta_{ij} $$
	for $1 \leq i, j \leq 3$. Then \eqref{ii-eq} and \eqref{ij-equ} indicate that $f_{ij}$ satisfies
	\begin{equation}\label{fij-eq}
		\begin{aligned}
			\partial_{t} f_{ij} + c_{2} \Lambda e_{1} \cdot \nabla_{x} f_{ij} = 0 \,.
		\end{aligned}
	\end{equation}
	Multiplying the equation \eqref{fij-eq} by $f_{ij}$ and integrating by parts over $x \in \mathbb{R}^{3}$, one has
	\begin{equation*}
		\begin{aligned}
			\tfrac{1}{2} \tfrac{\d}{\d t} \int_{\mathbb{R}^{3}} |f_{ij}|^{2} \d x & = - c_{2} \int_{\mathbb{R}^{3}} \Lambda e_{1} \cdot \nabla_{x} f_{ij} \cdot f_{ij} \d x \\
			& = \tfrac{c_{2}}{2} \int_{\mathbb{R}^{3}} ( \nabla_{x} \cdot (\Lambda e_{1}) ) |f_{ij}|^{2} \d x \\
			& \leq \tfrac{c_{2}}{2} \| \nabla_{x} \cdot (\Lambda e_{1}) \|_{L^{\infty}} \int_{\mathbb{R}^{3}} |f_{ij}|^{2} \d x \,.
		\end{aligned}
	\end{equation*}
	Then the Gr\"onwall inequality implies
	\begin{equation}\label{Gronwall-ineq}
		\begin{aligned}
			0 \leq \| f_{ij} (t, \cdot) \|^2_{L^2} \leq \| f_{ij} (0, \cdot) \|^2_{L^2} \exp \left( \int_{0}^{t} c_{2} \| \nabla_{x} \cdot (\Lambda e_1) ( \tau, \cdot ) \|_{L^\infty} \d \tau \right)
		\end{aligned}
	\end{equation}
	for all $t \geq 0$. Due to $\Lambda^{in} (x) \in \mathrm{SO(3)}$ and $\Lambda (0, x) = \Lambda^{in} (x)$, one easily knows that
	\begin{equation*}
		\begin{aligned}
			f_{ij} (0, x) = \Lambda (0,x) e_{i} \cdot \Lambda (0,x) e_{j} - \delta_{ij} = \Lambda^{in} (x) e_{i} \cdot \Lambda^{in} (x) e_{j} - \delta_{ij} = 0 
		\end{aligned}
	\end{equation*}
	for $1 \leq i, j \leq 3$. Hence $\| f_{ij} (0, \cdot) \|^2_{L^2} = 0$. Then the inequality \eqref{Gronwall-ineq} shows that
	\begin{equation*}
		\begin{aligned}
			\| f_{ij} (t, \cdot) \|^2_{L^2} = 0
		\end{aligned}
	\end{equation*}
	for all $t \geq 0$, which means that $f_{ij} (t,x) \equiv 0$ for $1 \leq i, j, \leq 3$. Consequently, $\Lambda (t,x) \in \mathrm{SO(3)}$ for all $t \geq 0$. The proof of Lemma \ref{Initial} is thereby finished.
\end{proof}

\subsection{The stereographic projection transform for SOHB equations \eqref{SOHB}}

For convenience to study the SOHB system \eqref{SOHB}, we will first perform some proper transforms to deal with the geometric constraint $\Lambda \in \mathrm{SO(3)}$. As shown in Proposition 4.18 of \cite{DFM-2017-MMMAS}, one can further consider the SOHB system \eqref{SOHB} in the terms of the orthonormal basis given by 
$${ \lbrace \Omega = \Lambda e_{1} , \mathbf{u} = \Lambda e_{2} , \mathbf{v} = \Lambda e_{3} \rbrace} \,,$$ 
where $\{ e_1, e_2, e_3 \}$ is the canonical basis of $\R^3$. More precisely, the SOHB system \eqref{SOHB} can be equivalently expressed as
\begin{equation}\label{SOHB-1}
	\left\{
	\begin{aligned}
		\partial_{t} \rho+c_{1} \nabla_{x} \cdot ( \rho \Omega ) & = 0 \,, \\
		\rho D_{t} \Omega + P_{\Omega^{\perp}} ( c_{3} \nabla_{x} \rho + c_{4} \rho \mathbf{r} ) & = 0 \,, \\
		\rho D_{t} \mathbf{u} - \mathbf{u} \cdot ( c_{3} \nabla_{x} \rho + c_{4} \rho \mathbf{r} ) \Omega + c_{4} \rho \delta \mathbf{v} & = 0 \,, \\
		\rho D_{t} \mathbf{v} - \mathbf{v} \cdot ( c_{3} \nabla_{x} \rho + c_{4} \rho \mathbf{r} ) \Omega - c_{4} \rho \delta \mathbf{u} & = 0 \,,
	\end{aligned}
	\right.
\end{equation}
with the nonlinear constraints
\begin{equation}\label{Constraint-SO3}
	\begin{aligned}
		| \Omega | = | \mathbf{u} | = | \mathbf{v} | = 1 \,, \ \Omega \cdot \mathbf{u} = \Omega \cdot \mathbf{v} = \mathbf{u} \cdot \mathbf{v} = 0 \,,
	\end{aligned}
\end{equation}
where $ D_{t} : = \partial_{t} + c_{2} ( \Omega \cdot \nabla_{x} ) $,
\begin{equation*}
	\begin{aligned}
		\delta & = [ ( \Omega \cdot \nabla_{x} ) \mathbf{u} ] \cdot \mathbf{v} + [ ( u \cdot \nabla_{x} ) \mathbf{v} ] \cdot \Omega + [ ( \mathbf{v} \cdot\nabla_{x} ) \Omega ] \cdot \mathbf{u} \,, \\
		\mathbf{r} & = ( \nabla_{x} \cdot \Omega ) \Omega + ( \nabla_{x} \cdot \mathbf{u} ) \mathbf{u} + ( \nabla_{x} \cdot \mathbf{v} ) \mathbf{v} \,.
	\end{aligned}
\end{equation*}
Moreover, the operator $ P_{\Omega^{\perp}} = I - \Omega \otimes \Omega $ denotes the projection on the orthogonal of $\Omega$. 

We further need to consider the geometric constraints \eqref{Constraint-SO3}. Inspired by the work \cite{JLZ-ARMA-2020}, we will we adopt {\em stereographic projection transform}, which, compared with the spherical coordinates transform (see \cite{DLMP-2013-MAA,ZJ-NARWA-2017}, for instance), can avoid the coordinate singularity.

\begin{lemma}\label{Lmm-SOHB-spt}
	Let
\begin{equation}
	\begin{aligned}
		\Omega = ( \tfrac{2 \phi_{1}}{W_{1}} , \tfrac{2 \theta_{1}}{W_{1}} , \tfrac{\phi_{1}^{2} + \theta_{1}^{2} - 1}{W_{1}} )^\top \,,
		\mathbf{u} = ( \tfrac{2 \phi_{2}}{W_{2}} , \tfrac{2 \theta_{2}}{W_{2}} , \tfrac{\phi_{2}^{2} + \theta_{2}^{2} - 1}{W_{2}} )^\top \,,
		\mathbf{v} = ( \tfrac{2 \phi_{3}}{W_{3}} , \tfrac{2 \theta_{3}}{W_{3}} , \tfrac{\phi_{3}^{2} + \theta_{3}^{2} - 1}{W_{3}} )^\top \,,
	\end{aligned}
\end{equation}
where
\begin{equation}
	\begin{aligned}
		W_{1} = 1 + \phi_{1}^{2} + \theta_{1}^{2} \,, W_{2} = 1 + \phi_{2}^{2} + \theta_{2}^{2} \,, W_{3} = 1 + \phi_{3}^{2} + \theta_{3}^{2} ,
	\end{aligned}
\end{equation}
and $( \phi_{1} , \theta_{1} ) , ( \phi_{2} , \theta_{2} ) , ( \phi_{3} , \theta_{3} ) \in \mathbb{R}^{2}$ are the stereographic projection coordinates of $\Omega$, $\mathbf{u}$, $\mathbf{v}$, respectively. Under the constraints $\Omega \cdot \mathbf{u} = \Omega \cdot \mathbf{v} = \mathbf{u} \cdot \mathbf{v} = 0$, the system \eqref{SOHB-1} can be equivalently expressed by
\begin{equation}\label{SHOB-spt}
	\left\{
	\begin{aligned}
		\partial_{t} \rho + c_{1} \Omega \cdot \nabla_{x} \rho + c_{1} \rho \Omega_{\phi_{1}} \cdot \nabla_{x} \phi_{1} + c_{1} \rho \Omega_{\theta_{1}} \cdot \nabla_{x} \theta_{1} = 0 \,, \\
		4 \rho \partial_{t} \phi_{1} + c_{3} W_{1}^{2} \Omega_{\phi_{1}} \cdot \nabla_{x} \rho + 4 c_{2} \rho \Omega \cdot \nabla_{x} \phi_{1} = 0 \,, \\
		4 \rho \partial_{t} \theta_{1} + c_{3} W_{1}^{2} \Omega_{\theta_{1}} \cdot \nabla_{x} \rho + 4 c_{2} \rho \Omega \cdot \nabla_{x} \theta_{1} = 0 \,, \\
		4 \rho \partial_{t} \phi_{2} + 4 c_{2} \rho \Omega \cdot \nabla_{x} \phi_{2} = 0 \,, \\
		4 \rho \partial_{t} \theta_{2} + 4 c_{2} \rho \Omega \cdot \nabla_{x} \theta_{2} = 0 \,, \\
		4 \rho \partial_{t} \phi_{3} + 4 c_{2} \rho \Omega \cdot \nabla_{x} \phi_{3} = 0 \,, \\
		4 \rho \partial_{t} \theta_{3} + 4 c_{2} \rho \Omega \cdot \nabla_{x} \theta_{3} = 0 \,,
	\end{aligned}
	\right.
\end{equation}
where the vector $\Omega_{\phi_{1}}$ and $\Omega_{\theta_{1}}$ are the partial derivatives of $\Omega$ with respect to the variables $\phi_{1}$ and $\theta_{1}$, which can be explicitly represented as
\begin{equation*}
	\begin{aligned}
		\Omega_{\phi_{1}} = ( \tfrac{2 ( 1 - \phi_{1}^{2} + \theta_{1}^{2} ) }{ W_{1}^{2} } , - \tfrac{ 4 \phi_{1} \theta_{1} }{ W_{1}^{2} } , \tfrac{4 \phi_{1} }{ W_{1}^{2} } )^\top \,, \quad
		\Omega_{\theta_{1}} = ( - \tfrac{ 4 \phi_{1} \theta_{1} }{ W_{1}^{2} } , \tfrac{ 2 ( 1 + \phi_{1}^{2} - \theta_{1}^{2} ) }{ W_{1}^{2} } , \tfrac{ 4 \theta_{1} }{ W_{1}^{2} } )^\top \,,
	\end{aligned}
\end{equation*}
respectively. 
\end{lemma}

\begin{proof}
	Note that
	\begin{equation}
		\begin{aligned}
			\Omega_{\phi_{1}} = & ( \tfrac{ 2 ( 1 - \phi_{1}^{2} + \theta_{1}^{2} ) }{ W_{1}^{2} } , - \tfrac{ 4 \phi_{1} \theta_{1} }{ W_{1}^{2} } , \tfrac{ 4 \phi_{1} }{ W_{1}^{2} } )^\top \,, \quad
			\Omega_{\theta_{1}} = & ( - \tfrac{ 4 \phi_{1} \theta_{1} }{ W_{1}^{2} } , \tfrac{ 2 ( 1 + \phi_{1}^{2} - \theta_{1}^{2} ) }{ W_{1}^{2} } , \tfrac{ 4 \theta_{1} }{ W_{1}^{2} } )^\top \,, \\
			\mathbf{u}_{\phi_{2}} = & ( \tfrac{ 2 ( 1 - \phi_{2}^{2} + \theta_{2}^{2} ) }{ W_{2}^{2} } , - \tfrac{ 4 \phi_{2} \theta_{2} }{ W_{2}^{2} } , \tfrac{ 4 \phi_{2} }{ W_{2}^{2} } )^\top \,, \quad
			\mathbf{u}_{\theta_{2}} = & ( - \tfrac{ 4 \phi_{2} \theta_{2} }{ W_{2}^{2} } , \tfrac{ 2 ( 1 + \phi_{2}^{2} - \theta_{2}^{2} ) }{ W_{2}^{2} } , \tfrac{ 4 \theta_{2} }{ W_{2}^{2} } )^\top \,, \\
			\mathbf{v}_{\phi_{3}} = & ( \tfrac{ 2 ( 1 - \phi_{3}^{2} + \theta_{3}^{2} ) }{ W_{3}^{2} } , - \tfrac{ 4 \phi_{3} \theta_{3} }{ W_{3}^{2} } , \tfrac{ 4 \phi_{3} }{ W_{3}^{2} } )^\top \,, \quad
			\mathbf{v}_{\theta_{3}} = & ( - \tfrac{ 4 \phi_{3} \theta_{3} }{ W_{3}^{2} } , \tfrac{ 2 ( 1 + \phi_{3}^{2} - \theta_{3}^{2} ) }{ W_{3}^{2} } , \tfrac{ 4 \theta_{3} }{ W_{3}^{2} } )^\top \,.
		\end{aligned}
	\end{equation}
	It is easy to see that
	\begin{equation}\label{R1}
		\begin{aligned}
			& | \Omega_{\phi_{1}} |^{2} = | \Omega_{\theta_{1}} |^{2} = \tfrac{4}{W_{1}^{2}} \,, \ | \mathbf{u}_{\phi_{2}} |^{2} = | \mathbf{u}_{\theta_{2}} |^{2} = \tfrac{4}{W_{2}^{2}} \,, \ | \mathbf{v}_{\phi_{3}} |^{2} = | \mathbf{v}_{\theta_{3}} |^{2} = \tfrac{4}{W_{3}^{2}} \,, \\
			& \Omega_{\phi_{1}} \cdot \Omega_{\theta_{1}} = \mathbf{u}_{\phi_{2}} \cdot \mathbf{u}_{\theta_{2}} = \mathbf{v}_{\phi_{3}} \cdot \mathbf{v}_{\theta_{3}} = 0 \,.
		\end{aligned}
	\end{equation}
	From the constraints $| \Omega | = | \mathbf{u} | = | \mathbf{v} | = 1$, it follows that 
	\begin{equation}\label{R2}
		\begin{aligned}
			\Omega_{\phi_{1}} \cdot \Omega = \Omega_{\theta_{1}} \cdot \Omega = \mathbf{u}_{\phi_{2}} \cdot \mathbf{u} = \mathbf{u}_{\theta_{2}} \cdot \mathbf{u} = \mathbf{v}_{\phi_{3}} \cdot \mathbf{v} = \mathbf{v}_{\theta_{3}} \cdot \mathbf{v} = 0 \,.
		\end{aligned}
	\end{equation}
	The constraints $\Omega \cdot \mathbf{u} = \Omega \cdot \mathbf{v} = \mathbf{u} \cdot \mathbf{v} = 0$ further indicate that
	\begin{equation}\label{R3}
		\begin{aligned}
			& \Omega_{\phi_{1}} \cdot \mathbf{u} = \Omega_{\theta_{1}} \cdot \mathbf{u} = 0 \,, \qquad \Omega \cdot \mathbf{u}_{\phi_{2}} = \Omega \cdot \mathbf{u}_{\theta_{2}} =0 \,, \\
			& \Omega_{\phi_{1}} \cdot \mathbf{v} = \Omega_{\theta_{1}} \cdot \mathbf{v} = 0 \,, \qquad \Omega \cdot \mathbf{v}_{\phi_{3}} = \Omega \cdot \mathbf{v}_{\theta_{3}} =0 \,, \\
			& \mathbf{u} \cdot \mathbf{v}_{\phi_{3}} = \mathbf{u} \cdot \mathbf{v}_{\theta_{3}} = 0 \,, \qquad \mathbf{v} \cdot \mathbf{u}_{\phi_{2}} = \mathbf{v} \cdot \mathbf{u}_{\theta_{2}} = 0 \,.
		\end{aligned}
	\end{equation}
	Moreover, by (A.1) of Appendix A in \cite{JLZ-ARMA-2020}, one has
	\begin{equation}\label{R4}
		\begin{aligned}
			P_{\Omega^{\perp}} a = \tfrac{1}{4} W_{1}^{2} ( \Omega_{\phi_{1} } \cdot a ) \Omega_{\phi_{1}} + \tfrac{1}{4} W_{1}^{2} ( \Omega_{\theta_{1}} \cdot a ) \Omega_{\theta_{1}}
		\end{aligned}
	\end{equation}
	for all $a\in\mathbb{R}^{3}$. 
	
	Based on the all above relations, we now derive the system \eqref{SHOB-spt}.
	
	{\bf Step 1. Derivation of $\rho$-equation.}
	
	For the second term $ \nabla_{x} \cdot ( \rho \Omega ) $, we derive that
	\begin{equation*}
		\begin{aligned}
			\nabla_{x} \cdot ( \rho \Omega ) & = \Omega \cdot \nabla_{x} \rho + \rho \nabla_{x} \cdot \Omega \\
			& = \Omega \cdot \nabla_{x} \rho + \rho ( \Omega_{\phi_{1}} \cdot \nabla_{x} \phi_{1} + \Omega_{\theta_{1}} \cdot \nabla_{x} \theta_{1} ) \,.
		\end{aligned}
	\end{equation*}
	Then, the first equation of the system \eqref{SHOB-spt} is immediately derived from the first $\rho$-equation of the system \eqref{SOHB-1} and the above equality.
	
	{\bf Step 2. Derivation of $\Omega$-equation.}
	
	For the second $\Omega$-equation, we calculate that by the chain rules of differentiation
	\begin{equation*}
		\begin{aligned}
			& \rho ( \partial_{t} \Omega + c_{2} ( \Omega \cdot \nabla_{x} ) \Omega ) + c_{3} P_{\Omega^{\perp}} \nabla_{x} \rho \\
			= & \rho \Omega_{\phi_{1}} \partial_{t} \phi_{1} + c_{2} \rho ( \Omega \cdot \nabla_{x} \phi_{1} ) \Omega_{\phi_{1}} + \tfrac{ c_{3} W_{1}^{2} }{4} ( \Omega_{\phi_{1}} \cdot \nabla_{x} \rho ) \Omega_{\phi_{1}} + \tfrac{ c_{3} W_{1}^{2} }{4} ( \Omega_{\theta_{1}} \cdot \nabla_{x} \rho ) \Omega_{\theta_{1}} \\
			& + \rho \Omega_{\theta_{1}} \partial_{t} \theta_{1} + c_{2} \rho ( \Omega \cdot \nabla_{x} \theta_{1} ) \Omega_{\theta_{1}} \,,
		\end{aligned}
	\end{equation*}
    and
	\begin{equation*}
		\begin{aligned}
			c_{4} P_{\Omega^{\perp}} \rho \mathbf{r} = & \tfrac{1}{4} W_{1}^{2} ( \Omega_{\phi_{1}} \cdot \rho \mathbf{r} ) \Omega_{\phi_{1}} + \tfrac{1}{4} W_{1}^{2} ( \Omega_{\theta_{1}} \cdot \rho \mathbf{r} ) \Omega_{\theta_{1}} \\
			= & \tfrac{1}{4} W_{1}^{2} \rho \lbrack( \Omega_{\phi_{1}} \cdot \mathbf{u} ) \Omega_{\phi_{1}} + ( \Omega_{\theta_{1}} \cdot \mathbf{u} ) \Omega_{\theta_{1}} \rbrack( \mathbf{u}_{\phi_{2}} \cdot \nabla_{x} \phi_{2} + \mathbf{u}_{\theta_{2}} \cdot \nabla_{x} \theta_{2} ) \\
			& + \tfrac{1}{4} W_{1}^{2} \rho \lbrack( \Omega_{\phi_{1}} \cdot \mathbf{v} ) \Omega_{\phi_{1}} + ( \Omega_{\theta_{1}} \cdot \mathbf{v} ) \Omega_{\theta_{1}} \rbrack( \mathbf{v}_{\phi_{3}} \cdot \nabla_{x} \phi_{3} + \mathbf{v}_{\theta_{3}} \cdot \nabla_{x} \theta_{3} ) \,,
		\end{aligned}
	\end{equation*}
	where the relation \eqref{R4} has been utilized. By the previous two relations, it infers
	\begin{equation}\label{Q2}
		\begin{aligned}
			& \rho ( \partial_{t} \Omega + c_{2} ( \Omega \cdot \nabla_{x} ) \Omega ) + c_{3} P_{\Omega^{\perp}} \nabla_{x} \rho + c_{4} P_{\Omega^{\perp}} \rho \mathbf{r} \\
			= & \rho \Omega_{\phi_{1}} \partial_{t} \phi_{1} + c_{2} \rho ( \Omega \cdot \nabla_{x} \phi_{1} ) \Omega_{\phi_{1}} + \tfrac{ c_{3} W_{1}^{2} }{4} ( \Omega_{\phi_{1}} \cdot \nabla_{x} \rho ) \Omega_{\phi_{1}} \\
			& + \tfrac{ c_{3} W_{1}^{2} }{4} ( \Omega_{\theta_{1}} \cdot \nabla_{x} \rho ) \Omega_{\theta_{1}} + \rho \Omega_{\theta_{1}} \partial_{t} \theta_{1} + c_{2} \rho ( \Omega \cdot \nabla_{x} \theta_{1} ) \Omega_{\theta_{1}} = 0 \,.
		\end{aligned}
	\end{equation}
	Then, by the relation $ \Omega_{\phi_{1}} \cdot \Omega_{\theta_{1}} = 0 $, we derive from the dot product with $ \Omega_{\phi_{1}} $ and $ \Omega_{\theta_{1}} $ in the above equation \eqref{Q2} that the second $\phi_{1}$-equation and the third $\theta_{1}$-equation of the system \eqref{SHOB-spt} hold.
	
	{\bf Step 3. Derivation of $ \mathbf{u} $-equation.}
	
	We need to compute the term $ - c_{4} \rho ( \mathbf{u} \cdot \mathbf{r} ) \Omega $ and the term $ c_{4} \rho \delta \mathbf{v} $. By the chain rules of differentiation, it is easy to see that
	\begin{equation*}
		\begin{aligned}
			-c_{4} \rho ( \mathbf{u} \cdot \mathbf{r} ) \Omega & = - c_{4} \rho \left( \mathbf{u} \cdot ( \nabla_{x} \cdot \Omega ) \Omega + \mathbf{u} \cdot ( \nabla_{x} \cdot \mathbf{u} ) \mathbf{u} + \mathbf{u} \cdot ( \nabla_{x} \cdot \mathbf{v} ) \mathbf{v} \right) \Omega \\
			& = - c_{4} \rho ( \nabla_{x} \cdot \mathbf{u} ) | \mathbf{u} |^{2} \Omega \\
			& = - c_{4} \rho \Omega ( \mathbf{u}_{\phi_{2}} \cdot \nabla_{x} \phi_{2} + \mathbf{u}_{\theta_{2}} \cdot \nabla_{x} \theta_{2} ) \,,
		\end{aligned}
	\end{equation*}
    and
	\begin{equation*}
		\begin{aligned}
			c_{4} \rho \delta \mathbf{v} = & c_{4} \left( \lbrack( \Omega \cdot \nabla_{x} ) \mathbf{u} \rbrack \cdot \mathbf{v} + \lbrack ( \mathbf{u} \cdot \nabla_{x} ) \mathbf{v} \rbrack \cdot \Omega + \lbrack ( \mathbf{v} \cdot \nabla_{x} ) \Omega \rbrack \cdot \mathbf{u} \right) \cdot \mathbf{v} \\
			= & c_{4} \rho \lbrack ( \Omega \cdot \nabla_{x} \phi_{2} ) \mathbf{u}_{\phi_{2}} \cdot \mathbf{v} + ( \Omega \cdot \nabla_{x} \theta_{2} ) \mathbf{u}_{\theta_{2}} \cdot \mathbf{v} + ( \mathbf{u} \cdot \nabla_{x} \phi_{3} ) \mathbf{v}_{\phi_{3}} \cdot \Omega + ( \mathbf{u} \cdot \nabla_{x} \theta_{3} ) \mathbf{v}_{\theta_{3}} \cdot \Omega \rbrack \mathbf{v} \\
			& + c_{4} \rho \lbrack ( \mathbf{v} \cdot \nabla_{x} \phi_{1} ) \Omega_{\phi_{1}} \cdot \mathbf{u} + ( \mathbf{v} \cdot \nabla_{x} \theta_{1} ) \Omega_{\theta_{1}} \cdot \mathbf{u} \rbrack \mathbf{v} \,.
		\end{aligned}
	\end{equation*}
	It therefore implies that
	\begin{equation}
		\begin{aligned}
			&\rho D_{t}\mathbf{u}-\mathbf{u}\cdot (c_{3}\nabla_{x}\rho+c_{4}\rho \mathbf{r})\Omega+c_{4}\rho \delta \mathbf{v}\\
			=&\rho(\mathbf{u}_{\phi_{2}}\partial_{t}\phi_{2}+\mathbf{u}_{\theta_{2}}\partial_{t}\theta_{2})+c_{2}\rho\lbrack(\Omega\cdot\nabla_{x}\phi_{2})\mathbf{u}_{\phi_{2}}+(\Omega\cdot\nabla_{x}\theta_{2})\mathbf{u}_{\theta_{2}}\rbrack-c_{3}\Omega(\mathbf{u}\cdot\nabla_{x}\rho)\\
			&-c_{4}\rho\Omega(\mathbf{u}_{\phi_{2}}\cdot\nabla_{x}\phi_{2}+\mathbf{u}_{\theta_{2}}\cdot\nabla_{x}\theta_{2})=0.
		\end{aligned}
	\end{equation}
	Then, by the relations \eqref{R1}-\eqref{R2}-\eqref{R3} and $ \mathbf{u}_{\phi_{2}} \cdot \mathbf{u}_{\theta_{2}} = 0 $, we derive from the dot product with $\mathbf{u}_{\phi_{2}}$ and 
	$\mathbf{u}_{\theta_{2}}$ in the above equation that the fourth $\phi_{2}$-equation and the fifth $\theta_{2}$-equation of the system \eqref{SHOB-spt} hold.
	
	{\bf Step 4. Derivation of $\mathbf{v}$-equation.}
	By the similar arguments in Step 3, one has
	\begin{equation}
		\begin{aligned}
			& \rho D_{t} \mathbf{v} - \mathbf{v} \cdot ( c_{3} \nabla_{x} \rho + c_{4} \rho \mathbf{r} ) \Omega - c_{4} \rho \delta \mathbf{u} \\
			= & \rho ( \mathbf{v}_{\phi_{3}} \partial_{t} \phi_{3} + \mathbf{v}_{\theta_{3}} \partial_{t} \theta_{3} ) + c_{2} \rho \lbrack ( \Omega \cdot \nabla_{x} \phi_{3} ) \mathbf{v}_{\phi_{3}} + ( \Omega \cdot \nabla_{x} \theta_{3} ) \mathbf{v}_{\theta_{3}} \rbrack - c_{3} \Omega ( \mathbf{v} \cdot \nabla_{x} \rho ) \\
			& - c_{4} \rho \Omega ( \mathbf{v}_{\phi_{3}} \cdot \nabla_{x} \phi_{3} + \mathbf{v}_{\theta_{3}} \cdot \nabla_{x} \theta_{3} ) = 0 \,.
		\end{aligned}
	\end{equation}
	Consequently, based on the relations \eqref{R1}-\eqref{R2}-\eqref{R3}, the sixth $\phi_{3}$-equation and the seventh $\theta_{3}$-equation of the system \eqref{SHOB-spt} can be derived from doting with $\mathbf{v}_{\phi_3}$ and $\mathbf{v}_{\theta_3}$ in the above equation, respectively. Then the proof of Lemma \ref{Lmm-SOHB-spt} is finished.
\end{proof}

\subsection{Well-posedness of the SOHB equations \eqref{SOHB}: Proof of Theorem \ref{LWP}}\label{subsec:Sym}

Denote by
\begin{equation*}
	\begin{aligned}
		\mathnormal{U} = ( \rho , \phi_{1} , \theta_{1} , \phi_{2} , \theta_{2} , \phi_{3} , \theta_{3} )^\top \,.
	\end{aligned}
\end{equation*}
Then the system \eqref{SHOB-spt} can be rewritten as the matrix form
\begin{equation}\label{SHS}
	\begin{aligned}
		\mathbf{A}_{0} ( \mathnormal{U} ) \partial_{t} \mathnormal{U} + \sum_{i=1}^{3} \mathbf{A}^{i} ( \mathnormal{U} ) \partial_{x_{i}} \mathnormal{U} = 0 \,,
	\end{aligned}
\end{equation}
where 
\begin{equation}
	\begin{aligned}
		\mathbf{A}_{0}=
		\begin{pmatrix}
			1 & & & & & & \\
			& 4 \rho & & & & & \\
			& & 4 \rho & & & & \\
			& & & 4 \rho & & & \\
			& & & & 4 \rho & & \\
			& & & & & 4 \rho & \\
			& & & & & & 4 \rho
		\end{pmatrix} \,,
	\end{aligned}
\end{equation}
and
\begin{equation}
	\begin{aligned}
		\mathbf{A}^{i} =
		\begin{pmatrix}
			\mathbf{B}_{3}^{i} & 0 \\
			0 & 4 c_{2} \rho \Omega^{i} \mathbf{I}_{4}
		\end{pmatrix} =
		\begin{pmatrix}
			c_{1} \Omega^{i} & c_{1} \rho \Omega_{\phi_{1}}^{i} & c_{1} \rho \Omega_{\theta_{1}}^{i} & 0 & 0 & 0 & 0 \\
			c_{3} W_{1}^{2} \Omega_{\phi_{1}}^{i} & 4 c_{2} \rho \Omega^{i} & 0 & 0 & 0 & 0 & 0 \\
			c_{3} W_{1}^{2} \Omega_{\theta_{1}}^{i} & 0 & 4 c_{2} \rho \Omega^{i} & 0 & 0 & 0 & 0 \\
			0 & 0 & 0 & 4 c_{2} \rho \Omega^{i} & 0 & 0 & 0 \\
			0 & 0 & 0 & 0 & 4 c_{2} \rho \Omega^{i} & 0 & 0 \\
			0 & 0 & 0 & 0 & 0 & 4 c_{2} \rho \Omega^{i} & 0 \\
			0 & 0 & 0 & 0 & 0 & 0 & 4 c_{2} \rho \Omega^{i}
		\end{pmatrix} 
	\end{aligned}
\end{equation}
for $i = 1,2,3$. Here the symbol $X^i$ stands for the $i$-th component of the vector $X$ for $X = \Omega$, $\Omega_{\phi_{1}}$ and $\Omega_{\theta_{1}}$. Moreover, the symbol $\mathbf{I}_{4}$ represents the $4 \times 4$ identity matrix.

Let
\begin{equation}
	\begin{aligned}
		\mathbf{K}_{3} =
		\begin{pmatrix}
			1 & & \\
			& \tfrac{ c_{1} \rho }{ c_{3} W_{1}^{2} } & \\
			& & \tfrac{ c_{1} \rho }{ c_{3} W_{1}^{2} }
		\end{pmatrix} \,.
	\end{aligned}
\end{equation}
It is easy to verify that
\begin{equation}
	\begin{aligned}
		\widetilde{\mathbf{B}}_{3}^{i} : = \mathbf{K}_{3} \mathbf{B}_{3}^{i} =
		\begin{pmatrix}
			c_{1} \Omega^{i} & c_{1} \rho \Omega_{\phi_{1}}^{i} & c_{1} \rho \Omega_{\theta_{1}}^{i} \\
			c_{1} \rho \Omega_{\phi_{1}}^{i} & \tfrac{ 4 c_{1} c_{2} \rho^{2} }{ c_{3} W_{1}^{2} } \Omega^{i} & 0 \\
			c_{1} \rho \Omega_{\theta_{1}}^{i} & 0 & \tfrac{ 4 c_{1} c_{2} \rho^{2} }{ c_{3} W_{1}^{2} } \Omega^{i}
		\end{pmatrix}
	\end{aligned}
\end{equation}
is a $3 \times 3$ symmetric matrix. Denote by
\begin{equation*}
	\begin{aligned}
		\widetilde{\mathbf{K}}_{3} =
		\begin{pmatrix}
			\mathbf{K}_{3} & \\
			& \mathbf{I}_{4}
		\end{pmatrix} \,.
	\end{aligned}
\end{equation*}
Then the $7 \times 7$ matrices
\begin{equation}
	\begin{aligned}
		\widetilde{\mathbf{A}}^{i} : = \widetilde{ \mathbf{K} }_{3} \mathbf{A}^{i} =
		\begin{pmatrix}
			\widetilde{\mathbf{B}}_{3}^{i}& \\
			&4c_{2}\rho\Omega^{i}\mathbf{I}_{4}
		\end{pmatrix} \quad (i = 1,2,3)
	\end{aligned}
\end{equation}
and
\begin{equation}
	\begin{aligned}
		\widetilde{\mathbf{A}}_{0} : = \widetilde{ \mathbf{K} }_{3} \mathbf{A}_{0} =
		\begin{pmatrix}
			1 & & & & & & \\
			& \tfrac{ 4 c_{1} \rho^{2} }{ c_{3} W_{1}^{2} } & & & & & \\
			& & \tfrac{ 4 c_{1} \rho^{2} }{ c_{3} W_{1}^{2} } & & & & \\
			& & & 4 \rho & & & \\
			& & & & 4 \rho & & \\
			& & & & & 4 \rho & \\
			& & & & & & 4 \rho
		\end{pmatrix}
	\end{aligned}
\end{equation}
are all symmetrical. Consequently, from left multiplying \eqref{SHS} by $\widetilde{\mathbf{K}}_{3}$, it follows the following symmetric hyperbolic system
\begin{equation}\label{U-equ}
	\begin{aligned}
		\widetilde{ \mathbf{A} }_{0} \partial_{t} \mathnormal{U} + \sum_{i=1}^{3} \widetilde{ \mathbf{A} }^{i} \partial_{x_{i}} \mathnormal{U} = 0 \,.
	\end{aligned}
\end{equation}

Then by applying Proposition 2.1 in Page 425 of \cite{MET-2011-AMSS}, one can immediately obtain the following conclusion:

{\em
Let the integer $m \geq 3$. Assume the initial data 
$$U (0, x) = U^{in} (x) : = ( \rho^{in} (x) \,, \phi_{1}^{in} (x) \,, \theta_{1}^{in} (x) \,, \phi_{2}^{in} (x) \,, \theta_{2}^{in} (x) \,, \phi_{3}^{in} (x) \,, \theta_{3}^{in} (x) )^\top $$ 
of the symmetric hyperbolic system \eqref{U-equ} satisfies $U^{in} \in H^m_x$ and $\inf_{x \in \R^3} \rho^{in} (x) > 0$. Then there exists a time $T>0$ such that the system \eqref{U-equ} with initial data $U (0, x) = U^{in} (x)$ admits a unique solution $U (t,x) \in L^{\infty} ( [0,T] \,, H^m_x ) ) \cap H^{1} ( [0,T] \,, H^{m-1}_x ) $ with $\inf_{ (t,x) \in [0, T] \times \R^3} \rho (t,x) > 0$.
}

Therefore, the above conclusion immediately implies the results of Theorem \ref{LWP}.

\section{A priori estimate uniform-in-$\epsilon$: Proof of Theorem \ref{Hydrodynamic}}\label{Sec:UnifEst}

The proof of Theorem \ref{Hydrodynamic} is based on the following key lemma, which represents {\em a priori} estimates for the remainder equation \eqref{RE} uniformly in $\epsilon$: 

\begin{lemma}[Uniform-in-$\epsilon$ Estimate]\label{UEE}
	Assume that $d_\star = d - \frac{25 \sqrt[4]{3} \nu_0 }{ c_1 \lambda_0 } > 0$, where $\lambda_0 > 0$ is the Poincar\'e constant given in Lemma \ref{CE}. Let $T > 0$ be given in Theorem \ref{LWP} and integer $s \geq 2$. For $k = 0, 1, \cdots, s$, define the energy functionals as follows:
	\begin{equation}\label{Ekt}
		\begin{aligned}
			& E_{k} (t) = \sum_{ |\ss| = k } \left\Arrowvert \partial_x^{\ss} ( \tfrac{ f_{R}^{\epsilon} }{ M_{\Lambda_{0}} } - \rho_{R}^{\epsilon} ) \right\Arrowvert_{ L_{x,A}^{2} (  M_{ \Lambda_{0} } ) }^{2} + \sum_{|\ss| = k} \left\Arrowvert \partial_{x}^{\ss} \rho_{R}^{\epsilon} \right\Arrowvert_{L_{x}^{2}}^{2} \,,\\
			& D_{k} (t) = \frac{1}{\epsilon} \sum_{|\ss| = k} \left\Arrowvert \nabla_{A} \partial_{x}^{\ss} \left( \tfrac{ f_{R}^{\epsilon} }{ M_{ \Lambda_{0} } } \right) \right\Arrowvert_{L_{x,A}^{2}( M_{ \Lambda_{0} } ) }^{2} \,.
		\end{aligned}
	\end{equation}
	Then there are some positive constants $C$, $c'_k$, $c''_k$, $\epsilon_{0}$ such that for any $\epsilon\in (0,\epsilon_{0})$ and $t\in [0,T]$,
	\begin{equation}\label{UEE-ineq}
		\begin{aligned}
			\tfrac{\d}{\d t} \mathcal{E} (t) + \mathcal{D} (t) \leq C ( 1 + \mathcal{E} (t) ) + C \epsilon \mathcal{E}^{2} (t) \,,
		\end{aligned}
	\end{equation}
	where
	\begin{equation}
		\begin{aligned}
			\mathcal{E} (t) = \sum_{0 \leq k \leq s} c'_k E_{k} (t) \,, \ \mathcal{D} (t) = \sum_{0 \leq k \leq s} c''_k D_{k} (t) \,.
		\end{aligned}
	\end{equation}
	Here the all existed constants depend on, $s$, the all coefficients in \eqref{SOHB} and \eqref{KM}, the norms $ \| \nabla_x \rho^{in} \|_{H^{s+3}_x} $ and $\| \nabla_x \Lambda_0^{in} \|_{H^{s+3}_x}$.
\end{lemma}

The majority goal of this section is to justify Lemma \ref{UEE}.

\subsection{Preparations}

First, we introduce the coercivity estimate of the linear operator $\L_{M_{\Lambda_0}}$ which provides the dissipative mechanism with parameter singularity $\frac{1}{\epsilon}$ of the remainder equation. 
\begin{lemma}[Coercivity estimates of linear operator $\L_{M_{\Lambda_0}}$]\label{CE}
    Recall that the linear operator $\mathcal{L}_{M_{\Lambda_{0}}}$ reads
    \begin{equation}\label{LMfR}
        \begin{aligned}
            \mathcal{L}_{M_{\Lambda_{0}}} f_{R}^{\epsilon} = d \nabla_{A} \cdot \left( M_{\Lambda_{0}} \nabla_{A} \left( \tfrac{ f_{R}^{\epsilon} }{ M_{ \Lambda_{0} } } \right) \right) \,.
        \end{aligned}
    \end{equation}
Then
\begin{equation}
    \begin{aligned}
        \int_{\mathrm{SO(3)}} -\mathcal{L}_{M_{\Lambda_{0}}}f_{R}^{\epsilon} ( \tfrac{f_{R}^{\epsilon}}{M_{\Lambda_{0}}}- \rho_R^\epsilon ) \d A = d \left\| \nabla_{A} \left( \tfrac{ f_{R}^{\epsilon} }{ M_{ \Lambda_{0} } } \right) \right\|^2_{L^2_A (M_{\Lambda_0})} \geq d \lambda_{0} \left\| \tfrac{ f^{\epsilon}_{R} }{ M_{ \Lambda_{0} } } - \rho_R^\epsilon \right\|^{2}_{ L^2_A ( M_{ \Lambda_{0} } )} 
    \end{aligned}
\end{equation}
with $\lambda_{0} > 0$ being the Poincar{\'e} constant with respect to the Fokker-Planck operator, where 
$$\rho_R^\epsilon = \int_{\mathrm{SO}(3)} f_R^\epsilon \d A \,.$$
\end{lemma}

\begin{proof}
    Observe that by integration by parts over $A \in \mathrm{SO(3)}$,
    \begin{equation}\label{P1}
        \begin{aligned}
            \int_{\mathrm{SO(3)}} -\mathcal{L}_{M_{\Lambda_{0}}}f_{R}^{\epsilon} ( \tfrac{f_{R}^{\epsilon}}{M_{\Lambda_{0}}}- \rho_R^\epsilon ) \d A
            & = d \int_{\mathrm{SO}(3)} M_{ \Lambda_{0} } \nabla_{A} \left( \tfrac{ f_{R}^{\epsilon} }{ M_{ \Lambda_{0} } } \right) \cdot \nabla_{A} \left( \tfrac{ f_{R}^{\epsilon} }{ M_{ \Lambda_{0} } } - \rho_{R} \right) \d A \\
            & = \left\| \nabla_{A} \left( \tfrac{ f_{R}^{\epsilon} }{ M_{ \Lambda_{0} } } \right) \right\|^2_{L^2_A (M_{\Lambda_0})} \,,
        \end{aligned}
    \end{equation}
    where we have used the fact $ \nabla_{A} \rho_{R}^{\epsilon} = 0 $. Due to
    \begin{equation*}
    	\begin{aligned}
    		\int_{\mathrm{SO}(3)} \big( \tfrac{ f_{R}^{\epsilon} }{ M_{ \Lambda_{0} } } - \rho_R^\epsilon \big) M_{\Lambda_0} \d A = 0 \,,
    	\end{aligned}
    \end{equation*}
    the Poincar{\'e} inequality shows that there is a positive constant $\lambda_0 > 0$ such that
    \begin{equation}\label{P2}
    	\begin{aligned}
    		\lambda_{0} \left\| \tfrac{ f^{\epsilon}_{R} }{ M_{ \Lambda_{0} } } - \rho_R^\epsilon \right\|^{2}_{ L^2_A ( M_{ \Lambda_{0} } )} \leq \left\| \nabla_{A} \left( \tfrac{ f_{R}^{\epsilon} }{ M_{ \Lambda_{0} } } \right) \right\|^2_{L^2_A (M_{\Lambda_0})} \,.
    	\end{aligned}
    \end{equation}
    Then the inequalities \eqref{P1} and \eqref{P2} finish the proof of Lemma \ref{CE}.
\end{proof}

Second, we focus on the bounds for expanded term $f_1$. More precisely, we aim at estimating the bounds of the quantities $\| \tfrac{f_1}{ M_{ \Lambda_0 } } \|_{H^s_x L^2_A ( M_{ \Lambda_0 } )}$ and $ \| \tfrac{ R ( f_1 ) }{ M_{ \Lambda_0 } } \|_{H^s_x L^2_A ( M_{ \Lambda_0 } )} $, which will utilized later.

\begin{lemma}[Bounds for expanded term $f_1$]\label{Lmm-f1}
	Let $s \geq 2$ and the function $f_1 (t,x, A)$ be determined in \eqref{f1}. Let $(\rho_0^{in}, \Lambda_0^{in})$ be the initial data of $(\rho_0, \Lambda_0)$ which subjects to the SOHB system \eqref{SOHB}. Assume that $(\rho_0^{in}, \Lambda_0^{in})$ obeys the hypotheses given in Theorem \ref{LWP} with $m = s + 4$. Then we have
	\begin{equation}\label{f1-Rf1-bnd}
		\begin{aligned}
			& \| \tfrac{f_1}{ M_{ \Lambda_0 } } \|_{ H^s_x L^2_A ( M_{ \Lambda_0 } ) }^2 + \| \nabla_A ( \tfrac{f_1}{ M_{ \Lambda_0 } } ) \|^2_{H^s_x L^2_A ( M_{ \Lambda_0 } ) } \leq C ( \| \nabla_x \rho_0^{in} \|_{H^{s+2}_x}, \| \nabla_x \Lambda_0^{in} \|_{H^{s+2}_x} ) \,, \\
			& \| \tfrac{ R ( f_1 ) }{ M_{ \Lambda_0 } } \|_{H^s_x L^2_A ( M_{ \Lambda_0 } )} \leq C ( \| \nabla_x \rho_0^{in} \|_{H^{s+3}_x}, \| \nabla_x \Lambda_0^{in} \|_{H^{s+3}_x} ) \,.
		\end{aligned}
	\end{equation}
\end{lemma}

The proof will be given in Section \ref{Sec:f1} later.

\subsection{Uniform-in-$\epsilon$ estimates: Proof of Lemma \ref{UEE}}

The next items are devoted to the proof of Lemma \ref{UEE}. As shown in Lemma \ref{CE}, the dissipative structure comes from the operator $\L_{M_{ \Lambda_0 }} f_R^\epsilon$ by multiplying the unknowns $ \tfrac{ f_{R}^{\epsilon} }{ M_{ \Lambda_{0} } } - \rho_R^\epsilon$. In the remainder equation \eqref{RE}, the operators $\L_{ M_{ \Lambda_0 } } f_R^\epsilon $, $ L_R f_R^\epsilon $ and $\widetilde{Q} (f_R^\epsilon)$ are all divergence form with respect to the variable $A \in \mathrm{SO(3)}$, which means that
\begin{equation*}
	\begin{aligned}
		\int_{\mathrm{SO}(3)} \big( \tfrac{1}{\epsilon} \L_{ M_{ \Lambda_0 } } f_R^\epsilon - \tfrac{1}{\epsilon} L_R f_R^\epsilon + \widetilde{Q} (f_R^\epsilon) \big) \d A = 0 \,.
	\end{aligned}
\end{equation*}
In other words, $\tfrac{1}{\epsilon} \L_{ M_{ \Lambda_0 } } f_R^\epsilon - \tfrac{1}{\epsilon} L_R f_R^\epsilon + \widetilde{Q} (f_R^\epsilon) \in L^2_A (M_{ \Lambda_0 }) / \mathrm{Ker} (\L_{M_{ \Lambda_0 }})$. It then follows that from projecting the remainder equation \eqref{RE} into $\mathrm{Ker} (\L_{M_{ \Lambda_0 }})$ that
\begin{equation}\label{macro-eq}
	\begin{aligned}
		\partial_{t} \rho_{R}^{\epsilon} + \int_{\mathrm{SO(3)}} A e_{1} \cdot \nabla_{x} f_{R}^{\epsilon} \d A = \int_{\mathrm{SO}(3)} R (f_1) \d A \,.
	\end{aligned}
\end{equation}
We call \eqref{macro-eq} the macro-equation of the remainder equation \eqref{RE}. Note that $f_1 = \P_{\L} f_1$ as given in \eqref{f1}. It thereby infers that
\begin{equation}\label{R-f1-int}
	\begin{aligned}
		\int_{\mathrm{SO}(3)} R (f_1) \d A = \nabla_x \cdot \int_{\mathrm{SO}(3)} A e_1 f_1 \d A \,.
	\end{aligned}
\end{equation}

Note that 
\begin{equation*}
	\begin{aligned}
		\lambda [f_R^\epsilon - \rho_R^\epsilon M_{\Lambda_0}] = \lambda [f_R^\epsilon] - \rho_R^\epsilon \lambda [ M_{\Lambda_0} ] = \lambda [f_R^\epsilon] - c_1 \rho_R^\epsilon \Lambda_0 \,,
	\end{aligned}
\end{equation*}
where the last equality is derived from Lemma \ref{CRFM}. Then one has
\begin{equation*}
	\begin{aligned}
		P_{T_{\Lambda_0}} ( \lambda [f_R^\epsilon - \rho_R^\epsilon M_{\Lambda_0}] ) = P_{T_{\Lambda_0}} ( \lambda [f_R^\epsilon] ) \,.
	\end{aligned}
\end{equation*}
Recalling the definition of $L_R f_R^\epsilon$ in \eqref{LR-opt}, one sees
\begin{equation*}
	\begin{aligned}
		L_R f_R^\epsilon = L_R (f_R^\epsilon - \rho_R^\epsilon M_{\Lambda_0}) \,.
	\end{aligned}
\end{equation*}
Then the micro form of the remainder equation \eqref{RE} can be expressed by
\begin{equation}\label{micro-eq}
	\begin{aligned}
		(\partial_t + A e_1 \cdot \nabla_x ) (f_R^\epsilon - \rho_R^\epsilon M_{ \Lambda_0 }) + (\partial_t + A e_1 \cdot \nabla_x ) (\rho_R^\epsilon M_{ \Lambda_0 }) - \tfrac{1}{\epsilon} \L_{ M_{ \Lambda_0 } } f_R^\epsilon \\
		= - \tfrac{1}{\epsilon} L_R (f_R^\epsilon - \rho_R^\epsilon M_{\Lambda_0}) + \widetilde{Q} (f_R^\epsilon) + R (f_1) \,.
	\end{aligned}
\end{equation}

\subsubsection{Basic $L^2$-estimates}

In this subsection, we mainly derive the $L^2$-estimates of the remainder equation \eqref{RE}. The basic energy and energy dissipation structures in $L^2$-estimates is important for higher order derivatives estimates. Moreover, some key cancellations in basic estimates will also instruct how to close higher order derivatives bounds. This idea has been clearly illustrated in \cite{JL-SIMA-2019}.

Via multiplying the micro equation \eqref{micro-eq} by $\tfrac{ f_{R}^{\epsilon} }{ M_{\Lambda_{0}}} - \rho_{R}^\epsilon$ and taking integration over $(x, A) \in \R^3 \times \mathrm{SO(3)}$, it implies that
	\begin{align}\label{L2-1}
		\no & \iint_{\R^3 \times \mathrm{SO(3)}} [ (\partial_t + A e_1 \cdot \nabla_x ) (f_R^\epsilon - \rho_R^\epsilon M_{ \Lambda_0 }) ] (\tfrac{ f_{R}^{\epsilon} }{ M_{\Lambda_{0}}} - \rho_{R}^\epsilon) \d A \d x \\
		\no & + \iint_{\R^3 \times \mathrm{SO(3)}} [ (\partial_t + A e_1 \cdot \nabla_x ) (\rho_R^\epsilon M_{ \Lambda_0 }) ] (\tfrac{ f_{R}^{\epsilon} }{ M_{\Lambda_{0}}} - \rho_{R}^\epsilon) \d A \d x \\
		\no & + \iint_{\R^3 \times \mathrm{SO(3)}} - \frac{1}{\epsilon} \L_{ M_{ \Lambda_0 } } f_R^\epsilon (\tfrac{ f_{R}^{\epsilon} }{ M_{\Lambda_{0}}} - \rho_{R}^\epsilon) \d A \d x \\
		\no & + \iint_{\R^3 \times \mathrm{SO(3)}} \frac{1}{\epsilon} (\tfrac{ f_{R}^{\epsilon} }{ M_{\Lambda_{0}}} - \rho_{R}^\epsilon) L_R ( f_{R}^{\epsilon} - \rho_{R}^\epsilon M_{\Lambda_{0}} ) \d A \d x \\
		= & \iint_{\R^3 \times \mathrm{SO(3)}} \widetilde{Q} ( f_R^\epsilon ) (\tfrac{ f_{R}^{\epsilon} }{ M_{\Lambda_{0}}} - \rho_{R}^\epsilon) \d A \d x + \iint_{\R^3 \times \mathrm{SO(3)}} R (f_1) (\tfrac{ f_{R}^{\epsilon} }{ M_{\Lambda_{0}}} - \rho_{R}^\epsilon) \d A \d x \,.
	\end{align}

{\bf Step 1. Estimates of transport effect for micro part: $(\partial_t + A e_1 \cdot \nabla_x ) (f_R^\epsilon - \rho_R^\epsilon M_{ \Lambda_0 })$.} It is easy to see that
\begin{equation*}
	\begin{aligned}
		& \iint_{\R^3 \times \mathrm{SO(3)}} [ (\partial_t + A e_1 \cdot \nabla_x ) (f_R^\epsilon - \rho_R^\epsilon M_{ \Lambda_0 }) ] (\tfrac{ f_{R}^{\epsilon} }{ M_{\Lambda_{0}}} - \rho_{R}^\epsilon) \d A \d x \\
		= & \iint_{\R^3 \times \mathrm{SO(3)}} [ (\partial_t + A e_1 \cdot \nabla_x ) (\tfrac{f_R^\epsilon}{ M_{ \Lambda_0 } } - \rho_R^\epsilon ) ] (\tfrac{f_R^\epsilon}{ M_{ \Lambda_0 } } - \rho_R^\epsilon ) M_{\Lambda_0} \d A \d x \\
		& + \iint_{\R^3 \times \mathrm{SO(3)}} (\tfrac{f_R^\epsilon}{ M_{ \Lambda_0 } } - \rho_R^\epsilon )^2 (\partial_t + A e_1 \cdot \nabla_x ) M_{ \Lambda_0 } \d A \d x \\
		= & \tfrac{1}{2} \tfrac{\d}{\d t} \big\| \tfrac{f_R^\epsilon}{ M_{ \Lambda_0 } } - \rho_R^\epsilon \big\|^2_{L^2_{x,A} (M_{\Lambda_0})} + \iint_{\R^3 \times \mathrm{SO(3)}} (\tfrac{f_R^\epsilon}{ M_{ \Lambda_0 } } - \rho_R^\epsilon )^2 (\partial_t + \tfrac{1}{2} A e_1 \cdot \nabla_x ) M_{ \Lambda_0 } \d A \d x \,.
	\end{aligned}
\end{equation*}
Recalling the expression of $M_{\Lambda_0}$ in \eqref{M-Lambda-0}, there holds
\begin{equation}\label{M-der}
	\begin{aligned}
		\partial M_{\Lambda_0} = \tfrac{\nu_0}{d} (A \cdot \partial \Lambda_0) M_{\Lambda_0} \,,
	\end{aligned}
\end{equation}
where the symbol $\partial$ stands for $\partial_t $ or $\nabla_x$. Observe that $|A| = \sqrt{A \cdot A} = \sqrt{3} $ and $|A e_1| = 1$ for $A \in \mathrm{SO(3)}$. Then
\begin{equation}\label{M0-der}
	\begin{aligned}
		| (\partial_t + \tfrac{1}{2} A e_1 \cdot \nabla_x ) M_{\Lambda_0} | \leq \tfrac{\sqrt{3} \nu_0}{d} ( |\partial_t \Lambda_0 | + |\nabla_x \Lambda_0| ) M_{\Lambda_0} : = \tfrac{\sqrt{3} \nu_0}{d} |\partial_{t,x} \Lambda_0 | M_{\Lambda_0} \,,
	\end{aligned}
\end{equation}
which means that
\begin{equation}\label{Bnd-L2-1}
	\begin{aligned}
		& | \iint_{\R^3 \times \mathrm{SO(3)}} (\tfrac{f_R^\epsilon}{ M_{ \Lambda_0 } } - \rho_R^\epsilon )^2 (\partial_t + \tfrac{1}{2} A e_1 \cdot \nabla_x ) M_{ \Lambda_0 } \d A \d x | \\
		\leq & \tfrac{\sqrt{3} \nu_0}{d} \| \partial_{t,x} \Lambda_0 \|_{L^\infty_x} \big\| \tfrac{f_R^\epsilon}{ M_{ \Lambda_0 } } - \rho_R^\epsilon \big\|^2_{L^2_{x,A} (M_{\Lambda_0})} \,.
	\end{aligned}
\end{equation}
It therefore holds that
	\begin{align}\label{MiTE-bnd}
		\no & \iint_{\R^3 \times \mathrm{SO(3)}} [ (\partial_t + A e_1 \cdot \nabla_x ) (f_R^\epsilon - \rho_R^\epsilon M_{ \Lambda_0 }) ] (\tfrac{ f_{R}^{\epsilon} }{ M_{\Lambda_{0}}} - \rho_{R}^\epsilon) \d A \d x \\
		\no \geq & \tfrac{1}{2} \tfrac{\d}{\d t} \big\| \tfrac{f_R^\epsilon}{ M_{ \Lambda_0 } } - \rho_R^\epsilon \big\|^2_{L^2_{x,A} (M_{\Lambda_0})} - \tfrac{\sqrt{3} \nu_0}{d} \| \partial_{t,x} \Lambda_0 \|_{L^\infty_x} \big\| \tfrac{f_R^\epsilon}{ M_{ \Lambda_0 } } - \rho_R^\epsilon \big\|^2_{L^2_{x,A} (M_{\Lambda_0})} \\
		\geq & \tfrac{1}{2} \tfrac{\d}{\d t} \big\| \tfrac{f_R^\epsilon}{ M_{ \Lambda_0 } } - \rho_R^\epsilon \big\|^2_{L^2_{x,A} (M_{\Lambda_0})} - \tfrac{\sqrt{3} \nu_0}{d} \| \partial_{t,x} \Lambda_0 \|_{H^2_x} E_0 (t) \,,
	\end{align}
where the functional $E_0 (t)$ can be found in \eqref{Ekt}.

{\bf Step 2. Estimates of transport effect for macro part: $(\partial_t + A e_1 \cdot \nabla_x ) (\rho_R^\epsilon M_{ \Lambda_0 })$.} A direct computation shows
\begin{equation}\label{MaTE-bnd-0}
	\begin{aligned}
		& \iint_{\R^3 \times \mathrm{SO(3)}} [ (\partial_t + A e_1 \cdot \nabla_x ) (\rho_R^\epsilon M_{ \Lambda_0 }) ] (\tfrac{ f_{R}^{\epsilon} }{ M_{\Lambda_{0}}} - \rho_{R}^\epsilon) \d A \d x \\
		= & \iint_{\R^3 \times \mathrm{SO(3)}} (\partial_t + A e_1 \cdot \nabla_x ) \rho_R^\epsilon \cdot (f_R^\epsilon - \rho_R^\epsilon M_{\Lambda_0}) \d A \d x \\
		& + \iint_{\R^3 \times \mathrm{SO(3)}} \rho_R^\epsilon (\partial_t + A e_1 \cdot \nabla_x ) M_{\Lambda_0} \cdot (\tfrac{ f_{R}^{\epsilon} }{ M_{\Lambda_{0}}} - \rho_{R}^\epsilon) \d A \d x \,.
	\end{aligned}
\end{equation}
Note that the cancellation
\begin{equation*}
	\begin{aligned}
		\int_{\mathrm{SO}(3)} (f_R^\epsilon - \rho_R^\epsilon M_{\Lambda_0}) \d A = 0 \,.
	\end{aligned}
\end{equation*}
Then
\begin{equation}\label{MaTE-bnd-1}
	\begin{aligned}
		\iint_{\R^3 \times \mathrm{SO(3)}} \partial_t \rho_R^\epsilon \cdot (f_R^\epsilon - \rho_R^\epsilon M_{\Lambda_0}) \d A \d x = 0 \,.
	\end{aligned}
\end{equation}
Moreover, the integration by parts over $x \in \R^3$ indicates that
\begin{equation*}
	\begin{aligned}
		& \iint_{\R^3 \times \mathrm{SO(3)}} A e_1 \cdot \nabla_x \rho_R^\epsilon ( f_R^\epsilon - \rho_R^\epsilon M_{\Lambda_0} ) \d A \d x \\
		= & - \iint_{\R^3 \times \mathrm{SO(3)}} \rho_R^\epsilon A e_1 \cdot \nabla_x f_R^\epsilon \d A \d x + \tfrac{1}{2} \iint_{\R^3 \times \mathrm{SO(3)}} (\rho_R^\epsilon)^2 A e_1 \cdot \nabla_x M_{\Lambda_0} \d A \d x \,.
	\end{aligned}
\end{equation*}
It then follows from the macro-equation \eqref{macro-eq} that
\begin{equation*}
	\begin{aligned}
		- \iint_{\R^3 \times \mathrm{SO(3)}} \rho_R^\epsilon A e_1 \cdot \nabla_x f_R^\epsilon \d A \d x = & \int_{\R^3} \rho_R^\epsilon \Big( \partial_t \rho_R^\epsilon - \int_{\mathrm{SO}(3)} R (f_1) \d A \Big) \d x \\
		= & \tfrac{1}{2} \tfrac{\d}{\d t} \| \rho_R^\epsilon \|^2_{L^2_x} - \iint_{\R^3 \times \mathrm{SO(3)}} \rho_R^\epsilon R (f_1) \d A \d x \,.
	\end{aligned}
\end{equation*}
Hence,
\begin{equation*}
	\begin{aligned}
		& \iint_{\R^3 \times \mathrm{SO(3)}} A e_1 \cdot \nabla_x \rho_R^\epsilon ( f_R^\epsilon - \rho_R^\epsilon M_{\Lambda_0} ) \d A \d x \\
		= & \tfrac{1}{2} \tfrac{\d}{\d t} \| \rho_R^\epsilon \|^2_{L^2_x} - \iint_{\R^3 \times \mathrm{SO(3)}} \rho_R^\epsilon R (f_1) \d A \d x + \tfrac{1}{2} \iint_{\R^3 \times \mathrm{SO(3)}} (\rho_R^\epsilon)^2 A e_1 \cdot \nabla_x M_{\Lambda_0} \d A \d x \,.
	\end{aligned}
\end{equation*}
The H\"older inequality implies
\begin{equation*}
	\begin{aligned}
		| \iint_{\R^3 \times \mathrm{SO(3)}} \rho_R^\epsilon R (f_1) \d A \d x | \leq \| \int_{\mathrm{SO}(3)} R (f_1) \d A \|_{L^2_x} \| \rho_R^\epsilon \|_{L^2_x} \,.
	\end{aligned}
\end{equation*}
Recalling \eqref{R-f1-int}, one gains
\begin{equation*}
	\begin{aligned}
		| \int_{\mathrm{SO}(3)} R (f_1) \d A | = | \int_{\mathrm{SO}(3)} A e_1 \cdot \tfrac{\nabla_x f_1}{ M_{ \Lambda_0 } } M_{ \Lambda_0 } \d A | \leq \| \tfrac{\nabla_x f_1}{ M_{ \Lambda_0 } } \|_{L^2_A ( M_{ \Lambda_0 } )} \,,
	\end{aligned}
\end{equation*}
which means that
\begin{equation*}
	\begin{aligned}
		\| \int_{\mathrm{SO}(3)} R (f_1) \d A \|_{L^2_x} \leq \| \tfrac{\nabla_x f_1}{ M_{ \Lambda_0 } } \|_{L^2_{x,A} ( M_{ \Lambda_0 } )} \,.
	\end{aligned}
\end{equation*}
Then
\begin{equation*}
	\begin{aligned}
		| \iint_{\R^3 \times \mathrm{SO(3)}} \rho_R^\epsilon R (f_1) \d A \d x | \leq \| \tfrac{\nabla_x f_1}{ M_{ \Lambda_0 } } \|_{L^2_{x,A} ( M_{ \Lambda_0 } )} \| \rho_R^\epsilon \|_{L^2_x} \,.
	\end{aligned}
\end{equation*}
It further infers from the similar arguments in \eqref{Bnd-L2-1} that
\begin{equation*}
	\begin{aligned}
		\tfrac{1}{2} | \iint_{\R^3 \times \mathrm{SO(3)}} (\rho_R^\epsilon)^2 A e_1 \cdot \nabla_x M_{\Lambda_0} \d A \d x | \leq \tfrac{\sqrt{3} \nu_0}{2 d} \| \partial_{t,x} \Lambda_0 \|_{L^\infty_x} \| \rho_R^\epsilon \|^2_{L^2_x} \,.
	\end{aligned}
\end{equation*}
It thereby holds
\begin{equation}\label{MaTE-bnd-2}
	\begin{aligned}
		& \iint_{\R^3 \times \mathrm{SO(3)}} A e_1 \cdot \nabla_x \rho_R^\epsilon ( f_R^\epsilon - \rho_R^\epsilon M_{\Lambda_0} ) \d A \d x \\
		\geq & \tfrac{1}{2} \tfrac{\d}{\d t} \| \rho_R^\epsilon \|^2_{L^2_x} - \tfrac{\sqrt{3} \nu_0}{2 d} \| \partial_{t,x} \Lambda_0 \|_{L^\infty_x} \| \rho_R^\epsilon \|^2_{L^2_x} - \| \tfrac{\nabla_x f_1}{ M_{ \Lambda_0 } } \|_{L^2_{x,A} ( M_{ \Lambda_0 } )} \| \rho_R^\epsilon \|_{L^2_x} \,.
	\end{aligned}
\end{equation}
Moreover, the bound \eqref{M0-der} reduces to
\begin{equation}\label{MaTE-bnd-3}
	\begin{aligned}
		& | \iint_{\R^3 \times \mathrm{SO(3)}} \rho_R^\epsilon (\partial_t + A e_1 \cdot \nabla_x ) M_{\Lambda_0} \cdot (\tfrac{ f_{R}^{\epsilon} }{ M_{\Lambda_{0}}} - \rho_{R}^\epsilon) \d A \d x | \\ 
		\leq & \tfrac{\sqrt{3} \nu_0}{d} \| \partial_{t,x} \Lambda_0 \|_{L^\infty_x} \| \rho_R^\epsilon \|_{L^2_x} \big\| \tfrac{f_R^\epsilon}{ M_{ \Lambda_0 } } - \rho_R^\epsilon \big\|_{L^2_{x,A} (M_{\Lambda_0}) } \,.
	\end{aligned}
\end{equation}
Plugging \eqref{MaTE-bnd-1}, \eqref{MaTE-bnd-2} and \eqref{MaTE-bnd-3} into \eqref{MaTE-bnd-0}, one gains
\begin{equation}\label{MaTE-bnd}
	\begin{aligned}
		& \iint_{\R^3 \times \mathrm{SO(3)}} [ (\partial_t + A e_1 \cdot \nabla_x ) (\rho_R^\epsilon M_{ \Lambda_0 }) ] (\tfrac{ f_{R}^{\epsilon} }{ M_{\Lambda_{0}}} - \rho_{R}^\epsilon) \d A \d x \\
		\geq & \tfrac{1}{2} \tfrac{\d}{\d t} \| \rho_R^\epsilon \|^2_{L^2_x} - \tfrac{\sqrt{3} \nu_0}{2 d} \| \partial_{t,x} \Lambda_0 \|_{L^\infty_x} \| \rho_R^\epsilon \|^2_{L^2_x} - \| \tfrac{\nabla_x f_1}{ M_{ \Lambda_0 } } \|_{L^2_{x,A} ( M_{ \Lambda_0 } )} \| \rho_R^\epsilon \|_{L^2_x} \\
		& - \tfrac{\sqrt{3} \nu_0}{d} \| \partial_{t,x} \Lambda_0 \|_{L^\infty_x} \| \rho_R^\epsilon \|_{L^2_x} \big\| \tfrac{f_R^\epsilon}{ M_{ \Lambda_0 } } - \rho_R^\epsilon \big\|_{L^2_{x,A} (M_{\Lambda_0}) } \\
		\geq & \tfrac{1}{2} \tfrac{\d}{\d t} \| \rho_R^\epsilon \|^2_{L^2_x} - \tfrac{3 \sqrt{3} \nu_0}{2 d} \| \partial_{t,x} \Lambda_0 \|_{H^2_x} E_0 (t) - \| \tfrac{\nabla_x f_1}{ M_{ \Lambda_0 } } \|_{L^2_{x,A} ( M_{ \Lambda_0 } )} E_0^\frac{1}{2} (t) \,,
	\end{aligned}
\end{equation}
where the functional $E_0 (t)$ is defined in \eqref{Ekt}.

{\bf Step 3. Estimates of coercivity: $ - \frac{1}{\epsilon} \L_{ M_{ \Lambda_0 } } f_R^\epsilon $.} By Lemma \ref{CE}, one directly obtains
\begin{equation}\label{Coer-Bnd}
	\begin{aligned}
		\iint_{\R^3 \times \mathrm{SO(3)}} - \frac{1}{\epsilon} \L_{ M_{ \Lambda_0 } } f_R^\epsilon (\tfrac{ f_{R}^{\epsilon} }{ M_{\Lambda_{0}}} - \rho_{R}^\epsilon) \d A \d x = \frac{d}{\epsilon} \left\| \nabla_{A} \left( \tfrac{ f_{R}^{\epsilon} }{ M_{ \Lambda_{0} } } \right) \right\|^2_{L^2_{x,A} (M_{\Lambda_0})} = d D_0 (t) \,,
	\end{aligned}
\end{equation}
where the functional $D_0 (t)$ is mentioned in \eqref{Ekt}.

{\bf Step 4. Estimates of error linear operator: $ \frac{1}{\epsilon} L_R ( f_{R}^{\epsilon} - \rho_{R}^\epsilon M_{\Lambda_{0}} ) $.} Recalling the definition of $L_R ( f_R^\epsilon - \rho_R^\epsilon M_{\Lambda_0} )$ in \eqref{LR-opt}, one has
\begin{equation}\label{EE-0}
	\begin{aligned}
		& - \iint_{\R^3 \times \mathrm{SO(3)}} \frac{1}{\epsilon} (\tfrac{ f_{R}^{\epsilon} }{ M_{\Lambda_{0}}} - \rho_{R}^\epsilon) L_R ( f_{R}^{\epsilon} - \rho_{R}^\epsilon M_{\Lambda_{0}} ) \d A \d x \\
		= & - \frac{1}{\epsilon} \iint_{\R^3 \times \mathrm{SO(3)}} (\tfrac{ f_{R}^{\epsilon} }{ M_{\Lambda_{0}}} - \rho_{R}^\epsilon) \tfrac{\nu_0}{c_1} \nabla_A \cdot \big[ M_{\Lambda_0} \nabla_A \big( A \cdot P_{T_{\Lambda_0}} ( \lambda [ f_R^\epsilon - \rho_R^\epsilon M_{\Lambda_0} ] ) \big) \big] \d A \d x \\
		= & \tfrac{\nu_0}{c_1} \tfrac{1}{\epsilon} \iint_{\R^3 \times \mathrm{SO(3)}} \nabla_A \left( \tfrac{f_R^\epsilon}{M_{\Lambda_0}} \right) \cdot \nabla_A \big( A \cdot P_{T_{\Lambda_0}} ( \lambda [ f_R^\epsilon - \rho_R^\epsilon M_{\Lambda_0} ] ) \big) M_{\Lambda_0} \d A \d x \\
		\leq & \tfrac{\nu_0}{c_1} \tfrac{1}{\epsilon} \iint_{\R^3 \times \mathrm{SO(3)}} \left| \nabla_A \left( \tfrac{f_R^\epsilon}{M_{\Lambda_0}} \right) \right| \big| \nabla_A \big( A \cdot P_{T_{\Lambda_0}} ( \lambda [ f_R^\epsilon - \rho_R^\epsilon M_{\Lambda_0} ] ) \big) \big| M_{\Lambda_0} \d A \d x \,.
	\end{aligned}
\end{equation}
By Lemma \ref{Lmm-Proj-SO3},
\begin{equation*}
	\begin{aligned}
		\nabla_A \big( A \cdot P_{T_{\Lambda_0}} ( \lambda [ f_R^\epsilon - \rho_R^\epsilon M_{\Lambda_0} ] ) \big) = \tfrac{1}{2} \big\{ P_{T_{\Lambda_0}} ( \lambda [ f_R^\epsilon - \rho_R^\epsilon M_{\Lambda_0} ] ) - A P_{T_{\Lambda_0}} ( \lambda [ f_R^\epsilon - \rho_R^\epsilon M_{\Lambda_0} ] )^\top A \big\} \,.
	\end{aligned}
\end{equation*}
Note that for any $B \in \R^{3 \times 3}$ and $A \in \mathrm{SO(3)}$,
	\begin{align}\label{EE-1}
		\no |A B^\top A| = & \sqrt{ (A B^\top A) \cdot (A B^\top A) } = \sqrt{ \sum_{i,j=1}^3 (A B^\top A)_{ij} } = \sqrt{ \sum_{i,j=1}^3 \big( \sum_{k,l=1}^3 A_{ik} B_{lk} A_{lj} \big)^2 } \\
		\leq & 3 \sqrt{ \sum_{i,j=1}^3 \sum_{k,l=1}^3 A_{ik}^2 B_{lk}^2 A_{lj}^2 } \leq 3 \sqrt{\big( \sum_{i,j=1}^3 A_{ij}^2 \big)^2 \sum_{k,l=1}^3 B_{kl}^2 } = 9 \sqrt{\sum_{k,l=1}^3 B_{kl}^2 } = 9 |B| \,,
	\end{align}
where we have used the fact that $\sum_{i,j=1}^3 A_{ij}^2 = |A|^2 = 3$ for $A \in \mathrm{SO(3)}$. As a result,
\begin{equation}\label{EE-2}
	\begin{aligned}
		| \nabla_A \big( A \cdot P_{T_{\Lambda_0}} ( \lambda [ f_R^\epsilon - \rho_R^\epsilon M_{\Lambda_0} ] ) \big) | \leq 5 | P_{T_{\Lambda_0}} ( \lambda [ f_R^\epsilon - \rho_R^\epsilon M_{\Lambda_0} ] ) | \,.
	\end{aligned}
\end{equation}
Moreover, by Lemma \ref{Lmm-Proj-SO3},
\begin{equation*}
	\begin{aligned}
		P_{T_{\Lambda_0}} ( \lambda [ f_R^\epsilon - \rho_R^\epsilon M_{\Lambda_0} ] ) = \tfrac{1}{2} \big\{ \lambda [ f_R^\epsilon - \rho_R^\epsilon M_{\Lambda_0} ] - \Lambda_0 \lambda [ f_R^\epsilon - \rho_R^\epsilon M_{\Lambda_0} ]^\top \Lambda_0 \big\} \,.
	\end{aligned}
\end{equation*}
Due to $\Lambda_0 \in \mathrm{SO(3)}$, \eqref{EE-1} implies that
\begin{equation}\label{EE-3}
	\begin{aligned}
		| P_{T_{\Lambda_0}} ( \lambda [ f_R^\epsilon - \rho_R^\epsilon M_{\Lambda_0} ] ) | \leq 5 | \lambda [ f_R^\epsilon - \rho_R^\epsilon M_{\Lambda_0} ] | \,.
	\end{aligned}
\end{equation}
Recalling the definition of $\lambda [f]$ in \eqref{KM}, one has
\begin{equation}\label{EE-4}
	\begin{aligned}
		| \lambda [ f_R^\epsilon - \rho_R^\epsilon M_{\Lambda_0} ] | \leq & \int_{\mathrm{SO}(3)} |A| \big| \tfrac{f_R^\epsilon}{M_{\Lambda_0}} - \rho_R^\epsilon \big| M_{\Lambda_0} \d A \\
		\leq & \Big( \int_{\mathrm{SO}(3)} |A| M_{\Lambda_0} \d A \Big)^\frac{1}{2} \Big( \int_{\mathrm{SO}(3)} \big| \tfrac{f_R^\epsilon}{M_{\Lambda_0}} - \rho_R^\epsilon \big|^2 M_{\Lambda_0} \d A \Big)^\frac{1}{2} \\
		\leq & \tfrac{ \sqrt[4]{3} }{ \sqrt{\lambda_0} } \| \nabla_A ( \tfrac{f_R^\epsilon}{M_{\Lambda_0}} ) \|_{L^2_A (M_{\Lambda_0})} \,,
	\end{aligned}
\end{equation}
where the last inequality is derived from the facts $|A| = \sqrt{3}$, $\int_{\mathrm{SO}(3)} M_{\Lambda_0} \d A = 1$, and the Poincar\'e inequality in Lemma \ref{CE}. Consequently, \eqref{EE-2}-\eqref{EE-3}-\eqref{EE-4} tell us that
\begin{equation}\label{EE-5}
	\begin{aligned}
		| \nabla_A \big( A \cdot P_{T_{\Lambda_0}} ( \lambda [ f_R^\epsilon - \rho_R^\epsilon M_{\Lambda_0} ] ) \big) | \leq \tfrac{25 \sqrt[4]{3}}{\sqrt{\lambda_0}} \| \nabla_A ( \tfrac{f_R^\epsilon}{M_{\Lambda_0}} ) \|_{L^2_A (M_{\Lambda_0})} \,.
	\end{aligned}
\end{equation}
Plugging \eqref{EE-5} into \eqref{EE-0}, one has
\begin{equation}\label{EE-bnd}
	\begin{aligned}
		& - \iint_{\R^3 \times \mathrm{SO(3)}} \frac{1}{\epsilon} (\tfrac{ f_{R}^{\epsilon} }{ M_{\Lambda_{0}}} - \rho_{R}^\epsilon) L_R ( f_{R}^{\epsilon} - \rho_{R}^\epsilon M_{\Lambda_{0}} ) \d A \d x \\
		\leq & \tfrac{\nu_0}{c_1} \tfrac{1}{\epsilon} \tfrac{25 \sqrt[4]{3}}{\sqrt{\lambda_0}} \iint_{\R^3 \times \mathrm{SO(3)}} | \nabla_A ( \tfrac{f_R^\epsilon}{M_{\Lambda_0}} ) | \cdot \| \nabla_A ( \tfrac{f_R^\epsilon}{M_{\Lambda_0}} ) \|_{L^2_A (M_{\Lambda_0})} M_{\Lambda_0} \d A \d x \\
		\leq & \tfrac{25 \sqrt[4]{3} \nu_0 }{c_1 \sqrt{\lambda_0}} \tfrac{1}{\epsilon} \int_{\R^3} \| \nabla_A ( \tfrac{f_R^\epsilon}{M_{\Lambda_0}} ) \|^2_{L^2_A (M_{\Lambda_0})} \Big( \int_{\mathrm{SO}(3)} M_{\Lambda_0} \d A \Big)^\frac{1}{2} \d x \\
		= & \tfrac{25 \sqrt[4]{3} \nu_0 }{c_1 \sqrt{\lambda_0}} \tfrac{1}{\epsilon} \| \nabla_A ( \tfrac{f_R^\epsilon}{M_{\Lambda_0}} ) \|^2_{ L^2_{x,A} ( M_{ \Lambda_0 } ) } = \tfrac{25 \sqrt[4]{3} \nu_0 }{c_1 \sqrt{\lambda_0}} D_0 (t) \,,
	\end{aligned}
\end{equation}
where the functional $D_0 (t)$ is defined in \eqref{Ekt}.

{\bf Step 5. Estimates of source term: $R (f_1)$.} By the H\"older inequality, it is easy to see that
\begin{equation}\label{EST-bnd}
	\begin{aligned}
		| \iint_{\R^3 \times \mathrm{SO(3)}} R (f_1) (\tfrac{ f_{R}^{\epsilon} }{ M_{\Lambda_{0}}} - \rho_{R}^\epsilon) \d A \d x | \leq & \| \tfrac{R (f_1)}{M_{ \Lambda_0 }} \|_{L^2_{x,A} ( M_{ \Lambda_0 } ) } \| \tfrac{ f_{R}^{\epsilon} }{ M_{\Lambda_{0}}} - \rho_{R}^\epsilon \|_{L^2_{x,v} ( M_{ \Lambda_0 } ) } \\
		\leq & \| \tfrac{R (f_1)}{M_{ \Lambda_0 }} \|_{L^2_{x,A} ( M_{ \Lambda_0 } ) } E_0^\frac{1}{2} (t) \,,
	\end{aligned}
\end{equation}
where the functional $E_0 (t)$ is given in \eqref{Ekt}.

{\bf Step 6. Estimates of nonlinear term: $ \widetilde{Q} (f_R^\epsilon) $.} Recall the definition of $ \widetilde{Q} (f_R^\epsilon) $ in \eqref{RWQ}, hence,
\begin{equation*}
	\begin{aligned}
		\widetilde{Q} ( f_{R}^{\epsilon} ) = & - \frac{1}{\epsilon} \nu_{0} \nabla_{A} \cdot [ f_{R}^{\epsilon} \nabla_{A} ( A \cdot \Lambda [ f_{0} + \epsilon f_{1} + \epsilon f_{R}^{\epsilon} ] - A \cdot \Lambda [f_0] ) ] \\
		&- \frac{1}{\epsilon} \nu_0 \nabla_A \cdot [ f_1 \nabla_A ( A \cdot \Lambda [f_0 + \epsilon f_1 + \epsilon f_R^\epsilon] - A \cdot \Lambda [f_0] ) ] \\ 
		& - \frac{1}{\epsilon^2} \nu_0 \nabla_A \cdot \big[ f_0 \nabla_A \big( A \cdot \{ \Lambda [f_0 + \epsilon f_1 + \epsilon f_R^\epsilon] - \Lambda [f_0] - \epsilon \tfrac{\d}{\d \epsilon} |_{\epsilon = 0} \Lambda [f_0 + \epsilon g]_{g = f_1 + f_R^\epsilon} \} \big) \big] \,.
	\end{aligned}
\end{equation*}
By employing the Taylor expansion for $\Lambda [ f_0 + \epsilon f_1 + \epsilon f_R^\epsilon ]$ with respect to the parameter $\epsilon$, and applying Lemma \ref{Lmm-Lambda}, it is implied that
\begin{equation*}
	\begin{aligned}
		\Lambda [ f_{0} + \epsilon f_{1} + \epsilon f_{R}^{\epsilon} ] - \Lambda [f_0] = & \epsilon \tfrac{\d}{\d \epsilon} |_{\epsilon = 0} \Lambda [ f_0 + \epsilon g ]_{g = f_1 + f_R^\epsilon} + \mathcal{O} (\epsilon^2) \\
		= & \epsilon (c_1 \rho_0)^{-1} P_{T_{ \Lambda_0 }} ( \lambda [ f_1 + f_R^\epsilon ] ) + \mathcal{O} (\epsilon^2) \,,
	\end{aligned}
\end{equation*}
and
	\begin{align*}
		& \Lambda [f_0 + \epsilon f_1 + \epsilon f_R^\epsilon] - \Lambda [f_0] - \epsilon \tfrac{\d}{\d \epsilon} |_{\epsilon = 0} \Lambda [f_0 + \epsilon g]_{g = f_1 + f_R^\epsilon} \\
		= & \epsilon^2 \tfrac{\d^2}{\d \epsilon^2} \Lambda [ f_0 + \epsilon g ]_{g = f_1 + f_R^\epsilon} + \mathcal{O} (\epsilon^3) \\
		= & - \tfrac{\epsilon^2}{2 ( c_1 \rho_0 )^2} \lambda [ f_1 + f_R^\epsilon ] ( \lambda [ f_1 + f_R^\epsilon ]^\top \Lambda_0 + \Lambda_0^\top \lambda [ f_1 + f_R^\epsilon ] ) \\
		& - \tfrac{\epsilon^2}{( c_1 \rho_0 )^3} \lambda [ f_1 + f_R^\epsilon ]^\top \lambda [ f_1 + f_R^\epsilon ] + \tfrac{3 \epsilon^2}{4 ( c_1 \rho_0 )^3} \big( \lambda [ f_1 + f_R^\epsilon ]^\top \Lambda_0 + \Lambda_0^\top \lambda [ f_1 + f_R^\epsilon ] \big)^2 + \mathcal{O} (\epsilon^3) \,.
	\end{align*}
As a consequence,
\begin{equation}\label{Q-decomp}
	\begin{aligned}
		\widetilde{Q} (f_R^\epsilon) = & \underbrace{ - \tfrac{\nu_0}{c_1 \rho_0} \nabla_A \cdot \Big[ ( f_1 + f_R^\epsilon ) \nabla_A \big( A \cdot P_{ T_{ \Lambda_0 } } ( \lambda [ f_1 + f_R^\epsilon ] ) \big) \Big] }_{: = Q_1} \\
		& + \underbrace{ \tfrac{\nu_0}{ 2 c_1^2 \rho_0 } \nabla_A \cdot \Big\{ M_{ \Lambda_0 } \nabla_A \Big[ A \cdot \lambda [ f_1 + f_R^\epsilon ] \big( \lambda [ f_1 + f_R^\epsilon ]^\top \Lambda_0 + \Lambda_0^\top \lambda [ f_1 + f_R^\epsilon ] \big) \Big] \Big\} }_{:= Q_2} \\
		& + \underbrace{ \tfrac{\nu_0}{c_1^3 \rho_0^2} \nabla_A \cdot \Big\{ M_{ \Lambda_0 } \nabla_A \Big[ A \cdot \lambda [ f_1 + f_R^\epsilon ]^\top \lambda [ f_1 + f_R^\epsilon ] \Big] \Big\} }_{:= Q_3} \\
		& \underbrace{ - \tfrac{3 \nu_0}{4 c_1^3 \rho_0^2} \nabla_A \cdot \Big\{ M_{ \Lambda_0 } \nabla_A \Big[ A \cdot \big( \lambda [ f_1 + f_R^\epsilon ]^\top \Lambda_0 + \Lambda_0^\top \lambda [ f_1 + f_R^\epsilon ] \big)^2 \Big] \Big\} }_{:= Q_4} + \mathcal{O} (\epsilon) \,.
	\end{aligned}
\end{equation}
Notice that the goal of this paper is to justify the limit of the SOKB system \eqref{KM} as $\epsilon \to 0$. Then the term $\mathcal{O} (\epsilon)$ in \eqref{Q-decomp} is an infinitesimal small quantity as $\epsilon \to 0$, the effect of which is very small. The majority parts in \eqref{Q-decomp} are the quantities $Q_i$ ($i = 1,2,3,4$). For simplicity and convenience of computations, we will neglect the infinitesimal small quantity $\mathcal{O} (\epsilon)$ in the following calculations.

\underline{\em Case 6.1. Control of $Q_1$.} The integration by parts over $A \in \mathrm{SO(3)}$ reduces to
\begin{equation}\label{Q1-1}
	\begin{aligned}
		& \iint_{\R^3 \times \mathrm{SO(3)}} Q_1 ( \tfrac{ f_R^\epsilon }{ M_{ \Lambda_0 } } - \rho_R^\epsilon ) \d A \d x \\
		= & \iint_{\R^3 \times \mathrm{SO(3)}} - \tfrac{\nu_0}{c_1 \rho_0} \nabla_A \cdot \Big[ ( f_1 + f_R^\epsilon ) \nabla_A \big( A \cdot P_{ T_{ \Lambda_0 } } ( \lambda [ f_1 + f_R^\epsilon ] ) \big) \Big] ( \tfrac{ f_R^\epsilon }{ M_{ \Lambda_0 } } - \rho_R^\epsilon ) \d A \d x \\
		= & \iint_{\R^3 \times \mathrm{SO(3)}} \tfrac{\nu_0}{c_1 \rho_0} \Big[ ( f_1 + f_R^\epsilon ) \nabla_A \big( A \cdot P_{ T_{ \Lambda_0 } } ( \lambda [ f_1 + f_R^\epsilon ] ) \big) \Big] \cdot \nabla_A ( \tfrac{ f_R^\epsilon }{ M_{ \Lambda_0 } } ) \d A \d x \\
		\lesssim & \iint_{\R^3 \times \mathrm{SO(3)}} ( |\tfrac{f_1}{M_{\Lambda_0}}| + | \tfrac{f_R^\epsilon}{ M_{ \Lambda_0 } } - \rho_R^\epsilon | + | \rho_R^\epsilon | ) \big| \nabla_A \big( A \cdot P_{ T_{ \Lambda_0 } } ( \lambda [ f_1 + f_R^\epsilon ] ) \big) \big| | \nabla_A ( \tfrac{ f_R^\epsilon }{ M_{ \Lambda_0 } } ) | M_{ \Lambda_0 } \d A \d x \,,
	\end{aligned}
\end{equation}
due to $\inf_{(t,x)\in [0, T] \times \R^3} \rho_0 (t,x) > 0$ by Theorem \ref{LWP}. It is easy to follow from Lemma \ref{Lmm-Proj-SO3} that
\begin{equation*}
	\begin{aligned}
		& \nabla_A \big( A \cdot P_{ T_{ \Lambda_0 } } ( \lambda [ f_1 + f_R^\epsilon ] ) \big) \\
		= & \tfrac{1}{4} \Big\{ \lambda [ f_1 + f_R^\epsilon ] - \Lambda_0 \lambda [ f_1 + f_R^\epsilon ]^\top \Lambda_0 - A \lambda [ f_1 + f_R^\epsilon ]^\top A + A \Lambda_0^\top \lambda [ f_1 + f_R^\epsilon ] \Lambda_0^\top A  \Big\} \,,
	\end{aligned}
\end{equation*}
which, together with $|A| = |\Lambda_0| = \sqrt{3}$ for $A, \Lambda_0 \in \mathrm{SO(3)}$, infers that
\begin{equation}\label{Q1-2}
	\begin{aligned}
		\big| \nabla_A \big( A \cdot P_{ T_{ \Lambda_0 } } ( \lambda [ f_1 + f_R^\epsilon ] ) \big) \big| \lesssim | \lambda [ f_1 + f_R^\epsilon ] | \,.
	\end{aligned}
\end{equation}
Recall the definition of $\lambda [f]$ in \eqref{KM}, one has
\begin{equation}\label{Q1-3}
	\begin{aligned}
		| \lambda [ f_1 + f_R^\epsilon ] | = & | \int_{\mathrm{SO}(3)} A ( f_1 + f_R^\epsilon ) \d A | \\
		= & | \int_{\mathrm{SO}(3)} A ( \tfrac{f_1}{ M_{ \Lambda_0 } } + \tfrac{f_R^\epsilon}{ M_{ \Lambda_0 } } - \rho_R^\epsilon + \rho_R^\epsilon ) M_{ \Lambda_0 } \d A | \\
		\leq & \int_{\mathrm{SO}(3)} |A| \big( | \tfrac{f_1}{ M_{ \Lambda_0 } } | + | \tfrac{f_R^\epsilon}{ M_{ \Lambda_0 } } - \rho_R^\epsilon | + | \rho_R^\epsilon | \big) M_{ \Lambda_0 } \d A \\
		\lesssim & \| \tfrac{f_1}{ M_{ \Lambda_0 } } \|_{L^2_A ( M_{ \Lambda_0 } )} + \| \tfrac{f_R^\epsilon}{ M_{ \Lambda_0 } } - \rho_R^\epsilon \|_{ L^2_A ( M_{ \Lambda_0 } ) } + | \rho_R^\epsilon | \,.
	\end{aligned}
\end{equation}
Then the relations \eqref{Q1-1}, \eqref{Q1-2} and \eqref{Q1-3} show that
	\begin{align}\label{Q1-bnd}
		\no & | \iint_{\R^3 \times \mathrm{SO(3)}} Q_1 ( \tfrac{ f_R^\epsilon }{ M_{ \Lambda_0 } } - \rho_R^\epsilon ) \d A \d x | \\
		\no \lesssim & \int_{\R^3} ( \| \tfrac{f_1}{ M_{ \Lambda_0 } } \|_{L^2_A ( M_{ \Lambda_0 } )} + \| \tfrac{f_R^\epsilon}{ M_{ \Lambda_0 } } - \rho_R^\epsilon \|_{ L^2_A ( M_{ \Lambda_0 } ) } + | \rho_R^\epsilon | )^2 \| \nabla_A ( \tfrac{f_R^\epsilon}{ M_{ \Lambda_0 } } ) \|_{L^2_A ( M_{ \Lambda_0 } )} \d x \\
		\no \lesssim & ( \| \tfrac{f_1}{ M_{ \Lambda_0 } } \|_{L^\infty_x L^2_A ( M_{ \Lambda_0 } )} + \| \tfrac{f_R^\epsilon}{ M_{ \Lambda_0 } } - \rho_R^\epsilon \|_{ L^\infty_x L^2_A ( M_{ \Lambda_0 } ) } + \| \rho_R^\epsilon \|_{ L^\infty_x } ) \\
		\no & \times ( \| \tfrac{f_1}{ M_{ \Lambda_0 } } \|_{L^2_{x,A} ( M_{ \Lambda_0 } )} + \| \tfrac{f_R^\epsilon}{ M_{ \Lambda_0 } } - \rho_R^\epsilon \|_{ L^2_{x,A} ( M_{ \Lambda_0 } ) } + \| \rho_R^\epsilon \|_{L^2_x} ) \| \nabla_A ( \tfrac{f_R^\epsilon}{ M_{ \Lambda_0 } } ) \|_{L^2_{x,A} ( M_{ \Lambda_0 } )} \\
		\lesssim & ( \| \tfrac{f_1}{ M_{ \Lambda_0 } } \|_{H^2_x L^2_A ( M_{ \Lambda_0 } )} + \| \tfrac{f_R^\epsilon}{ M_{ \Lambda_0 } } - \rho_R^\epsilon \|_{ H^2_x L^2_A ( M_{ \Lambda_0 } ) } + \| \rho_R^\epsilon \|_{ H^2_x } ) \\
		\no & \times ( \| \tfrac{f_1}{ M_{ \Lambda_0 } } \|_{L^2_{x,A} ( M_{ \Lambda_0 } )} + \| \tfrac{f_R^\epsilon}{ M_{ \Lambda_0 } } - \rho_R^\epsilon \|_{ L^2_{x,A} ( M_{ \Lambda_0 } ) } + \| \rho_R^\epsilon \|_{L^2_x} ) \| \nabla_A ( \tfrac{f_R^\epsilon}{ M_{ \Lambda_0 } } ) \|_{L^2_{x,A} ( M_{ \Lambda_0 } )} \\
		\no \lesssim & \sqrt{\epsilon} \big( \| \tfrac{f_1}{M_{\Lambda_0}} \|_{H^2_x L^2_A ( M_{\Lambda_0} )} + E_2^\frac{1}{2} (t) \big) \big( \| \tfrac{f_1}{M_{\Lambda_0}} \|_{L^2_{x,A} ( M_{\Lambda_0} )} + E_0^\frac{1}{2} (t) \big) D_0^\frac{1}{2} (t) \,,
	\end{align}
where the last second inequality is derived from the Sobolev embedding $H^2_x \hookrightarrow L^\infty_x$, and the functionals $E_k (t)$, $D_k (t)$ are defined in \eqref{Ekt}.

\underline{\em Case 6.2. Control of $Q_2$.} Following the similar arguments in \eqref{Q1-1}, one has
\begin{equation*}
	\begin{aligned}
		& | \iint_{\R^3 \times \mathrm{SO(3)}} Q_2 ( \tfrac{ f_R^\epsilon }{ M_{ \Lambda_0 } } - \rho_R^\epsilon ) \d A \d x | \\
		\lesssim & \iint_{\R^3 \times \mathrm{SO(3)}} \Big| \nabla_A \Big[ A \cdot \lambda [ f_1 + f_R^\epsilon ] \big( \lambda [ f_1 + f_R^\epsilon ]^\top \Lambda_0 + \Lambda_0^\top \lambda [ f_1 + f_R^\epsilon ] \big) \Big] \Big| \cdot | \nabla_A ( \tfrac{f_R^\epsilon}{ M_{ \Lambda_0 } } ) | M_{\Lambda_0} \d A \d x \,.
	\end{aligned}
\end{equation*}
From the similar arguments in \eqref{Q1-2}, there holds
\begin{equation*}
	\begin{aligned}
		\Big| \nabla_A \Big[ A \cdot \lambda [ f_1 + f_R^\epsilon ] \big( \lambda [ f_1 + f_R^\epsilon ]^\top \Lambda_0 + \Lambda_0^\top \lambda [ f_1 + f_R^\epsilon ] \big) \Big] \Big| \lesssim | \lambda [ f_1 + f_R^\epsilon ] |^2 \,,
	\end{aligned}
\end{equation*}
which, together with \eqref{Q1-3} and similar estimates in \eqref{Q1-bnd}, means that
\begin{equation}\label{Q2-bnd}
	\begin{aligned}
		& | \iint_{\R^3 \times \mathrm{SO(3)}} Q_2 ( \tfrac{ f_R^\epsilon }{ M_{ \Lambda_0 } } - \rho_R^\epsilon ) \d A \d x | \\
		\lesssim & \iint_{\R^3 \times \mathrm{SO(3)}} \big( \| \tfrac{f_1}{ M_{ \Lambda_0 } } \|_{L^2_A ( M_{ \Lambda_0 } )} + \| \tfrac{f_R^\epsilon}{ M_{ \Lambda_0 } } - \rho_R^\epsilon \|_{ L^2_A ( M_{ \Lambda_0 } ) } + | \rho_R^\epsilon | \big)^2 | \nabla_A ( \tfrac{f_R^\epsilon}{ M_{ \Lambda_0 } } ) | M_{\Lambda_0} \d A \d x \\
		\lesssim & \sqrt{\epsilon} \big( \| \tfrac{f_1}{M_{\Lambda_0}} \|_{H^2_x L^2_A ( M_{\Lambda_0} )} + E_2^\frac{1}{2} (t) \big) \big( \| \tfrac{f_1}{M_{\Lambda_0}} \|_{L^2_{x,A} ( M_{\Lambda_0} )} + E_0^\frac{1}{2} (t) \big) D_0^\frac{1}{2} (t) \,.
	\end{aligned}
\end{equation}

\underline{\em Case 6.3. Control of $Q_3$.} Following the similar arguments in \eqref{Q2-bnd}, one has
\begin{equation}\label{Q3-bnd}
	\begin{aligned}
		& | \iint_{\R^3 \times \mathrm{SO(3)}} Q_3 ( \tfrac{ f_R^\epsilon }{ M_{ \Lambda_0 } } - \rho_R^\epsilon ) \d A \d x | \\
		\lesssim & \sqrt{\epsilon} \big( \| \tfrac{f_1}{M_{\Lambda_0}} \|_{H^2_x L^2_A ( M_{\Lambda_0} )} + E_2^\frac{1}{2} (t) \big) \big( \| \tfrac{f_1}{M_{\Lambda_0}} \|_{L^2_{x,A} ( M_{\Lambda_0} )} + E_0^\frac{1}{2} (t) \big) D_0^\frac{1}{2} (t) \,.
	\end{aligned}
\end{equation}

\underline{\em Case 6.4. Control of $Q_4$.} Following the similar arguments in \eqref{Q2-bnd}, one has
\begin{equation}\label{Q4-bnd}
	\begin{aligned}
		& | \iint_{\R^3 \times \mathrm{SO(3)}} Q_4 ( \tfrac{ f_R^\epsilon }{ M_{ \Lambda_0 } } - \rho_R^\epsilon ) \d A \d x | \\
		\lesssim & \sqrt{\epsilon} \big( \| \tfrac{f_1}{M_{\Lambda_0}} \|_{H^2_x L^2_A ( M_{\Lambda_0} )} + E_2^\frac{1}{2} (t) \big) \big( \| \tfrac{f_1}{M_{\Lambda_0}} \|_{L^2_{x,A} ( M_{\Lambda_0} )} + E_0^\frac{1}{2} (t) \big) D_0^\frac{1}{2} (t) \,.
	\end{aligned}
\end{equation}

In summary, the relations \eqref{Q-decomp}, \eqref{Q1-bnd}, \eqref{Q2-bnd}, \eqref{Q3-bnd} and \eqref{Q4-bnd} tell us that
\begin{equation}\label{ENT-bnd}
	\begin{aligned}
		& | \iint_{\R^3 \times \mathrm{SO(3)}} \widetilde{Q} (f_R^\epsilon) ( \tfrac{ f_R^\epsilon }{ M_{ \Lambda_0 } } - \rho_R^\epsilon ) \d A \d x | \\
		\lesssim & \sqrt{\epsilon} \big( \| \tfrac{f_1}{M_{\Lambda_0}} \|_{H^2_x L^2_A ( M_{\Lambda_0} )} + E_2^\frac{1}{2} (t) \big) \big( \| \tfrac{f_1}{M_{\Lambda_0}} \|_{L^2_{x,A} ( M_{\Lambda_0} )} + E_0^\frac{1}{2} (t) \big) D_0^\frac{1}{2} (t) \,.
	\end{aligned}
\end{equation}

It therefore follows from substituting the bounds \eqref{MiTE-bnd}, \eqref{MaTE-bnd}, \eqref{Coer-Bnd}, \eqref{EE-bnd}, \eqref{EST-bnd} and \eqref{ENT-bnd} into the equation \eqref{L2-1} that
\begin{equation*}
	\begin{aligned}
		& \tfrac{1}{2} \tfrac{\d}{\d t} \big( \| \tfrac{f_R^\epsilon}{ M_{ \Lambda_0 } } - \rho_R^\epsilon \|^2_{L^2_{x,A} ( M_{ \Lambda_0 } )} + \| \rho_R^\epsilon \|^2_{L^2_x} \big) + d_\star D_0 (t) \\
		\lesssim & \| \partial_{t,x} \Lambda_0 \|_{H^2_x} E_0 (t) + \big( \| \tfrac{\nabla_x f_1}{ M_{ \Lambda_0 } } \|_{L^2_{x,A} ( M_{ \Lambda_0 } )} + \| \tfrac{ R ( f_1 ) }{ M_{ \Lambda_0 } } \|_{L^2_{x,A} ( M_{ \Lambda_0 } )} \big) E_0^\frac{1}{2} (t) \\
		& + \sqrt{\epsilon} \big( \| \tfrac{f_1}{M_{\Lambda_0}} \|_{H^2_x L^2_A ( M_{\Lambda_0} )} + E_2^\frac{1}{2} (t) \big) \big( \| \tfrac{f_1}{M_{\Lambda_0}} \|_{L^2_{x,A} ( M_{\Lambda_0} )} + E_0^\frac{1}{2} (t) \big) D_0^\frac{1}{2} (t) \,,
	\end{aligned}
\end{equation*}
where the constant $d_\star = d - \frac{25 \sqrt[4]{3} \nu_0 }{ c_1 \lambda_0 } > 0$. Then, by the Young's inequality, there holds
\begin{equation}\label{L2-bnd}
	\begin{aligned}
		& \tfrac{\d}{\d t} E_0 (t) + d_\star D_0 (t) \\
		\lesssim & \| \partial_{t,x} \Lambda_0 \|_{H^2_x} E_0 (t) + \big( \| \tfrac{\nabla_x f_1}{ M_{ \Lambda_0 } } \|_{L^2_{x,A} ( M_{ \Lambda_0 } )} + \| \tfrac{ R ( f_1 ) }{ M_{ \Lambda_0 } } \|_{L^2_{x,A} ( M_{ \Lambda_0 } )} \big) E_0^\frac{1}{2} (t) \\
		& + \epsilon \big( \| \tfrac{f_1}{M_{\Lambda_0}} \|^2_{H^2_x L^2_A ( M_{\Lambda_0} )} + E_2 (t) \big) \big( \| \tfrac{ f_1 }{ M_{\Lambda_0} } \|^2_{L^2_{x,A} ( M_{\Lambda_0} )} + E_0 (t) \big) \\
		\lesssim & (1 + \| \nabla_x \Lambda_0 \|_{H^2_x} + \| \partial_t \Lambda_0 \|_{H^2_x} ) ( 1 + \| \tfrac{f_1}{ M_{ \Lambda_0 } } \|_{H^1_x L^2_A ( M_{ \Lambda_0 } )} + \| \tfrac{ R ( f_1 ) }{ M_{ \Lambda_0 } } \|_{L^2_{x,A} ( M_{ \Lambda_0 } )} ) \\
		& \times ( E_0 (t) + E_0^\frac{1}{2} (t) + \epsilon E_0 (t) E_2 (t) ) \\
		\lesssim & C( \| \nabla_x \rho_0^{in} \|_{H^3_x}, \| \nabla_x \Lambda_0^{in} \|_{H^3_x} ) ( E_0 (t) + E_0^\frac{1}{2} (t) + \epsilon E_0 (t) E_2 (t) ) \,,
	\end{aligned}
\end{equation}
where the last inequality is derived from Theorem \ref{LWP} and Lemma \ref{Lmm-f1}.

\subsubsection{Higher order derivatives estimates}\label{Subsec:DerEst}

We first represent the micro-equation \eqref{micro-eq} as
\begin{equation}\label{micro-eq-2}
	\begin{aligned}
		& ( \partial_t + A e_1 \cdot  \nabla_x ) ( \tfrac{ f_R^\epsilon }{ M_{ \Lambda_0 } } - \rho_R^\epsilon ) + ( \partial_t + A e_1 \cdot  \nabla_x )  \rho_R^\epsilon - \frac{1}{\epsilon} \tfrac{1}{ M_{ \Lambda_0 } } \L_{ M_{ \Lambda_0 } } f_R^\epsilon \\
		= & \frac{1}{\epsilon} \tfrac{1}{ M_{ \Lambda_0 } } L_R ( f_R^\epsilon - \rho_R^\epsilon M_{ \Lambda_0 } ) - \tfrac{ f_R^\epsilon }{ M_{ \Lambda_0 } } \tfrac{ ( \partial_t + A e_1 \cdot  \nabla_x ) M_{ \Lambda_0 } }{ M_{ \Lambda_0 } }  + \tfrac{1}{ M_{ \Lambda_0 } } \widetilde{Q} ( f_R^\epsilon ) + \tfrac{1}{ M_{ \Lambda_0 } } R ( f_1 ) \,.
	\end{aligned}
\end{equation}
For any multi-index $\ss = ( \ss_1, \ss_2, \ss_3 ) \in \mathbb{N}^3$ with $1 \leq |\ss| = k \leq s$ ($s \geq 2$), applying higher order spatial derivative operator $ \partial_x^{\ss} = \frac{\partial^{|\ss|}}{\partial x_1^{\ss_1} \partial x_2^{\ss_2} \partial x_3^{\ss_3} } $ on the micro-equation \eqref{micro-eq-2}, one gets
\begin{equation}\label{Hmicro-eq}
	\begin{aligned}
		& ( \partial_t + A e_1 \cdot  \nabla_x ) \partial_x^{\ss} ( \tfrac{ f_R^\epsilon }{ M_{ \Lambda_0 } } - \rho_R^\epsilon ) + ( \partial_t + A e_1 \cdot  \nabla_x ) \partial_x^{\ss} \rho_R^\epsilon - \frac{1}{\epsilon} \partial_x^{\ss} \big[ \tfrac{1}{ M_{ \Lambda_0 } } \L_{ M_{ \Lambda_0 } } f_R^\epsilon \big] \\
		= & \frac{1}{\epsilon} \partial_x^{\ss} \big[ \tfrac{1}{ M_{ \Lambda_0 } } L_R ( f_R^\epsilon - \rho_R^\epsilon M_{ \Lambda_0 } ) \big] - \partial_x^{\ss} \big[ \tfrac{ f_R^\epsilon }{ M_{ \Lambda_0 } } \tfrac{ ( \partial_t + A e_1 \cdot  \nabla_x ) M_{ \Lambda_0 } }{ M_{ \Lambda_0 } } \big]  \\
		& + \partial_x^{\ss} \big[ \tfrac{1}{ M_{ \Lambda_0 } } \widetilde{Q} ( f_R^\epsilon ) \big] + \partial_x^{\ss} \big[ \tfrac{1}{ M_{ \Lambda_0 } } R ( f_1 ) \big] \,.
	\end{aligned}
\end{equation}
Moreover, the macro-equation \eqref{macro-eq} indicates that
\begin{equation}\label{Hmacro-eq}
	\begin{aligned}
		\partial_{t} \partial_x^{\ss} \rho_{R}^{\epsilon} + \int_{\mathrm{SO(3)}} A e_{1} \cdot \nabla_{x} \partial_x^{\ss} f_{R}^{\epsilon} \d A = \int_{\mathrm{SO}(3)} \partial_x^{\ss} R (f_1) \d A \,,
	\end{aligned}
\end{equation}
where
\begin{equation}\label{HRf1-int}
	\begin{aligned}
		\int_{\mathrm{SO}(3)} \partial_x^{\ss} R (f_1) \d A = \int_{\mathrm{SO}(3)} A e_1 \cdot \nabla_x \partial_x^{\ss} f_1 \d A \,.
	\end{aligned}
\end{equation}
It follows from multiplying \eqref{Hmicro-eq} by $ \partial_x^{\ss} \big( \tfrac{ f_R^\epsilon }{ M_{ \Lambda_0 } } - \rho_R^\epsilon \big) M_{ \Lambda_0 } $ and integrating over $(x,A) \in \R^3 \times \mathrm{SO(3)}$ that
\begin{equation}\label{HigherOr-der}
 	\begin{aligned}
 		& \iint_{\R^3 \times \mathrm{SO(3)}} ( \partial_t + A e_1 \cdot  \nabla_x ) \partial_x^{\ss} ( \tfrac{ f_R^\epsilon }{ M_{ \Lambda_0 } } - \rho_R^\epsilon ) \cdot \partial_x^{\ss} \big( \tfrac{ f_R^\epsilon }{ M_{ \Lambda_0 } } - \rho_R^\epsilon \big) M_{ \Lambda_0 } \d A \d x \\
 		& + \iint_{\R^3 \times \mathrm{SO(3)}} ( \partial_t + A e_1 \cdot  \nabla_x ) \partial_x^{\ss} \rho_R^\epsilon \cdot \partial_x^{\ss} \big( \tfrac{ f_R^\epsilon }{ M_{ \Lambda_0 } } - \rho_R^\epsilon \big) M_{ \Lambda_0 } \d A \d x \\
 		& + \iint_{\R^3 \times \mathrm{SO(3)}} - \frac{1}{\epsilon} \partial_x^{\ss} \big[ \tfrac{1}{ M_{ \Lambda_0 } } \L_{ M_{ \Lambda_0 } } f_R^\epsilon \big] \cdot \partial_x^{\ss} \big( \tfrac{ f_R^\epsilon }{ M_{ \Lambda_0 } } - \rho_R^\epsilon \big) M_{ \Lambda_0 } \d A \d x \\
 		= & - \iint_{\R^3 \times \mathrm{SO(3)}} \frac{1}{\epsilon} \partial_x^{\ss} \big[ \tfrac{1}{ M_{ \Lambda_0 } } L_R ( f_R^\epsilon - \rho_R^\epsilon M_{ \Lambda_0 } ) \big] \cdot \partial_x^{\ss} \big( \tfrac{ f_R^\epsilon }{ M_{ \Lambda_0 } } - \rho_R^\epsilon \big) M_{ \Lambda_0 } \d A \d x \\
 		& + \iint_{\R^3 \times \mathrm{SO(3)}} - \partial_x^{\ss} \big[ \tfrac{ f_R^\epsilon }{ M_{ \Lambda_0 } } \tfrac{ ( \partial_t + A e_1 \cdot  \nabla_x ) M_{ \Lambda_0 } }{ M_{ \Lambda_0 } } \big] \cdot \partial_x^{\ss} \big( \tfrac{ f_R^\epsilon }{ M_{ \Lambda_0 } } - \rho_R^\epsilon \big) M_{ \Lambda_0 } \d A \d x \\
 		& + \iint_{\R^3 \times \mathrm{SO(3)}} \partial_x^{\ss} \big[ \tfrac{1}{ M_{ \Lambda_0 } } R ( f_1 ) \big] \cdot \partial_x^{\ss} \big( \tfrac{ f_R^\epsilon }{ M_{ \Lambda_0 } } - \rho_R^\epsilon \big) M_{ \Lambda_0 } \d A \d x \\
 		& + \iint_{\R^3 \times \mathrm{SO(3)}} \partial_x^{\ss} \big[ \tfrac{1}{ M_{ \Lambda_0 } } \widetilde{Q} ( f_R^\epsilon ) \big] \cdot \partial_x^{\ss} \big( \tfrac{ f_R^\epsilon }{ M_{ \Lambda_0 } } - \rho_R^\epsilon \big) M_{ \Lambda_0 } \d A \d x \,.
 	\end{aligned}
 \end{equation}

{\bf Step 1. Estimates of transport effect for micro part: $ ( \partial_t + A e_1 \cdot  \nabla_x ) \partial_x^{\ss} ( \tfrac{ f_R^\epsilon }{ M_{ \Lambda_0 } } - \rho_R^\epsilon ) $.} Following the similar arguments in \eqref{MiTE-bnd}, one has
\begin{equation}\label{HETEMi-bnd}
	\begin{aligned}
		& \iint_{\R^3 \times \mathrm{SO(3)}} ( \partial_t + A e_1 \cdot  \nabla_x ) \partial_x^{\ss} ( \tfrac{ f_R^\epsilon }{ M_{ \Lambda_0 } } - \rho_R^\epsilon ) \cdot \partial_x^{\ss} \big( \tfrac{ f_R^\epsilon }{ M_{ \Lambda_0 } } - \rho_R^\epsilon \big) M_{ \Lambda_0 } \d A \d x \\
		= & \tfrac{1}{2} \tfrac{\d}{\d t} \| \partial_x^{\ss} ( \tfrac{ f_R^\epsilon }{ M_{ \Lambda_0 } } - \rho_R^\epsilon ) \|^2_{ L^2_{x,A} ( M_{ \Lambda_0 } ) } \\
		& - \iint_{\R^3 \times \mathrm{SO(3)}} \big[ \partial_x^{\ss} \big( \tfrac{ f_R^\epsilon }{ M_{ \Lambda_0 } } - \rho_R^\epsilon \big) \big]^2 ( \tfrac{1}{2} \partial_t + A e_1 \cdot \nabla_x ) M_{ \Lambda_0 } \d A \d x \\
		\geq & \tfrac{1}{2} \tfrac{\d}{\d t} \| \partial_x^{\ss} ( \tfrac{ f_R^\epsilon }{ M_{ \Lambda_0 } } - \rho_R^\epsilon ) \|^2_{ L^2_{x,A} ( M_{ \Lambda_0 } ) } - \tfrac{ \sqrt{3} \nu_0 }{d} \| \partial_{t,x} \Lambda_0 \|_{ H^2_x } E_k (t) \,,
	\end{aligned}
\end{equation}
where the functional $E_k (t)$ can be seen in \eqref{Ekt}.

{\bf Step 2. Estimates of transport effect for macro part: $ ( \partial_t + A e_1 \cdot  \nabla_x ) \partial_x^{\ss} \rho_R^\epsilon $.} Observe that
\begin{equation}\label{HETE-1}
	\begin{aligned}
		& \iint_{\R^3 \times \mathrm{SO(3)}} ( \partial_t + A e_1 \cdot  \nabla_x ) \partial_x^{\ss} \rho_R^\epsilon \cdot \partial_x^{\ss} \big( \tfrac{ f_R^\epsilon }{ M_{ \Lambda_0 } } - \rho_R^\epsilon \big) M_{ \Lambda_0 } \d A \d x \\
		= & \iint_{\R^3 \times \mathrm{SO(3)}} ( \partial_t + A e_1 \cdot  \nabla_x ) \partial_x^{\ss} \rho_R^\epsilon \cdot \partial_x^{\ss} \big( f_R^\epsilon - \rho_R^\epsilon M_{ \Lambda_0 } \big) \d A \d x \\
		& - \sum_{0 \neq \ss' \leq \ss} C_{\ss}^{\ss'} \iint_{\R^3 \times \mathrm{SO(3)}} ( \partial_t + A e_1 \cdot  \nabla_x ) \partial_x^{\ss} \rho_R^\epsilon \cdot \partial_x^{\ss - \ss'} ( \tfrac{ f_R^\epsilon }{ M_{ \Lambda_0 } } - \rho_R^\epsilon ) \partial_x^{\ss'} M_{ \Lambda_0 } \d A \d x \,.
	\end{aligned}
\end{equation}
Note that the cancellation
\begin{equation*}
	\begin{aligned}
		\int_{\mathrm{SO}(3)} \partial_x^{\ss} ( f_R^\epsilon - \rho_R^\epsilon M_{ \Lambda_0 } ) \d A = 0 \,,
	\end{aligned}
\end{equation*}
which means that
\begin{equation}\label{HETE-2}
	\begin{aligned}
		\iint_{\R^3 \times \mathrm{SO(3)}} \partial_t \partial_x^{\ss} \rho_R^\epsilon \cdot \partial_x^{\ss} ( f_R^\epsilon - \rho_R^\epsilon M_{ \Lambda_0 } ) \d A \d x = 0 \,.
	\end{aligned}
\end{equation}
Moreover, the integration by parts over $x \in \R^3$ shows that
\begin{equation}\label{HETE-3}
	\begin{aligned}
		& \iint_{\R^3 \times \mathrm{SO(3)}} A e_1 \cdot  \nabla_x \partial_x^{\ss} \rho_R^\epsilon \cdot \partial_x^{\ss} \big( f_R^\epsilon - \rho_R^\epsilon M_{ \Lambda_0 } \big) \d A \d x \\
		= & - \iint_{\R^3 \times \mathrm{SO(3)}} \partial_x^{\ss} \rho_R^\epsilon A e_1 \cdot  \nabla_x \partial_x^{\ss} f_R^\epsilon \d A \d x - \iint_{\R^3 \times \mathrm{SO(3)}} A e_1 \cdot  \nabla_x \partial_x^{\ss} \rho_R^\epsilon \partial_x^{\ss} \rho_R^\epsilon M_{ \Lambda_0 } \d A \d x \\
		& - \sum_{0 \neq \ss' \leq \ss} C_{\ss}^{\ss'} \iint_{\R^3 \times \mathrm{SO(3)}} A e_1 \cdot  \nabla_x \partial_x^{\ss} \rho_R^\epsilon \partial_x^{\ss - \ss' } \rho_R^\epsilon \partial_x^{\ss'} M_{ \Lambda_0 } \d A \d x \,.
	\end{aligned}
\end{equation}
By \eqref{Hmacro-eq}, one has
\begin{equation}\label{HETE-4}
	\begin{aligned}
		- \iint_{\R^3 \times \mathrm{SO(3)}} \partial_x^{\ss} \rho_R^\epsilon A e_1 \cdot  \nabla_x \partial_x^{\ss} f_R^\epsilon \d A \d x = & \int_{\R^3} \partial_x^{\ss} \rho_R^\epsilon ( \partial_t \partial_x^{\ss} \rho_R^\epsilon - \int_{\mathrm{SO}(3)} \partial_x^{\ss} R (f_1) \d A ) \d x \\
		= & \tfrac{1}{2} \tfrac{\d}{\d t} \| \partial_x^{\ss} \rho_R^\epsilon \|^2_{L^2_x} - \iint_{\R^3 \times \mathrm{SO(3)}} \partial_x^{\ss} \rho_R^\epsilon \partial_x^{\ss} R (f_1) \d A \d x \,.
	\end{aligned}
\end{equation}
Therefore, the equations \eqref{HETE-1}, \eqref{HETE-2}, \eqref{HETE-3} and \eqref{HETE-4} imply that
\begin{equation}\label{X1234}
	\begin{aligned}
		& \iint_{\R^3 \times \mathrm{SO(3)}} ( \partial_t + A e_1 \cdot  \nabla_x ) \partial_x^{\ss} \rho_R^\epsilon \cdot \partial_x^{\ss} \big( \tfrac{ f_R^\epsilon }{ M_{ \Lambda_0 } } - \rho_R^\epsilon \big) M_{ \Lambda_0 } \d A \d x \\ 
		= & \tfrac{1}{2} \tfrac{\d}{\d t} \| \partial_x^{\ss} \rho_R^\epsilon \|^2_{L^2_x} \underbrace{ - \iint_{\R^3 \times \mathrm{SO(3)}} \partial_x^{\ss} \rho_R^\epsilon \partial_x^{\ss} R (f_1) \d A \d x }_{:= \Xi_1} \\
		& \underbrace{ - \iint_{\R^3 \times \mathrm{SO(3)}} A e_1 \cdot  \nabla_x \partial_x^{\ss} \rho_R^\epsilon \partial_x^{\ss} \rho_R^\epsilon M_{ \Lambda_0 } \d A \d x }_{:= \Xi_2} \\
		& \underbrace{ - \sum_{0 \neq \ss' \leq \ss} C_{\ss}^{\ss'} \iint_{\R^3 \times \mathrm{SO(3)}} A e_1 \cdot  \nabla_x \partial_x^{\ss} \rho_R^\epsilon \partial_x^{\ss - \ss' } \rho_R^\epsilon \partial_x^{\ss'} M_{ \Lambda_0 } \d A \d x }_{:= \Xi_3} \\
		& \underbrace{ - \sum_{0 \neq \ss' \leq \ss} C_{\ss}^{\ss'} \iint_{\R^3 \times \mathrm{SO(3)}} ( \partial_t + A e_1 \cdot  \nabla_x ) \partial_x^{\ss} \rho_R^\epsilon \cdot \partial_x^{\ss - \ss'} ( \tfrac{ f_R^\epsilon }{ M_{ \Lambda_0 } } - \rho_R^\epsilon ) \partial_x^{\ss'} M_{ \Lambda_0 } \d A \d x }_{:= \Xi_4} \,.
	\end{aligned}
\end{equation}
For the quantity $\Xi_1$, the equality \eqref{HRf1-int} and the H\"older inequality implies
\begin{equation}\label{X1-bnd}
	\begin{aligned}
		| \Xi_1 | = & | \iint_{\R^3 \times \mathrm{SO(3)}} \partial_x^{\ss} \rho_R^\epsilon A e_1 \cdot \nabla_x \partial_x^{\ss} f_1 \d A \d x | \\
		\leq & \| \partial_x^{\ss} \rho_R^\epsilon \|_{L^2_x} \| \tfrac{ \nabla_x \partial_x^{\ss} f_1 }{ M_{ \Lambda_0 } } \|_{L^2_{x,A} ( M_{ \Lambda_0 } )} \leq \| \tfrac{ \nabla_x \partial_x^{\ss} f_1 }{ M_{ \Lambda_0 } } \|_{L^2_{x,A} ( M_{ \Lambda_0 } )} E_k^\frac{1}{2} (t) \,.
	\end{aligned}
\end{equation}
For the quantity $\Xi_2$, the integration by parts over $(x,A) \in \R^3 \times \mathrm{SO(3)}$ and \eqref{M-der} indicate that
\begin{equation}\label{X2-bnd}
	\begin{aligned}
		| \Xi_2 | = & | \tfrac{1}{2} \iint_{\R^3 \times \mathrm{SO(3)}} [ \partial_x^{\ss} \rho_R^\epsilon ]^2 A e_1 \cdot \nabla_x M_{ \Lambda_0 } \d A \d x | \\
		\lesssim & \| \nabla_x \Lambda_0 \|_{L^\infty_x} \| \partial_x^{\ss} \rho_R^\epsilon \|^2_{ L^2_x } \leq \| \nabla_x \Lambda_0 \|_{H^2_x} E_k (t) \,.
	\end{aligned}
\end{equation}
For the quantity $ \Xi_3 $, the integration by parts over $(x,A) \in \R^3 \times \mathrm{SO(3)}$ shows
\begin{equation}\label{X3-1}
	\begin{aligned}
		\Xi_3 = \sum_{0 \neq \ss' \leq \ss} C_{\ss}^{\ss'} \iint_{\R^3 \times \mathrm{SO(3)}} \partial_x^{\ss} \rho_R^\epsilon A e_1 \cdot \big( \nabla_x \partial_x^{\ss - \ss' } \rho_R^\epsilon \partial_x^{\ss'} M_{ \Lambda_0 } + \partial_x^{\ss - \ss' } \rho_R^\epsilon \nabla_x \partial_x^{\ss'} M_{ \Lambda_0 } \big) \d A \d x \,.
	\end{aligned}
\end{equation}
Following the similar arguments in Lemma 3.2 of \cite{JLT-M3AS-2019} or Lemma 2.2 of \cite{JLZ-ARMA-2020} and the Sobolev embedding $H^2_x \hookrightarrow L^\infty_x$, one has
\begin{equation}\label{X3-2}
	\begin{aligned}
		| \partial_x^{\ss} M_{ \Lambda_0 } | \lesssim \sum_{ \substack{ l |\tilde{\ss}| + k |\hat{\ss}| = |\ss| \\ |\tilde{\ss}|, |\hat{\ss}| \geq 1 , l, k \geq 0 \\ l^2 + k^2 \neq 0 } } | \partial_x^{\tilde{\ss}} \Lambda_0 |^l | \partial_x^{\hat{\ss}} \Lambda_0 |^k M_{ \Lambda_0 } \lesssim \| \nabla_x \Lambda_0 \|^{|\ss|}_{ H^{|\ss| + 1}_x } M_{ \Lambda_0 } 
	\end{aligned}
\end{equation}
for any multi-index $\ss \neq 0$. It then infers from \eqref{X3-1} and \eqref{X3-2} that 
\begin{equation}\label{X3-bnd}
	\begin{aligned}
		| \Xi_3 | \lesssim & \sum_{0 \neq \ss' \leq \ss} \| \partial_x^{\ss} \rho_R^\epsilon \|_{L^2_x} \big( \| \nabla_x \partial_x^{\ss - \ss'} \rho_R^\epsilon \|_{L^2_x} \| \nabla_x \Lambda_0 \|_{H^{|\ss'| + 1}_x}^{|\ss'|} + \| \partial_x^{\ss - \ss'} \rho_R^\epsilon \|_{L^2_x} \| \nabla_x \Lambda_0 \|_{H^{|\ss'| + 2}_x}^{|\ss'| + 1} \big) \\
		\lesssim & \| \nabla_x \Lambda_0 \|_{H^{s + 2}_x}^{s + 1} \| \partial_x^{\ss} \rho_R^\epsilon \|_{L^2_x} \sum_{ \ss' \leq \ss } \| \partial_x^{\ss'} \rho_R^\epsilon \|_{L^2_x} \\
		\lesssim & \| \nabla_x \Lambda_0 \|_{H^{s + 2}_x}^{s + 1} E_k^\frac{1}{2} (t) \sum_{0 \leq \mathfrak{j} \leq k} E_{ \mathfrak{j} }^\frac{1}{2} (t) \,.
	\end{aligned}
\end{equation}
For the quantity $\Xi_4$, it is easy to see
\begin{equation}\label{X4-1}
	\begin{aligned}
		\Xi_4 = & \underbrace{ - \sum_{0 \neq \ss' \leq \ss} C_{\ss}^{\ss'} \iint_{\R^3 \times \mathrm{SO(3)}} \partial_t \partial_x^{\ss} \rho_R^\epsilon \cdot \partial_x^{\ss - \ss'} ( \tfrac{ f_R^\epsilon }{ M_{ \Lambda_0 } } - \rho_R^\epsilon ) \partial_x^{\ss'} M_{ \Lambda_0 } \d A \d x }_{:= \Xi_{41}} \\
		& + \underbrace{ \sum_{0 \neq \ss' \leq \ss} C_{\ss}^{\ss'} \iint_{\R^3 \times \mathrm{SO(3)}}  \partial_x^{\ss} \rho_R^\epsilon \cdot ( A e_1 \cdot  \nabla_x ) \big[ \partial_x^{\ss - \ss'} ( \tfrac{ f_R^\epsilon }{ M_{ \Lambda_0 } } - \rho_R^\epsilon ) \partial_x^{\ss'} M_{ \Lambda_0 } \big] \d A \d x }_{:= \Xi_{42} } \,.
	\end{aligned}
\end{equation}
Recalling the equation \eqref{Hmacro-eq} and \eqref{HRf1-int}, it follows
\begin{equation*}
	\begin{aligned}
		\Xi_{41} = & \sum_{0 \neq \ss' \leq \ss} C_{\ss}^{\ss'} \iint_{\R^3 \times \mathrm{SO(3)}} \Big( \nabla_x \cdot \int_{\mathrm{SO}(3)} A e_1 \partial_x^{\ss} ( f_R^\epsilon - f_1 ) \d A \Big) \\
		& \qquad \qquad \qquad \qquad \qquad \cdot \nabla_x \partial_x^{\ss - \ss'} ( \tfrac{ f_R^\epsilon }{ M_{ \Lambda_0 } } - \rho_R^\epsilon ) \partial_x^{\ss'} M_{ \Lambda_0 } \d A \d x \\
		= & - \sum_{0 \neq \ss' \leq \ss} C_{\ss}^{\ss'} \iint_{\R^3 \times \mathrm{SO(3)}} \Big( \int_{\mathrm{SO}(3)} A e_1 \partial_x^{\ss} ( f_R^\epsilon - f_1 ) \d A \Big) \\
		& \qquad \qquad \qquad \qquad \qquad \cdot \nabla_x \big[ \partial_x^{\ss - \ss'} ( \tfrac{ f_R^\epsilon }{ M_{ \Lambda_0 } } - \rho_R^\epsilon ) \partial_x^{\ss'} M_{ \Lambda_0 } \big] \d A \d x \,,
	\end{aligned}
\end{equation*}
which, together with \eqref{X3-2}, means that
\begin{equation*}
	\begin{aligned}
		| \Xi_{41} | \lesssim & \sum_{0 \neq \ss' \leq \ss} \| \int_{\mathrm{SO}(3)} A e_1 \partial_x^{\ss} ( f_R^\epsilon - f_1 ) \d A \|_{L^2_x} \\
		& \qquad \times \big( \| \nabla_x \partial_x^{\ss - \ss'} ( \tfrac{ f_R^\epsilon }{ M_{ \Lambda_0 } } - \rho_R^\epsilon ) \|_{L^2_{x, A} ( M_{ \Lambda_0 } )} \| \nabla_x \Lambda_0 \|_{H^{|\ss'| + 1}_x}^{|\ss'|} \\
		& \qquad \qquad \qquad + \| \partial_x^{\ss - \ss'} ( \tfrac{ f_R^\epsilon }{ M_{ \Lambda_0 } } - \rho_R^\epsilon ) \|_{L^2_{x, A} ( M_{ \Lambda_0 } )} \| \nabla_x \Lambda_0 \|_{H^{|\ss'| + 2}_x}^{|\ss'| + 1} \big) \\
		\lesssim & \| \nabla_x \Lambda_0 \|_{H^{s + 2}_x}^{s + 1} \| \int_{\mathrm{SO}(3)} A e_1 \partial_x^{\ss} ( f_R^\epsilon - f_1 ) \d A \|_{L^2_x} \sum_{ \ss' \leq \ss } \| \partial_x^{\ss'} ( \tfrac{ f_R^\epsilon }{ M_{ \Lambda_0 } } - \rho_R^\epsilon ) \|_{L^2_{x, A} ( M_{ \Lambda_0 } )} \\
		\lesssim & \| \nabla_x \Lambda_0 \|_{H^{s + 2}_x}^{s + 1} \| \int_{\mathrm{SO}(3)} A e_1 \partial_x^{\ss} ( f_R^\epsilon - f_1 ) \d A \|_{L^2_x} \sum_{0 \leq \mathfrak{j} \leq k} E_{\mathfrak{j}}^\frac{1}{2} (t) \,.
	\end{aligned}
\end{equation*}
Observe that
\begin{equation*}
	\begin{aligned}
		\int_{\mathrm{SO}(3)} A e_1 \partial_x^{\ss} ( f_R^\epsilon - f_1 ) \d A = \sum_{\ss' \leq \ss} C_{\ss}^{\ss'} \int_{\mathrm{SO}(3)} A e_1 \partial_x^{\ss - \ss'} \big[ ( \tfrac{ f_R^\epsilon }{ M_{ \Lambda_0 } } - \rho_R^\epsilon ) + \rho_R^\epsilon - \tfrac{f_1}{ M_{ \Lambda_0 } } \big] \partial_x^{\ss'} M_{ \Lambda_0 } \d A \,.
	\end{aligned}
\end{equation*}
Together with \eqref{X3-2}, it is easy to obtain
\begin{equation*}
	\begin{aligned}
		& \| \int_{\mathrm{SO}(3)} A e_1 \partial_x^{\ss} ( f_R^\epsilon - f_1 ) \d A \|_{L^2_x} \\
		\lesssim & ( 1 + \| \nabla_x \Lambda_0 \|_{ H^{s+1}_x }^s ) \sum_{\ss' \leq \ss} \big( \| \partial_x^{\ss'} ( \tfrac{ f_R^\epsilon }{ M_{ \Lambda_0 } } - \rho_R^\epsilon ) \|_{L^2_{x,A} ( M_{ \Lambda_0 } ) } + \| \partial_x^{\ss'} \rho_R^\epsilon \|_{L^2_x} + \| \partial_x^{\ss'} ( \tfrac{ f_1 }{ M_{ \Lambda_0 } } ) \|_{L^2_{x,A} ( M_{ \Lambda_0 } ) } \big) \\
		\lesssim & ( 1 + \| \nabla_x \Lambda_0 \|_{ H^{s+1}_x }^s ) \big( \sum_{0 \leq \mathfrak{j} \leq k} E_{\mathfrak{j}}^\frac{1}{2} (t) + \| \tfrac{ f_1 }{ M_{ \Lambda_0 } } \|_{H^s_x L^2_A ( M_{ \Lambda_0 } )} \big) \,.
	\end{aligned}
\end{equation*}
As a result,
\begin{equation}\label{X4-2}
	\begin{aligned}
		| \Xi_{41} | \lesssim ( 1 + \| \nabla_x \Lambda_0 \|_{ H^{s+2}_x }^{s+1} ) \sum_{0 \leq \mathfrak{j} \leq k} \big( E_{\mathfrak{j}} (t) + \| \tfrac{ f_1 }{ M_{ \Lambda_0 } } \|_{H^s_x L^2_A ( M_{ \Lambda_0 } )} E_{\mathfrak{j}}^\frac{1}{2} (t) \big) \,.
	\end{aligned}
\end{equation}
Following the similar arguments in \eqref{X3-bnd}, the term $\Xi_{42}$ can be bounded by
\begin{equation}\label{X4-3}
	\begin{aligned}
		| \Xi_{42} | \lesssim & \| \nabla_x \Lambda_0 \|_{H^{s + 2}_x}^{s + 1} \| \partial_x^{\ss} ( \tfrac{ f_R^\epsilon }{ M_{ \Lambda_0 } } - \rho_R^\epsilon ) \|_{L^2_{x,A} ( M_{ \Lambda_0 } ) } \sum_{ \ss' \leq \ss } \| \partial_x^{\ss'} ( \tfrac{ f_R^\epsilon }{ M_{ \Lambda_0 } } - \rho_R^\epsilon ) \|_{L^2_{x,A} ( M_{ \Lambda_0 } ) } \\
		\lesssim & \| \nabla_x \Lambda_0 \|_{H^{s + 2}_x}^{s + 1} E_k^\frac{1}{2} (t) \sum_{0 \leq \mathfrak{j} \leq k} E_{ \mathfrak{j} }^\frac{1}{2} (t)
	\end{aligned}
\end{equation}
It thereby follows from \eqref{X4-1}, \eqref{X4-2} and \eqref{X4-3} that
\begin{equation}\label{X4-bnd}
	\begin{aligned}
		| \Xi_4 | \lesssim ( 1 + \| \nabla_x \Lambda_0 \|_{ H^{s+2}_x }^{s+1} ) \sum_{0 \leq \mathfrak{j} \leq k} \big( E_{\mathfrak{j}} (t) + \| \tfrac{ f_1 }{ M_{ \Lambda_0 } } \|_{H^s_x L^2_A ( M_{ \Lambda_0 } )} E_{\mathfrak{j}}^\frac{1}{2} (t) \big) \,.
	\end{aligned}
\end{equation}

Note that from the similar arguments in \eqref{X3-bnd}, 
\begin{equation*}
	\begin{aligned}
		\| \tfrac{ \nabla_x \partial_x^{\ss} f_1 }{ M_{ \Lambda_0 } } \|_{L^2_{x,A} ( M_{ \Lambda_0 } )} \lesssim \| \nabla_x \Lambda_0 \|_{ H^{s+2}_x }^{s+1} \| \tfrac{ f_1 }{ M_{ \Lambda_0 } } \|_{H^s_x L^2_A ( M_{ \Lambda_0 } )} \,.
	\end{aligned}
\end{equation*}
Substituting \eqref{X1-bnd}, \eqref{X2-bnd}, \eqref{X3-bnd} and \eqref{X4-bnd} into \eqref{X1234}, one therefore has
\begin{equation}\label{HETEMa-bnd}
	\begin{aligned}
		& \iint_{\R^3 \times \mathrm{SO(3)}} ( \partial_t + A e_1 \cdot  \nabla_x ) \partial_x^{\ss} \rho_R^\epsilon \cdot \partial_x^{\ss} \big( \tfrac{ f_R^\epsilon }{ M_{ \Lambda_0 } } - \rho_R^\epsilon \big) M_{ \Lambda_0 } \d A \d x \\
		\geq & \tfrac{1}{2} \tfrac{\d}{\d t} \| \partial_x^{\ss} \rho_R^\epsilon \|^2_{L^2_x} - C ( 1 + \| \nabla_x \Lambda_0 \|_{ H^{s+2}_x }^{s+1} ) \sum_{0 \leq \mathfrak{j} \leq k} \big( E_{\mathfrak{j}} (t) + \| \tfrac{ f_1 }{ M_{ \Lambda_0 } } \|_{H^s_x L^2_A ( M_{ \Lambda_0 } )} E_{\mathfrak{j}}^\frac{1}{2} (t) \big) \,.
	\end{aligned}
\end{equation}

{\bf Step 3. Estimates of coercivity: $ - \frac{1}{\epsilon} \partial_x^{\ss} \big[ \tfrac{1}{ M_{ \Lambda_0 } } \L_{ M_{ \Lambda_0 } } f_R^\epsilon \big] $.} Obviously,
	\begin{align}\label{C1C2C3}
		\no & \iint_{\R^3 \times \mathrm{SO(3)}} - \frac{1}{\epsilon} \partial_x^{\ss} \big[ \tfrac{1}{ M_{ \Lambda_0 } } \L_{ M_{ \Lambda_0 } } f_R^\epsilon \big] \cdot \partial_x^{\ss} \big( \tfrac{ f_R^\epsilon }{ M_{ \Lambda_0 } } - \rho_R^\epsilon \big) M_{ \Lambda_0 } \d A \d x \\
		\no = & \underbrace{ \iint_{\R^3 \times \mathrm{SO(3)}} - \frac{1}{\epsilon} \nabla_A \cdot \big[ M_{ \Lambda_0 } \nabla_A \partial_x^{\ss} ( \tfrac{ f_R^\epsilon }{ M_{ \Lambda_0 } } ) \big] \cdot \partial_x^{\ss} \big( \tfrac{ f_R^\epsilon }{ M_{ \Lambda_0 } } - \rho_R^\epsilon \big) \d A \d x }_{:= \C_1} \\
		\no & \underbrace{ - \frac{1}{\epsilon} \sum_{0 \neq \ss' \leq \ss} d C_{\ss}^{\ss'} \iint_{\R^3 \times \mathrm{SO(3)}} \nabla_A \cdot [ \partial_x^{\ss'} M_{ \Lambda_0 } \nabla_A \partial_x^{\ss - \ss'} ( \tfrac{ f_R^\epsilon }{ M_{ \Lambda_0 } } ) \cdot \partial_x^{\ss} \big( \tfrac{ f_R^\epsilon }{ M_{ \Lambda_0 } } - \rho_R^\epsilon \big) \d A \d x }_{:= \C_2} \\
		& \left.
		\begin{aligned}
		& - \frac{1}{\epsilon} \sum_{0 \neq \ss' \leq \ss} \sum_{\ss'' \leq \ss - \ss'} d C_{\ss}^{\ss'} C_{\ss - \ss'}^{\ss''} M_{ \Lambda_0 } \partial_x^{\ss'} ( \tfrac{1}{ M_{ \Lambda_0 } } ) \\
		& \qquad \qquad \qquad \times \nabla_A \cdot [ \partial_x^{\ss''} M_{ \Lambda_0 } \nabla_A \partial_x^{\ss - \ss' - \ss''} ( \tfrac{ f_R^\epsilon }{ M_{ \Lambda_0 } } ) \cdot \partial_x^{\ss} \big( \tfrac{ f_R^\epsilon }{ M_{ \Lambda_0 } } - \rho_R^\epsilon \big) \d A \d x 
		\end{aligned}
	    \right\} {:= \C_3} \,.
	\end{align}
Lemma \ref{CE} shows 
\begin{equation}\label{C1-bnd}
	\begin{aligned}
		\C_1 = \tfrac{d}{\epsilon} \| \nabla_A \partial_x^{\ss} ( \tfrac{ f_R^\epsilon }{ M_{ \Lambda_0 } } ) \|^2_{L^2_{x,A} ( M_{ \Lambda_0 } ) } \,.
	\end{aligned}
\end{equation}
By integration by parts over $( x, A ) \in \R^3 \times \mathrm{SO(3)}$, it infers
\begin{equation*}
	\begin{aligned}
		| \C_2 | \lesssim \frac{1}{\epsilon} \sum_{0 \neq \ss' \leq \ss} \iint_{\R^3 \times \mathrm{SO(3)}} | \partial_x^{\ss'} M_{ \Lambda_0 } | \cdot | \nabla_A \partial_x^{\ss - \ss'} ( \tfrac{ f_R^\epsilon }{ M_{ \Lambda_0 } } ) | \cdot | \nabla_A  \partial_x^{\ss} \big( \tfrac{ f_R^\epsilon }{ M_{ \Lambda_0 } } \big) | \d A \d x  \,.
	\end{aligned}
\end{equation*}
Together with \eqref{X3-2},
\begin{equation}\label{C2-bnd}
	\begin{aligned}
		| \C_2 | \lesssim & \frac{1}{\epsilon} \| \nabla_A  \partial_x^{\ss} \big( \tfrac{ f_R^\epsilon }{ M_{ \Lambda_0 } } \big) \|_{ L^2_{x,A} ( M_{ \Lambda_0 } ) } \sum_{0 \neq \ss' \leq \ss} \| \nabla_x \Lambda_0 \|^{|\ss'|}_{ H^{|\ss'| + 1}_x } \| \nabla_A \partial_x^{\ss - \ss'} ( \tfrac{ f_R^\epsilon }{ M_{ \Lambda_0 } } ) \|_{ L^2_{x,A} ( M_{ \Lambda_0 } ) } \\
		\lesssim & \frac{1}{\epsilon} \| \nabla_x \Lambda_0 \|^s_{ H^{s + 1}_x } \| \nabla_A  \partial_x^{\ss} \big( \tfrac{ f_R^\epsilon }{ M_{ \Lambda_0 } } \big) \|_{ L^2_{x,A} ( M_{ \Lambda_0 } ) } \sum_{0 \neq \ss' \leq \ss} \| \nabla_A \partial_x^{\ss - \ss'} ( \tfrac{ f_R^\epsilon }{ M_{ \Lambda_0 } } ) \|_{ L^2_{x,A} ( M_{ \Lambda_0 } ) } \\
		\lesssim & \| \nabla_x \Lambda_0 \|^s_{ H^{s + 1}_x } \tfrac{1}{\epsilon} \| \nabla_A \partial_x^{\ss} ( \tfrac{ f_R^\epsilon }{ M_{ \Lambda_0 } } ) \|_{L^2_{x,A} ( M_{ \Lambda_0 } ) } \sum_{0 \leq \mathfrak{j} \leq k - 1} D_{\mathfrak{j}}^\frac{1}{2} (t) \,.
	\end{aligned}
\end{equation}
By integration by parts over $(x,A) \in \R^3 \times \mathrm{SO(3)}$, the quantity $\C_3$ can be bounded by
\begin{equation*}
	\begin{aligned}
		| \C_3 | \lesssim & \frac{1}{\epsilon} \sum_{0 \neq \ss' \leq \ss} \sum_{\ss'' \leq \ss - \ss'} \Big\{ \iint_{\R^3 \times \mathrm{SO(3)}} | \nabla_A \partial_x^{\ss} ( \tfrac{ f_R^\epsilon }{ M_{ \Lambda_0 } } ) | \, | \nabla_A \partial_x^{ \ss - \ss' - \ss'' } ( \tfrac{ f_R^\epsilon }{ M_{ \Lambda_0 } } ) | \\
		& \qquad \qquad \qquad \qquad \qquad \qquad \times | \partial_x^{\ss'} ( \tfrac{1}{ M_{ \Lambda_0 } } ) | \, | \partial_x^{\ss''} M_{ \Lambda_0 } | M_{ \Lambda_0 } \d A \d x \\
		& + \iint_{\R^3 \times \mathrm{SO(3)}} | \nabla_A \partial_x^{ \ss - \ss' - \ss'' } ( \tfrac{ f_R^\epsilon }{ M_{ \Lambda_0 } } ) | \, | \partial_x^{\ss} ( \tfrac{ f_R^\epsilon }{ M_{ \Lambda_0 } } - \rho_R^\epsilon ) | \\
		& \qquad \qquad \qquad \qquad \quad \times | \nabla_A ( M_{ \Lambda_0 } \partial_x^{\ss'} ( \tfrac{1}{ M_{ \Lambda_0 } } ) ) | \, | \partial_x^{\ss''} M_{ \Lambda_0 } | \d A \d x \Big\} \,.
	\end{aligned}
\end{equation*}
Following the similar arguments in \eqref{X3-2}, one can derive that
\begin{equation}\label{C3-0}
	\begin{aligned}
		& | \partial_x^{\ss''} M_{ \Lambda_0 } | \lesssim ( 1 + \| \nabla_x \Lambda_0 \|^{\ss''}_{ H_x^{ |\ss''| + 1 } } ) M_{ \Lambda_0 } \,, \\
		& | \partial_x^{\ss'} ( \tfrac{1}{ M_{ \Lambda_0 } } ) | \lesssim \| \nabla_x \Lambda_0 \|^{|\ss'|}_{ H_x^{ |\ss'| + 1 } } \tfrac{1}{ M_{ \Lambda_0 } } \,, \\
		& | \nabla_A ( M_{ \Lambda_0 } \partial_x^{\ss'} ( \tfrac{1}{ M_{ \Lambda_0 } } ) ) | \lesssim \| \nabla_x \Lambda_0 \|^{|\ss'|}_{ H_x^{ |\ss'| + 1 } } \,. 
	\end{aligned}
\end{equation}
Then,
\begin{equation}\label{C3-1}
	\begin{aligned}
		| \C_3 | \lesssim & \frac{1}{\epsilon} \sum_{ \substack{ 0 \neq \ss' \leq \ss \\ \ss'' \leq \ss - \ss' } } ( 1 + \| \nabla_x \Lambda_0 \|^{\ss''}_{ H_x^{ |\ss''| + 1 } } ) \| \nabla_x \Lambda_0 \|^{|\ss'|}_{ H_x^{ |\ss'| + 1 } } \| \nabla_A \partial_x^{\ss - \ss' - \ss''} ( \tfrac{ f_R^\epsilon }{ M_{ \Lambda_0 } } ) \|_{ L^2_{x,A} ( M_{ \Lambda_0 } ) } \\
		& \qquad \qquad \times \big( \| \nabla_A \partial_x^{\ss} ( \tfrac{ f_R^\epsilon }{ M_{ \Lambda_0 } } ) \|_{ L^2_{x,A} ( M_{ \Lambda_0 } ) } + \| \partial_x^{\ss} ( \tfrac{ f_R^\epsilon }{ M_{ \Lambda_0 } } - \rho_R^\epsilon ) \|_{ L^2_{x,A} ( M_{ \Lambda_0 } ) } \big) \,.
	\end{aligned}
\end{equation}
Here we want to dominate the quantity $ \| \partial_x^{\ss} ( \tfrac{ f_R^\epsilon }{ M_{ \Lambda_0 } } - \rho_R^\epsilon ) \|_{ L^2_{x,A} ( M_{ \Lambda_0 } ) } $ in terms of the norm $ \| \nabla_A \partial_x^{\ss} ( \tfrac{ f_R^\epsilon }{ M_{ \Lambda_0 } } ) \|_{ L^2_{x,A} ( M_{ \Lambda_0 } ) } $. However the Poincar\'e inequality given in Lemma \ref{CE} cannot be applied here for $\ss \neq 0$ due to $ \int_{\mathrm{SO}(3)} \partial_x^{\ss} ( \tfrac{ f_R^\epsilon }{ M_{ \Lambda_0 } } - \rho_R^\epsilon ) M_{\Lambda_0} \d A \neq 0 $. Fortunately, one has
	\begin{align*}
		\lambda_0 \| \partial_x^{\ss} ( \tfrac{ f_R^\epsilon }{ M_{ \Lambda_0 } } - \rho_R^\epsilon ) \|_{ L^2_{x,A} ( M_{ \Lambda_0 } ) }^2 \leq & \| \nabla_A \partial_x^{\ss} ( \tfrac{ f_R^\epsilon }{ M_{ \Lambda_0 } } ) \|_{ L^2_{x,A} ( M_{ \Lambda_0 } ) }^2 \\
		& + C | \int_{\mathrm{SO}(3)} \partial_x^{\ss} ( \tfrac{ f_R^\epsilon }{ M_{ \Lambda_0 } } - \rho_R^\epsilon ) M_{\Lambda_0} \d A |^2 
	\end{align*}
for any $ \ss \neq 0 $. Thanks to $ \int_{\mathrm{SO}(3)} \partial_x^{\ss} [ ( \tfrac{ f_R^\epsilon }{ M_{ \Lambda_0 } } - \rho_R^\epsilon ) M_{\Lambda_0} ] \d A = 0 $, it holds
\begin{equation*}
	\begin{aligned}
		\int_{\mathrm{SO}(3)} \partial_x^{\ss} ( \tfrac{ f_R^\epsilon }{ M_{ \Lambda_0 } } - \rho_R^\epsilon ) M_{\Lambda_0} \d A = - \sum_{0 \neq \ss' \leq \ss} C_{\ss}^{\ss'} \int_{\mathrm{SO}(3)} \partial_x^{\ss - \ss'} ( \tfrac{ f_R^\epsilon }{ M_{ \Lambda_0 } } - \rho_R^\epsilon ) \partial_x^{\ss'} M_{\Lambda_0} \d A \,.
	\end{aligned}
\end{equation*}
It then follow from \eqref{X3-2} that
\begin{equation*}
	\begin{aligned}
		| \int_{\mathrm{SO}(3)} \partial_x^{\ss} ( \tfrac{ f_R^\epsilon }{ M_{ \Lambda_0 } } - \rho_R^\epsilon ) M_{\Lambda_0} \d A |^2 \lesssim \| \nabla_x \Lambda_0 \|^{2s}_{H^{s+1}_x} \sum_{0 \neq \ss' \leq \ss} \| \partial_x^{\ss - \ss'} ( \tfrac{ f_R^\epsilon }{ M_{ \Lambda_0 } } - \rho_R^\epsilon ) \|^2_{ L^2_A ( M_{ \Lambda_0 } ) } \,,
	\end{aligned}
\end{equation*}
which means that
\begin{equation*}
	\begin{aligned}
		\lambda_0 \| \partial_x^{\ss} ( \tfrac{ f_R^\epsilon }{ M_{ \Lambda_0 } } - \rho_R^\epsilon ) \|_{ L^2_{x,A} ( M_{ \Lambda_0 } ) }^2 \leq & \| \nabla_A \partial_x^{\ss} ( \tfrac{ f_R^\epsilon }{ M_{ \Lambda_0 } } ) \|_{ L^2_{x,A} ( M_{ \Lambda_0 } ) }^2 \\
		& + C \| \nabla_x \Lambda_0 \|^{2s}_{H^{s+1}_x} \sum_{0 \neq \ss' \leq \ss} \| \partial_x^{\ss - \ss'} ( \tfrac{ f_R^\epsilon }{ M_{ \Lambda_0 } } - \rho_R^\epsilon ) \|^2_{ L^2_{x,A} ( M_{ \Lambda_0 } ) } 
	\end{aligned}
\end{equation*}
for any $\ss \in \mathbb{N}^3$ with $1 \leq |\ss| = k \leq s$. Repeating the above procedure $k$ times and combining with Lemma \ref{CE}, one can derived that
\begin{equation}\label{HPoincare}
	\begin{aligned}
		\lambda_0 \| \partial_x^{\ss} ( \tfrac{ f_R^\epsilon }{ M_{ \Lambda_0 } } - \rho_R^\epsilon ) \|_{ L^2_{x,A} ( M_{ \Lambda_0 } ) }^2 \leq & \| \nabla_A \partial_x^{\ss} ( \tfrac{ f_R^\epsilon }{ M_{ \Lambda_0 } } ) \|_{ L^2_{x,A} ( M_{ \Lambda_0 } ) }^2 \\
		& + C \| \nabla_x \Lambda_0 \|^{2s^2}_{H^{s+1}_x} \sum_{0 \neq \ss' \leq \ss} \| \nabla_A \partial_x^{\ss - \ss'} ( \tfrac{ f_R^\epsilon }{ M_{ \Lambda_0 } } ) \|^2_{ L^2_{x,A} ( M_{ \Lambda_0 } ) } \\
		\leq & \| \nabla_A \partial_x^{\ss} ( \tfrac{ f_R^\epsilon }{ M_{ \Lambda_0 } } ) \|_{ L^2_{x,A} ( M_{ \Lambda_0 } ) }^2 + C \epsilon \| \nabla_x \Lambda_0 \|^{2s^2}_{H^{s+1}_x} \sum_{0 \leq \mathfrak{j} \leq k - 1} D_{\mathfrak{j}} (t) 
	\end{aligned}
\end{equation}
for any $\ss \in \mathbb{N}^3$ with $1 \leq |\ss| = k \leq s$. From plugging \eqref{HPoincare} into \eqref{C3-1}, it then follows
\begin{equation}\label{C3-bnd}
	\begin{aligned}
		| \C_3 | \lesssim & \frac{1}{\epsilon} ( 1 + \| \nabla_x \Lambda_0 \|^s_{ H_x^{ s + 1 } } ) \| \nabla_x \Lambda_0 \|^s_{ H_x^{ s + 1 } } \sum_{0 \neq \ss' \leq \ss} \| \nabla_A \partial_x^{\ss - \ss'} ( \tfrac{ f_R^\epsilon }{ M_{ \Lambda_0 } } ) \|_{ L^2_{x,A} ( M_{ \Lambda_0 } ) } \\
		& \times \Big( \| \nabla_A \partial_x^{\ss} ( \tfrac{ f_R^\epsilon }{ M_{ \Lambda_0 } } ) \|_{ L^2_{x,A} ( M_{ \Lambda_0 } ) } + \| \nabla_x \Lambda_0 \|^{2s^2}_{H^{s+1}_x} \sum_{0 \neq \ss' \leq \ss} \| \nabla_A \partial_x^{\ss - \ss'} ( \tfrac{ f_R^\epsilon }{ M_{ \Lambda_0 } } ) \|^2_{ L^2_{x,A} ( M_{ \Lambda_0 } ) } \Big) \\
		\lesssim & ( 1 + \| \nabla_x \Lambda_0 \|^s_{ H_x^{ s + 1 } } ) \| \nabla_x \Lambda_0 \|^s_{ H_x^{ s + 1 } } \sum_{0 \leq \mathfrak{j} \leq k - 1} D_{\mathfrak{j}}^\frac{1}{2} (t) \\
		& \qquad \qquad \qquad \times  \big( \tfrac{1}{\sqrt{\epsilon}} \| \nabla_A \partial_x^{\ss} ( \tfrac{ f_R^\epsilon }{ M_{ \Lambda_0 } } ) \|_{L^2_{x,A} ( M_{ \Lambda_0 } ) } + \| \nabla_x \Lambda_0 \|^{2s^2}_{H^{s+1}_x} \sum_{0 \leq \mathfrak{j} \leq k - 1} D_{\mathfrak{j}}^\frac{1}{2} (t) \big) \,.
	\end{aligned}
\end{equation}
It is therefore deduced from plugging the bounds \eqref{C1-bnd}, \eqref{C2-bnd} and \eqref{C3-bnd} into \eqref{C1C2C3} that
\begin{equation}\label{HEC-bnd}
	\begin{aligned}
		& \iint_{\R^3 \times \mathrm{SO(3)}} - \frac{1}{\epsilon} \partial_x^{\ss} \big[ \tfrac{1}{ M_{ \Lambda_0 } } \L_{ M_{ \Lambda_0 } } f_R^\epsilon \big] \cdot \partial_x^{\ss} \big( \tfrac{ f_R^\epsilon }{ M_{ \Lambda_0 } } - \rho_R^\epsilon \big) M_{ \Lambda_0 } \d A \d x \\ 
		\geq & \tfrac{d}{\epsilon} \| \nabla_A \partial_x^{\ss} ( \tfrac{ f_R^\epsilon }{ M_{ \Lambda_0 } } ) \|^2_{L^2_{x,A} ( M_{ \Lambda_0 } ) } - C ( 1 + \| \nabla_x \Lambda_0 \|^s_{ H_x^{ s + 1 } } ) \| \nabla_x \Lambda_0 \|^{s^2 + s}_{ H_x^{ s + 1 } } \sum_{0 \leq \mathfrak{j} \leq k - 1} D_{\mathfrak{j}} (t) \\
		& - C ( 1 + \| \nabla_x \Lambda_0 \|^s_{ H_x^{ s + 1 } } ) \| \nabla_x \Lambda_0 \|^s_{ H_x^{ s + 1 } } \tfrac{1}{\sqrt{\epsilon}} \| \nabla_A \partial_x^{\ss} ( \tfrac{ f_R^\epsilon }{ M_{ \Lambda_0 } } ) \|_{L^2_{x,A} ( M_{ \Lambda_0 } ) } \sum_{0 \leq \mathfrak{j} \leq k - 1} D_{\mathfrak{j}}^\frac{1}{2} (t) \,.
	\end{aligned}
\end{equation}

{\bf Step 4. Estimates of error linear operator: $ \frac{1}{\epsilon} \partial_x^{\ss} \big[ \tfrac{1}{ M_{ \Lambda_0 } } L_R ( f_R^\epsilon - \rho_R^\epsilon M_{ \Lambda_0 } ) \big] $.} Recalling the definition of $ L_R ( f_R^\epsilon - \rho_R^\epsilon M_{ \Lambda_0 } ) $ in \eqref{LR-opt} and integration by parts over $A \in \mathrm{SO(3)}$, it follows
	\begin{align}\label{F1F2F3}
		\no & - \iint_{\R^3 \times \mathrm{SO(3)}} \frac{1}{\epsilon} \partial_x^{\ss} \big[ \tfrac{1}{ M_{ \Lambda_0 } } L_R ( f_R^\epsilon - \rho_R^\epsilon M_{ \Lambda_0 } ) \big] \cdot \partial_x^{\ss} \big( \tfrac{ f_R^\epsilon }{ M_{ \Lambda_0 } } - \rho_R^\epsilon \big) M_{ \Lambda_0 } \d A \d x \\
		\no = & - \frac{\nu_0}{c_1} \frac{1}{\epsilon} \iint_{\R^3 \times \mathrm{SO(3)}} \partial_x^{\ss} ( \tfrac{ f_R^\epsilon }{ M_{ \Lambda_0 } } - \rho_R^\epsilon ) M_{ \Lambda_0 } \\
		\no & \qquad \qquad \qquad \times \partial_x^{\ss} \Big\{ \tfrac{1}{ M_{ \Lambda_0 } } \nabla_A \cdot \big[ M_{ \Lambda_0 } \nabla_A ( A \cdot P_{ T_{ \Lambda_0 } } ( \lambda [ f_R^\epsilon - \rho_R^\epsilon M_{ \Lambda_0 } ] ) ) \big] \Big\} \d A \d x \\
		\no = & \underbrace{ - \frac{\nu_0}{c_1} \frac{1}{\epsilon} \iint_{\R^3 \times \mathrm{SO(3)}} \nabla_A \partial_x^{\ss} ( \tfrac{ f_R^\epsilon }{ M_{ \Lambda_0 } } ) \nabla_A ( A \cdot \partial_x^{\ss} P_{ T_{ \Lambda_0 } } ( \lambda [ f_R^\epsilon - \rho_R^\epsilon M_{ \Lambda_0 } ] ) ) M_{ \Lambda_0 } \d A \d x }_{:= F_1} \\
		\no & \left.
		\begin{aligned}
		& - \frac{\nu_0}{c_1} \frac{1}{\epsilon} \sum_{0 \neq \ss' \leq \ss} C_{\ss}^{\ss'} \iint_{\R^3 \times \mathrm{SO(3)}} \nabla_A \partial_x^{\ss} ( \tfrac{ f_R^\epsilon }{ M_{ \Lambda_0 } } ) \\
		& \qquad \qquad \qquad \qquad \times \big[ \partial_x^{\ss'} M_{ \Lambda_0 } \nabla_A ( A \cdot \partial_x^{\ss - \ss'} P_{ T_{ \Lambda_0 } } ( \lambda [ f_R^\epsilon - \rho_R^\epsilon M_{ \Lambda_0 } ] ) ) \big] \d A \d x
		\end{aligned} 
	    \right\} {:= F_2}  \\
	    & \left.
	    \begin{aligned}
		& - \frac{\nu_0}{c_1} \frac{1}{\epsilon} \sum_{0 \neq \ss' \leq \ss} \sum_{\ss'' \leq \ss - \ss'} C_{\ss}^{\ss'} C_{\ss - \ss'}^{\ss''} \iint_{\R^3 \times \mathrm{SO(3)}} \nabla_A \big[ \partial_x^{\ss'} ( \tfrac{1}{ M_{ \Lambda_0 } } ) M_{ \Lambda_0 } \partial_x^{\ss} ( \tfrac{ f_R^\epsilon }{ M_{ \Lambda_0 } } - \rho_R^\epsilon ) \big] \\
		& \qquad \qquad \times \big[ \partial_x^{\ss''} M_{ \Lambda_0 } \nabla_A ( A \cdot \partial_x^{\ss - \ss' - \ss''} P_{ T_{ \Lambda_0 } } ( \lambda [ f_R^\epsilon - \rho_R^\epsilon M_{ \Lambda_0 } ] ) ) \big] \d A \d x 
		\end{aligned}
	    \right\} {:= F_3} \,.
	\end{align}

\underline{\em Case 4.1. Control of $F_1$.} By \eqref{EE-2}, there holds
\begin{equation}\label{F1-1}
	\begin{aligned}
		| \nabla_A ( A \cdot \partial_x^{\ss} P_{ T_{ \Lambda_0 } } ( \lambda [ f_R^\epsilon - \rho_R^\epsilon M_{ \Lambda_0 } ] ) ) | \leq 5 | \partial_x^{\ss} P_{ T_{ \Lambda_0 } } ( \lambda [ f_R^\epsilon - \rho_R^\epsilon M_{ \Lambda_0 } ] ) | \,.
	\end{aligned}
\end{equation}
Lemma \ref{Lmm-SOHB-spt} indicates that
\begin{equation*}
	\begin{aligned}
		\partial_x^{\ss} P_{ T_{ \Lambda_0 } } ( \lambda [ f_R^\epsilon - \rho_R^\epsilon M_{ \Lambda_0 } ] ) = & \tfrac{1}{2} \big\{ \partial_x^{\ss} \lambda [ f_R^\epsilon - \rho_R^\epsilon M_{ \Lambda_0 } ]  - \Lambda_0 \partial_x^{\ss} \lambda [ f_R^\epsilon - \rho_R^\epsilon M_{ \Lambda_0 } ]^\top \Lambda_0 \big\} \\
		& + \tfrac{1}{2} \sum_{0 \neq \ss' + \ss' \leq \ss} C_{\ss}^{\ss', \ss''} \partial_x^{\ss'} \Lambda_0 \partial_x^{\ss - \ss' - \ss''} \lambda [ f_R^\epsilon - \rho_R^\epsilon M_{ \Lambda_0 } ]^\top \partial_x^{\ss''} \Lambda_0 \,,
	\end{aligned}
\end{equation*}
which, together with the Sobolev embedding $H_x^2 \hookrightarrow L^\infty_x$ and \eqref{EE-3}, implies
\begin{equation}\label{F1-2}
	\begin{aligned}
		| \partial_x^{\ss} P_{ T_{ \Lambda_0 } } ( \lambda [ f_R^\epsilon - \rho_R^\epsilon & M_{ \Lambda_0 } ] ) | \leq 5 | \partial_x^{\ss} \lambda [ f_R^\epsilon - \rho_R^\epsilon M_{ \Lambda_0 } ] | \\
		& + C ( 1 + \| \nabla_x \Lambda_0 \|^s_{H^{s+1}_x} ) \| \nabla_x \Lambda_0 \|^s_{H^{s+1}_x} \sum_{0 \neq \ss' \leq \ss} | \partial_x^{\ss - \ss'} \lambda [ f_R^\epsilon - \rho_R^\epsilon M_{ \Lambda_0 } ] | \,.
	\end{aligned}
\end{equation}
Following the similar arguments in \eqref{EE-4} and employing \eqref{HPoincare}, one has
\begin{equation}\label{F1-3}
	\begin{aligned}
		| \partial_x^{\ss} \lambda [ f_R^\epsilon - \rho_R^\epsilon M_{ \Lambda_0 } ] | \leq & \sqrt[4]{3} \| \partial_x^{\ss} ( \tfrac{ f_R^\epsilon }{ M_{ \Lambda_0 } } - \rho_R^\epsilon ) \|_{L^2_A ( M_{ \Lambda_0 } )} \\
		\leq & \tfrac{ \sqrt[4]{3} }{ \sqrt{ \lambda_0 } } \| \nabla_A \partial_x^{\ss} ( \tfrac{ f_R^\epsilon }{ M_{ \Lambda_0 } } ) \|_{L^2_A ( M_{ \Lambda_0 } )} \\
		& + C \| \nabla_x \Lambda_0 \|^{s^2}_{H^{s+1}_x} \sum_{0 \neq \ss' \leq \ss} \| \nabla_A \partial_x^{\ss - \ss'} ( \tfrac{ f_R^\epsilon }{ M_{ \Lambda_0 } } ) \|_{ L^2_A ( M_{ \Lambda_0 } ) } 
	\end{aligned}
\end{equation}
for $\ss \in \mathbb{N}^3$ with $1 \leq |\ss| = k \leq s$. Moreover, together with Lemma \ref{CE} and \eqref{F1-3}, one has
\begin{equation}\label{F1-4}
	\begin{aligned}
		\sum_{0 \neq \ss' \leq \ss} | \partial_x^{\ss - \ss'} \lambda [ f_R^\epsilon - \rho_R^\epsilon M_{ \Lambda_0 } ] | \lesssim ( 1 + \| \nabla_x \Lambda_0 \|^{s^2}_{H^{s+1}_x} ) \sum_{0 \neq \ss' \leq \ss} \| \nabla_A \partial_x^{\ss - \ss'} ( \tfrac{ f_R^\epsilon }{ M_{ \Lambda_0 } } ) \|_{ L^2_A ( M_{ \Lambda_0 } ) } \,.
	\end{aligned}
\end{equation}
Consequently, the bounds \eqref{F1-1}, \eqref{F1-2}, \eqref{F1-3} and \eqref{F1-4} reduce to
\begin{equation}\label{F1-5}
	\begin{aligned}
		& | \nabla_A ( A \cdot \partial_x^{\ss} P_{ T_{ \Lambda_0 } } ( \lambda [ f_R^\epsilon - \rho_R^\epsilon M_{ \Lambda_0 } ] ) ) | \\
		\leq & \tfrac{ 25 \sqrt[4]{3} }{ \sqrt{ \lambda_0 } } \| \nabla_A \partial_x^{\ss} ( \tfrac{ f_R^\epsilon }{ M_{ \Lambda_0 } } ) \|_{L^2_A ( M_{ \Lambda_0 } )} + C ( 1 + \| \nabla_x \Lambda_0 \|^{s^2}_{H^{s+1}_x} ) \sum_{0 \neq \ss' \leq \ss} \| \nabla_A \partial_x^{\ss - \ss'} ( \tfrac{ f_R^\epsilon }{ M_{ \Lambda_0 } } ) \|_{ L^2_A ( M_{ \Lambda_0 } ) } \,.
	\end{aligned}
\end{equation}
As a consequence, the similar arguments in \eqref{EE-bnd} indicate that
	\begin{align}\label{F1-bnd}
		\no | F_1 | \leq & \tfrac{ 25 \sqrt[4]{3} \nu_0 }{ c_1 \sqrt{ \lambda_0 } } \tfrac{1}{\epsilon} \| \nabla_A \partial_x^{\ss} ( \tfrac{ f_R^\epsilon }{ M_{ \Lambda_0 } } ) \|^2_{L^2_{x,A} ( M_{ \Lambda_0 } )} \\
		\no & + C \tfrac{1}{\epsilon} ( 1 + \| \nabla_x \Lambda_0 \|^{s^2}_{H^{s+1}_x} ) \| \nabla_A \partial_x^{\ss} ( \tfrac{ f_R^\epsilon }{ M_{ \Lambda_0 } } ) \|_{L^2_{x,A} ( M_{ \Lambda_0 } )} \\
		\no & \qquad \times \sum_{0 \neq \ss' \leq \ss} \| \nabla_A \partial_x^{\ss - \ss'} ( \tfrac{ f_R^\epsilon }{ M_{ \Lambda_0 } } ) \|_{ L^2_{x,A} ( M_{ \Lambda_0 } ) } \\
		\no \leq & \tfrac{ 25 \sqrt[4]{3} \nu_0 }{ c_1 \sqrt{ \lambda_0 } } \tfrac{1}{\epsilon} \| \nabla_A \partial_x^{\ss} ( \tfrac{ f_R^\epsilon }{ M_{ \Lambda_0 } } ) \|^2_{L^2_{x,A} ( M_{ \Lambda_0 } )} \\
		& + C ( 1 + \| \nabla_x \Lambda_0 \|^{s^2}_{H^{s+1}_x} ) \tfrac{1}{\sqrt{\epsilon}} \| \nabla_A \partial_x^{\ss} ( \tfrac{ f_R^\epsilon }{ M_{ \Lambda_0 } } ) \|_{L^2_{x,A} ( M_{ \Lambda_0 } )} \sum_{0 \leq \mathfrak{j} \leq k - 1} D_{\mathfrak{j}}^\frac{1}{2} (t) \,.
	\end{align}

\underline{\em Case 4.2. Control of $F_2$.} It is derived from \eqref{X3-2} and \eqref{F1-5} that
\begin{equation}\label{F2-bnd}
	\begin{aligned}
		| F_2 | \lesssim & \tfrac{1}{\epsilon} \sum_{0 \neq \ss' \leq \ss} \iint_{\R^3 \times \mathrm{SO(3)}} | \nabla_A \partial_x^{\ss} ( \tfrac{ f_R^\epsilon }{ M_{ \Lambda_0 } } ) | \cdot | \partial_x^{\ss'} M_{ \Lambda_0 } | \\
		& \qquad \qquad \qquad \times | \nabla_A ( A \cdot \partial_x^{\ss - \ss'} P_{ T_{ \Lambda_0 } } ( \lambda [ f_R^\epsilon - \rho_R^\epsilon M_{ \Lambda_0 } ] ) ) | \d A \d x \\
		\lesssim & \tfrac{1}{\epsilon} ( 1 + \| \nabla_x \Lambda_0 \|^{s^2}_{H^{s+1}_x} ) \| \nabla_x \Lambda_0 \|^s_{H^{s+1}_x} \| \nabla_A \partial_x^{\ss} ( \tfrac{ f_R^\epsilon }{ M_{ \Lambda_0 } } ) \|_{L^2_{x,A} ( M_{ \Lambda_0 } )} \\
		& \qquad \qquad \qquad \times \sum_{0 \neq \ss' \leq \ss} \| \nabla_A \partial_x^{\ss - \ss'} ( \tfrac{ f_R^\epsilon }{ M_{ \Lambda_0 } } ) \|_{ L^2_{x,A} ( M_{ \Lambda_0 } ) } \\
		\lesssim & ( 1 + \| \nabla_x \Lambda_0 \|^{s^2}_{H^{s+1}_x} ) \| \nabla_x \Lambda_0 \|^s_{H^{s+1}_x} \tfrac{1}{\sqrt{\epsilon}} \| \nabla_A \partial_x^{\ss} ( \tfrac{ f_R^\epsilon }{ M_{ \Lambda_0 } } ) \|_{L^2_{x,A} ( M_{ \Lambda_0 } )} \sum_{0 \leq \mathfrak{j} \leq k - 1} D_{\mathfrak{j}}^\frac{1}{2} (t) \,.
	\end{aligned}
\end{equation}

\underline{\em Case 4.3. Control of $F_3$.} By employing the bounds \eqref{C3-0}, \eqref{HPoincare} and \eqref{F1-5}, one can derived from the similar arguments in \eqref{C3-bnd} that
\begin{equation}\label{F3-bnd}
	\begin{aligned}
		| F_3 | \lesssim & \tfrac{1}{\epsilon} ( 1 + \| \nabla_x \Lambda_0 \|^{s^2 + s}_{H^{s+1}_x} ) \| \nabla_x \Lambda_0 \|^s_{H^{s+1}_x} \sum_{0 \neq \ss' \leq \ss} \| \nabla_A \partial_x^{\ss - \ss'} ( \tfrac{ f_R^\epsilon }{ M_{ \Lambda_0 } } ) \|_{ L^2_{x,A} ( M_{ \Lambda_0 } ) } \\
		& \times \big( \| \nabla_A \partial_x^{\ss} ( \tfrac{ f_R^\epsilon }{ M_{ \Lambda_0 } } ) \|_{ L^2_{x,A} ( M_{ \Lambda_0 } ) } + \| \partial_x^{\ss} ( \tfrac{ f_R^\epsilon }{ M_{ \Lambda_0 } } - \rho_R^\epsilon ) \|_{ L^2_{x,A} ( M_{ \Lambda_0 } ) } \big) \\
		\lesssim & ( 1 + \| \nabla_x \Lambda_0 \|^{3 s^2 }_{H^{s+1}_x} ) \| \nabla_x \Lambda_0 \|^s_{H^{s+1}_x} \sum_{0 \leq \mathfrak{j} \leq k - 1} \big( \tfrac{1}{\sqrt{\epsilon}} \| \nabla_A \partial_x^{\ss} ( \tfrac{ f_R^\epsilon }{ M_{ \Lambda_0 } } ) \|_{L^2_{x,A} ( M_{ \Lambda_0 } )} D_{\mathfrak{j}}^\frac{1}{2} (t) + D_{\mathfrak{j}} (t) \big) \,.
	\end{aligned}
\end{equation}

Therefore, plugging \eqref{F1-bnd}, \eqref{F2-bnd} and \eqref{F3-bnd} into \eqref{F1F2F3}, one gains
\begin{equation}\label{HEELO-bnd}
	\begin{aligned}
		& | - \iint_{\R^3 \times \mathrm{SO(3)}} \frac{1}{\epsilon} \partial_x^{\ss} \big[ \tfrac{1}{ M_{ \Lambda_0 } } L_R ( f_R^\epsilon - \rho_R^\epsilon M_{ \Lambda_0 } ) \big] \cdot \partial_x^{\ss} \big( \tfrac{ f_R^\epsilon }{ M_{ \Lambda_0 } } - \rho_R^\epsilon \big) M_{ \Lambda_0 } \d A \d x | \\
		\leq & \tfrac{ 25 \sqrt[4]{3} \nu_0 }{ c_1 \sqrt{ \lambda_0 } } \tfrac{1}{\epsilon} \| \nabla_A \partial_x^{\ss} ( \tfrac{ f_R^\epsilon }{ M_{ \Lambda_0 } } ) \|^2_{L^2_{x,A} ( M_{ \Lambda_0 } )} \\
		& + C ( 1 + \| \nabla_x \Lambda_0 \|^{3 s^2 }_{H^{s+1}_x} ) \| \nabla_x \Lambda_0 \|^s_{H^{s+1}_x} \sum_{0 \leq \mathfrak{j} \leq k - 1} \big( \tfrac{1}{\sqrt{\epsilon}} \| \nabla_A \partial_x^{\ss} ( \tfrac{ f_R^\epsilon }{ M_{ \Lambda_0 } } ) \|_{L^2_{x,A} ( M_{ \Lambda_0 } )} D_{\mathfrak{j}}^\frac{1}{2} (t) + D_{\mathfrak{j}} (t) \big) \,.
	\end{aligned}
\end{equation}

{\bf Step 5. Estimates of transport effect for $M_{ \Lambda_0 }$: $ - \partial_x^{\ss} \big[ \tfrac{ f_R^\epsilon }{ M_{ \Lambda_0 } } \tfrac{ ( \partial_t + A e_1 \cdot  \nabla_x ) M_{ \Lambda_0 } }{ M_{ \Lambda_0 } } \big] $.} Note that \eqref{M-der} shows
\begin{equation*}
	\begin{aligned}
		\tfrac{ ( \partial_t + A e_1 \cdot  \nabla_x ) M_{ \Lambda_0 } }{ M_{ \Lambda_0 } } = \tfrac{\nu_0}{d} A \cdot ( \partial_t + A e_1 \cdot \nabla_x ) \Lambda_0 \,.
	\end{aligned}
\end{equation*}
Then, by the Sobolev embedding theory,
	\begin{align}\label{HETEM-bnd}
		\no & | \iint_{\R^3 \times \mathrm{SO(3)}} - \partial_x^{\ss} \big[ \tfrac{ f_R^\epsilon }{ M_{ \Lambda_0 } } \tfrac{ ( \partial_t + A e_1 \cdot  \nabla_x ) M_{ \Lambda_0 } }{ M_{ \Lambda_0 } } \big] \cdot \partial_x^{\ss} \big( \tfrac{ f_R^\epsilon }{ M_{ \Lambda_0 } } - \rho_R^\epsilon \big) M_{ \Lambda_0 } \d A \d x | \\
		\no = & \big| - \tfrac{\nu_0}{d} \sum_{\ss' \leq \ss} C_{\ss}^{\ss'} \iint_{\R^3 \times \mathrm{SO(3)}} \partial_x^{\ss} \big( \tfrac{ f_R^\epsilon }{ M_{ \Lambda_0 } } - \rho_R^\epsilon \big) M_{ \Lambda_0 } \\
		\no & \qquad \qquad \times \big[ \partial_x^{\ss'} ( \tfrac{ f_R^\epsilon }{ M_{ \Lambda_0 } } - \rho_R^\epsilon ) + \partial_x^{\ss'} \rho_R^\epsilon \big] \cdot \big[ A \cdot ( \partial_t + A e_1 \cdot \nabla_x ) \partial_x^{\ss - \ss'} \Lambda_0 \big] \d A \d x \big| \\
		\no \lesssim & \sum_{\ss' \leq \ss} \| \partial_x^{\ss} \big( \tfrac{ f_R^\epsilon }{ M_{ \Lambda_0 } } - \rho_R^\epsilon \big) \|_{ L^2_{x,A} ( M_{ \Lambda_0 } ) } \big( \| \partial_x^{\ss'} ( \tfrac{ f_R^\epsilon }{ M_{ \Lambda_0 } } - \rho_R^\epsilon ) \|_{ L^2_{x,A} ( M_{ \Lambda_0 } ) } + \| \partial_x^{\ss'} \rho_R^\epsilon \|_{ L^2_x } \big) \| \partial_{t,x} \Lambda_0 \|_{ L^\infty_x } \\
		\lesssim & \| \partial_{t,x} \Lambda_0 \|_{ H^2_x } E_k^\frac{1}{2} (t) \sum_{0 \leq \mathfrak{j} \leq k} E_{ \mathfrak{j} }^\frac{1}{2} (t) \,.
	\end{align}

{\bf Step 6. Estimates of source term: $ \partial_x^{\ss} \big[ \tfrac{1}{ M_{ \Lambda_0 } } R ( f_1 ) \big] $.} The H\"older inequality indicates that
\begin{equation}\label{HEST-bnd}
	\begin{aligned}
		& | \iint_{\R^3 \times \mathrm{SO(3)}} \partial_x^{\ss} \big[ \tfrac{1}{ M_{ \Lambda_0 } } R ( f_1 ) \big] \cdot \partial_x^{\ss} \big( \tfrac{ f_R^\epsilon }{ M_{ \Lambda_0 } } - \rho_R^\epsilon \big) M_{ \Lambda_0 } \d A \d x | \\
		\lesssim & \| \partial_x^{\ss} ( \tfrac{ R ( f_1 ) }{ M_{ \Lambda_0 } } ) \|_{ L^2_{x,A} ( M_{ \Lambda_0 } ) } \| \partial_x^{\ss} \big( \tfrac{ f_R^\epsilon }{ M_{ \Lambda_0 } } - \rho_R^\epsilon \big) \|_{ L^2_{x,A} ( M_{ \Lambda_0 } ) } \\
		\lesssim & \| \partial_x^{\ss} ( \tfrac{ R ( f_1 ) }{ M_{ \Lambda_0 } } ) \|_{ L^2_{x,A} ( M_{ \Lambda_0 } ) } E_k^\frac{1}{2} (t) \,.
	\end{aligned}
\end{equation}

{\bf Step 7. Estimates of nonlinear term: $ \partial_x^{\ss} \big[ \tfrac{1}{ M_{ \Lambda_0 } } \widetilde{Q} ( f_R^\epsilon ) \big] $.} As illustrated in \eqref{Q-decomp}, we neglect the infinitesimal small quantity $\mathcal{ \epsilon }$ involved in $ \widetilde{Q} ( f_R^\epsilon ) $. Then
\begin{equation}\label{W1234}
	\begin{aligned}
		& \iint_{\R^3 \times \mathrm{SO(3)}} \partial_x^{\ss} \big[ \tfrac{1}{ M_{ \Lambda_0 } } \widetilde{Q} ( f_R^\epsilon ) \big] \cdot \partial_x^{\ss} \big( \tfrac{ f_R^\epsilon }{ M_{ \Lambda_0 } } - \rho_R^\epsilon \big) M_{ \Lambda_0 } \d A \d x \\
		= & \sum_{i=1}^4 \underbrace{ \iint_{\R^3 \times \mathrm{SO(3)}} \partial_x^{\ss} \big[ \tfrac{1}{ M_{ \Lambda_0 } } Q_i \big] \cdot \partial_x^{\ss} \big( \tfrac{ f_R^\epsilon }{ M_{ \Lambda_0 } } - \rho_R^\epsilon \big) M_{ \Lambda_0 } \d A \d x }_{:= W_i} \,,
	\end{aligned}
\end{equation}
where $Q_i$ ($i= 1,2,3,4$) are defined in \eqref{Q-decomp}.

\underline{\em Case 7.1. Control of $W_1$.} Recalling the definition of $Q_1$ in \eqref{Q-decomp}, one easily gains
\begin{equation}\label{W1-0}
	\begin{aligned}
		W_1 = \tfrac{ \nu_0 }{c_1} \sum_{ \ss' + \ss'' + \ss^* + \ss^\sharp = \ss } C_{\ss}^{ \ss', \ss'', \ss^*, \ss^\sharp } \iint_{\R^3 \times \mathrm{SO(3)}} \nabla_A \big\{ \partial_x^{\ss} \big( \tfrac{ f_R^\epsilon }{ M_{ \Lambda_0 } } - \rho_R^\epsilon \big) [ \partial_x^{\ss'} ( \tfrac{1}{ M_{ \Lambda_0 } } ) M_{ \Lambda_0 } ] \big\} \\
		\times \partial_x^{\ss''} ( \tfrac{1}{\rho_0} ) \partial_x^{\ss^*} (f_1 + f_R^\epsilon) \nabla_A ( A \cdot \partial_x^{\ss^\sharp} P_{ T_{ \Lambda_0 } } ( \lambda [ f_1 + f_R^\epsilon ] ) ) \d A \d x \,.
	\end{aligned}
\end{equation}

By employing the similar arguments in \eqref{C3-1}, one easily deduces that
\begin{equation}\label{W1-1}
	\begin{aligned}
		& | \nabla_A \big\{ \partial_x^{\ss} \big( \tfrac{ f_R^\epsilon }{ M_{ \Lambda_0 } } - \rho_R^\epsilon \big) [ \partial_x^{\ss'} ( \tfrac{1}{ M_{ \Lambda_0 } } ) M_{ \Lambda_0 } ] \big\} | \\
		\lesssim & ( 1 + \| \nabla_x \Lambda_0 \|^s_{ H^{s+1}_x } ) \big( | \nabla_A \partial_x^{\ss} ( \tfrac{ f_R^\epsilon }{ M_{ \Lambda_0 } } ) | + | \partial_x^{\ss} ( \tfrac{ f_R^\epsilon }{ M_{ \Lambda_0 } } - \rho_R^\epsilon )  | \big)
	\end{aligned}
\end{equation}
It is also derived from applying the similar arguments in Lemma 3.2 of \cite{JLT-M3AS-2019} and the Sobolev embedding $H^2_x \hookrightarrow L^\infty_x$ that
\begin{equation}\label{W1-2}
	\begin{aligned}
		| \partial_x^{\ss''} ( \tfrac{1}{\rho_0} ) | \lesssim \tfrac{1}{\rho_0} ( 1 + \| \nabla_x \rho_0 \|^s_{H^{s+1}_x} ) \lesssim 1 + \| \nabla_x \rho_0 \|^s_{H^{s+1}_x} \,.
	\end{aligned}
\end{equation}
where the last inequality follows from the fact $\inf_{(t,x)\in [0, T] \times \R^3} \rho_0 (t,x) > 0$. Moreover, by \eqref{X3-2}, one gains
\begin{equation}\label{W1-3}
	\begin{aligned}
		| \partial_x^{\ss^*} (f_1 + f_R^\epsilon) | \lesssim \sum_{ \tilde{\ss}^* \leq \ss^* } | \partial_x^{\ss^* - \tilde{\ss}^*} M_{ \Lambda_0 } | \big( | \partial_x^{\tilde{\ss}^*} ( \tfrac{f_1}{ M_{ \Lambda_0 } } ) | + | \partial_x^{\tilde{\ss}^*} ( \tfrac{ f_R^\epsilon }{ M_{ \Lambda_0 } } - \rho_R^\epsilon ) | + | \partial_x^{\tilde{\ss}^*} \rho_R^\epsilon | \big) \\
		\lesssim ( 1 + \| \nabla_x \Lambda_0 \|^s_{ H^{s+1}_x } ) \sum_{ \tilde{\ss}^* \leq \ss^* } \big( | \partial_x^{\tilde{\ss}^*} ( \tfrac{f_1}{ M_{ \Lambda_0 } } ) | + | \partial_x^{\tilde{\ss}^*} ( \tfrac{ f_R^\epsilon }{ M_{ \Lambda_0 } } - \rho_R^\epsilon ) | + | \partial_x^{\tilde{\ss}^*} \rho_R^\epsilon | \big) M_{ \Lambda_0 } \,.
	\end{aligned}
\end{equation}
Furthermore, by the similar arguments in \eqref{F1-2}, \eqref{F1-3} and \eqref{F1-4},
\begin{equation}\label{W1-4}
	\begin{aligned}
		& | \nabla_A ( A \cdot \partial_x^{\ss^\sharp} P_{ T_{ \Lambda_0 } } ( \lambda [ f_1 + f_R^\epsilon ] ) ) | \\
		\lesssim & ( 1 + \| \nabla_x \Lambda_0 \|^s_{H^{s+1}_x} ) \| \nabla_x \Lambda_0 \|^s_{H^{s+1}_x} \sum_{\tilde{\ss}^\sharp \leq \ss^\sharp} | \partial_x^{\tilde{\ss}^\sharp} \lambda [ f_1 + f_R^\epsilon ] | \\
		\lesssim & ( 1 + \| \nabla_x \Lambda_0 \|^{2s}_{H^{s+1}_x} ) \| \nabla_x \Lambda_0 \|^{2s}_{H^{s+1}_x} \\
		& \times \sum_{\tilde{\ss}^\sharp \leq \ss^\sharp} \big( \| \partial_x^{\tilde{\ss}^\sharp} ( \tfrac{f_1}{ M_{ \Lambda_0 } } ) \|_{ L^2_A ( M_{ \Lambda_0 } ) } + \| \partial_x^{\tilde{\ss}^\sharp} ( \tfrac{ f_R^\epsilon }{ M_{ \Lambda_0 } } - \rho_R^\epsilon ) \|_{ L^2_A ( M_{ \Lambda_0 } ) } + | \partial_x^{\tilde{\ss}^\sharp} \rho_R^\epsilon | \big) \,.
	\end{aligned}
\end{equation}
Then we will use the following idea to estimate the quantity $W_1$. For the general form $\sum_{\ss^* + \ss^\sharp \leq \ss } \int_{\R^3} |\partial_x^{\ss^*} X | \, |\partial_x^{\ss^\sharp} Y | \, |\partial_x^{\ss} Z | \d x $, we dominate
\begin{equation*}
	\begin{aligned}
		& \int_{\R^3} |\partial_x^{\ss} X | \, | Y | \, |\partial_x^{\ss} Z | \d x + \int_{\R^3} | X | \, |\partial_x^{\ss} Y | \, |\partial_x^{\ss} Z | \d x \\
		\leq & \| \partial_x^{\ss} X \|_{L^2_x} \| Y \|_{L^\infty_x} \| \partial_x^{\ss} Z \|_{L^2_x} + \| X \|_{L^\infty_x} \| \partial_x^{\ss} Z \|_{L^2_x} \| \partial_x^{\ss} Z \|_{L^2_x} \\
		\leq & C \| \partial_x^{\ss} X \|_{L^2_x} \| Y \|_{H^2_x} \| \partial_x^{\ss} Z \|_{L^2_x} + \| X \|_{H^2_x} \| \partial_x^{\ss} Z \|_{L^2_x} \| \partial_x^{\ss} Z \|_{L^2_x}
	\end{aligned}
\end{equation*}
by employing the Sobolev embedding $H^2_x \hookrightarrow L^\infty_x$, and
\begin{equation*}
	\begin{aligned}
		\sum_{ \substack{ \ss^* + \ss^\sharp \leq \ss \\ \ss^*, \ss^\sharp \neq \ss } } \int_{\R^3} |\partial_x^{\ss^*} X | \, |\partial_x^{\ss^\sharp} Y | \, |\partial_x^{\ss} Z | \d x \leq \sum_{ \substack{ \ss^* + \ss^\sharp \leq \ss \\ \ss^*, \ss^\sharp \neq \ss } } \| \partial_x^{\ss^*} X \|_{L^4_x} \| \partial_x^{\ss^\sharp} Y \|_{L^4_x} \|\partial_x^{\ss} Z \|_{L^2_x} \\
		\leq C \| X \|_{H^{|\ss|}_x} \| Y \|_{H^{|\ss|}_x} \|\partial_x^{\ss} Z \|_{L^2_x} 
	\end{aligned}
\end{equation*}
by applying the Sobolev embedding $H^1_x \hookrightarrow L^4_x$. Following the previous ideas and combining with the bounds \eqref{HPoincare}-\eqref{W1-0}-\eqref{W1-1}-\eqref{W1-2}-\eqref{W1-3}-\eqref{W1-4}, one easily controls the quantity $W_1$ by
\begin{equation}\label{W1-bnd}
	\begin{aligned}
		|W_1| \lesssim & ( 1 + \| \nabla_x \rho_0 \|^s_{H^{s+1}_x} ) ( 1 + \| \nabla_x \Lambda_0 \|^{4s}_{H^{s+1}_x} ) \| \nabla_x \Lambda_0 \|^{2s}_{H^{s+1}_x} \\
		& \times \sum_{|\ss'| \leq s} \big( \| \partial_x^{\tilde{\ss'}} ( \tfrac{f_1}{ M_{ \Lambda_0 } } ) \|^2_{ L^2_{x,A} ( M_{ \Lambda_0 } ) } + \| \partial_x^{\tilde{\ss'}} ( \tfrac{ f_R^\epsilon }{ M_{ \Lambda_0 } } - \rho_R^\epsilon ) \|^2_{ L^2_{x,A} ( M_{ \Lambda_0 } ) } + \| \partial_x^{\tilde{\ss'}} \rho_R^\epsilon \|^2_{L^2_x} \big) \\
		& \times \big( \| \nabla_A \partial_x^{\ss} ( \tfrac{ f_R^\epsilon }{ M_{ \Lambda_0 } } ) \|_{ L^2_{x,A} ( M_{ \Lambda_0 } ) } + \| \partial_x^{\ss} ( \tfrac{ f_R^\epsilon }{ M_{ \Lambda_0 } } - \rho_R^\epsilon ) \|_{ L^2_{x,A} ( M_{ \Lambda_0 } ) } \big) \\
		\lesssim & \sqrt{\epsilon} ( 1 + \| \nabla_x \rho_0 \|^s_{H^{s+1}_x} ) ( 1 + \| \nabla_x \Lambda_0 \|^{2s^2 + 4s}_{H^{s+1}_x} ) \| \nabla_x \Lambda_0 \|^{2s}_{H^{s+1}_x} \\
		& \times ( \tfrac{1}{\sqrt{\epsilon}} \| \nabla_A \partial_x^{\ss} ( \tfrac{ f_R^\epsilon }{ M_{ \Lambda_0 } } ) \|_{L^2_{x,A} ( M_{ \Lambda_0 } )} + \sum_{0 \leq \mathfrak{j} \leq k-1} D_{\mathfrak{j}}^\frac{1}{2} (t) ) \sum_{ 0 \leq \mathfrak{j} \leq s } E_{\mathfrak{j}} (t) 
	\end{aligned}
\end{equation}
for $s \geq 2$.

\underline{\em Case 7.2. Control of $W_2$.} Recalling the definition of $Q_2$ in \eqref{Q-decomp}, it infers that
\begin{equation*}
	\begin{aligned}
		W_2 = & - \tfrac{\nu_0}{2 c_1^2} \sum_{\ss' + \ss'' \leq \ss} C_{\ss}^{\ss', \ss''} \iint_{\R^3 \times \mathrm{SO(3)}} \nabla_A \big\{ \partial_x^{\ss} \big( \tfrac{ f_R^\epsilon }{ M_{ \Lambda_0 } } - \rho_R^\epsilon \big) [ \partial_x^{\ss'} ( \tfrac{1}{ M_{ \Lambda_0 } } ) M_{ \Lambda_0 } ] \big\} \partial_x^{\ss''} ( \tfrac{1}{\rho_0} ) \\
		& \times \nabla_A \partial_x^{\ss - \ss' - \ss''} \big[ A \cdot \lambda [ f_1 + f_R^\epsilon ] \big( \lambda [ f_1 + f_R^\epsilon ]^\top \Lambda_0 + \Lambda_0^\top \lambda [ f_1 + f_R^\epsilon ] \big) \big] \d A \d x \,.
	\end{aligned}
\end{equation*}
The bounds \eqref{W1-1} and \eqref{W1-2} tell us
\begin{equation*}
	\begin{aligned}
		& | \nabla_A \big\{ \partial_x^{\ss} \big( \tfrac{ f_R^\epsilon }{ M_{ \Lambda_0 } } - \rho_R^\epsilon \big) [ \partial_x^{\ss'} ( \tfrac{1}{ M_{ \Lambda_0 } } ) M_{ \Lambda_0 } ] \big\} \partial_x^{\ss''} ( \tfrac{1}{\rho_0} ) | \\
		\lesssim & ( 1 + \| \nabla_x \rho_0 \|^s_{H^{s+1}_x} ) ( 1 + \| \nabla_x \Lambda_0 \|^s_{ H^{s+1}_x } ) \big( | \nabla_A \partial_x^{\ss} ( \tfrac{ f_R^\epsilon }{ M_{ \Lambda_0 } } ) | + | \partial_x^{\ss} ( \tfrac{ f_R^\epsilon }{ M_{ \Lambda_0 } } - \rho_R^\epsilon )  | \big) \,.
	\end{aligned}
\end{equation*}
Moreover, together with \eqref{W1-4} and $H^2_x \hookrightarrow L^\infty_x$,
	\begin{align*}
		& | \nabla_A \partial_x^{\ss - \ss' - \ss''} \big[ A \cdot \lambda [ f_1 + f_R^\epsilon ] \big( \lambda [ f_1 + f_R^\epsilon ]^\top \Lambda_0 + \Lambda_0^\top \lambda [ f_1 + f_R^\epsilon ] \big) \big] | \\
		\lesssim & ( 1 + \| \nabla_x \Lambda_0 \|_{ H^{s+1}_x } ) \sum_{ \tilde{\ss} \leq \ss - \ss' - \ss'' } | \partial_x^{\tilde{\ss}} \lambda [ f_1 + f_R^\epsilon ] | \, | \partial_x^{\ss - \ss' - \ss'' - \tilde{\ss}} \lambda [ f_1 + f_R^\epsilon ] | \\
		\lesssim & ( 1 + \| \nabla_x \Lambda_0 \|^{2s+1}_{ H^{s+1}_x } ) \| \nabla_x \Lambda_0 \|^{2s}_{ H^{s+1}_x } \\
		& \times \sum_{ \tilde{\ss} \leq \ss - \ss' - \ss'' } \bigg\{ \sum_{\tilde{\ss}' \leq \tilde{\ss}} \big( \| \partial_x^{\tilde{\ss}'} ( \tfrac{f_1}{ M_{ \Lambda_0 } } ) \|_{ L^2_A ( M_{ \Lambda_0 } ) } + \| \partial_x^{\tilde{\ss}'} ( \tfrac{ f_R^\epsilon }{ M_{ \Lambda_0 } } - \rho_R^\epsilon ) \|_{ L^2_A ( M_{ \Lambda_0 } ) } + | \partial_x^{\tilde{\ss}'} \rho_R^\epsilon | \big) \\
		& \times \sum_{\tilde{\ss}'' \leq \ss - \ss' - \ss'' - \tilde{\ss}} \big( \| \partial_x^{\tilde{\ss}''} ( \tfrac{f_1}{ M_{ \Lambda_0 } } ) \|_{ L^2_A ( M_{ \Lambda_0 } ) } + \| \partial_x^{\tilde{\ss}''} ( \tfrac{ f_R^\epsilon }{ M_{ \Lambda_0 } } - \rho_R^\epsilon ) \|_{ L^2_A ( M_{ \Lambda_0 } ) } + | \partial_x^{\tilde{\ss}''} \rho_R^\epsilon | \big) \bigg\} \,.
	\end{align*}
It therefore follows from the similar arguments in \eqref{W1-bnd} that
\begin{equation}\label{W2-bnd}
	\begin{aligned}
		| W_2 | \lesssim & \sqrt{\epsilon} ( 1 + \| \nabla_x \rho_0 \|^s_{H^{s+1}_x} ) ( 1 + \| \nabla_x \Lambda_0 \|^{2s^2 + 3s+1}_{ H^{s+1}_x } ) \| \nabla_x \Lambda_0 \|^{2s}_{ H^{s+1}_x } \\
		& \times ( \tfrac{1}{\sqrt{\epsilon}} \| \nabla_A \partial_x^{\ss} ( \tfrac{ f_R^\epsilon }{ M_{ \Lambda_0 } } ) \|_{L^2_{x,A} ( M_{ \Lambda_0 } )} + \sum_{0 \leq \mathfrak{j} \leq k-1} D_{\mathfrak{j}}^\frac{1}{2} (t) ) \sum_{ 0 \leq \mathfrak{j} \leq s } E_{\mathfrak{j}} (t) \,.
	\end{aligned}
\end{equation}

\underline{\em Case 7.3. Controls of $W_3$ and $W_4$.} By employing the similar estimates of the bound \eqref{W2-bnd}, one easily derives that
\begin{equation}\label{W34-bnd}
	\begin{aligned}
		|W_3| + |W_4| \lesssim & \sqrt{\epsilon} ( 1 + \| \nabla_x \rho_0 \|^{2s}_{H^{s+1}_x} ) ( 1 + \| \nabla_x \Lambda_0 \|^{2s^2 + 3s+1}_{ H^{s+1}_x } ) \| \nabla_x \Lambda_0 \|^{2s}_{ H^{s+1}_x } \\
		& \times ( \tfrac{1}{\sqrt{\epsilon}} \| \nabla_A \partial_x^{\ss} ( \tfrac{ f_R^\epsilon }{ M_{ \Lambda_0 } } ) \|_{L^2_{x,A} ( M_{ \Lambda_0 } )} + \sum_{0 \leq \mathfrak{j} \leq k-1} D_{\mathfrak{j}}^\frac{1}{2} (t) ) \sum_{ 0 \leq \mathfrak{j} \leq s } E_{\mathfrak{j}} (t) \,.
	\end{aligned}
\end{equation}

As a result, we derive from plugging \eqref{W1-bnd}, \eqref{W2-bnd} and \eqref{W34-bnd} into \eqref{W1234} that
\begin{equation}\label{HENT-bnd}
	\begin{aligned}
		& | \iint_{\R^3 \times \mathrm{SO(3)}} \partial_x^{\ss} \big[ \tfrac{1}{ M_{ \Lambda_0 } } \widetilde{Q} ( f_R^\epsilon ) \big] \cdot \partial_x^{\ss} \big( \tfrac{ f_R^\epsilon }{ M_{ \Lambda_0 } } - \rho_R^\epsilon \big) M_{ \Lambda_0 } \d A \d x | \\
		\lesssim & \sqrt{\epsilon} ( 1 + \| \nabla_x \rho_0 \|^{2s}_{H^{s+1}_x} ) ( 1 + \| \nabla_x \Lambda_0 \|^{2s^2+4s}_{ H^{s+1}_x } ) \| \nabla_x \Lambda_0 \|^{2s}_{ H^{s+1}_x } \\
		& \times ( \tfrac{1}{\sqrt{\epsilon}} \| \nabla_A \partial_x^{\ss} ( \tfrac{ f_R^\epsilon }{ M_{ \Lambda_0 } } ) \|_{L^2_{x,A} ( M_{ \Lambda_0 } )} + \sum_{0 \leq \mathfrak{j} \leq k-1} D_{\mathfrak{j}}^\frac{1}{2} (t) ) \sum_{ 0 \leq \mathfrak{j} \leq s } E_{\mathfrak{j}} (t) \,.
	\end{aligned}
\end{equation}

In summary, by substituting the bounds \eqref{HETEMi-bnd}, \eqref{HETEMa-bnd}, \eqref{HEC-bnd}, \eqref{HEELO-bnd}, \eqref{HETEM-bnd}, \eqref{HEST-bnd} and \eqref{HENT-bnd} into \eqref{HigherOr-der}, one immediately gains
\begin{equation}\label{Hder-1}
	\begin{aligned}
		& \tfrac{1}{2} \tfrac{\d}{\d t} \big( \| \partial_x^{\ss} ( \tfrac{ f_R^\epsilon }{ M_{ \Lambda_0 } } - \rho_R^\epsilon ) \|^2_{ L^2_{x,A} ( M_{ \Lambda_0 } ) } + \| \partial_x^{\ss} \rho_R^\epsilon \|^2_{L^2_x} \big) + d_\star \tfrac{1}{\epsilon} \| \nabla_A \partial_x^{\ss} ( \tfrac{ f_R^\epsilon }{ M_{ \Lambda_0 } } ) \|^2_{L^2_{x,A} ( M_{ \Lambda_0 } ) } \\
		\lesssim & \Phi_s (\rho_0, \Lambda_0) \Psi_s (f_1) \bigg\{ \tfrac{1}{\sqrt{\epsilon}} \| \nabla_A \partial_x^{\ss} ( \tfrac{ f_R^\epsilon }{ M_{ \Lambda_0 } } ) \|_{L^2_{x,A} ( M_{ \Lambda_0 } )} \Big( \sum_{0 \leq \mathfrak{j} \leq k - 1} D_{\mathfrak{j}}^\frac{1}{2} (t) + \sqrt{\epsilon} \sum_{ 0 \leq \mathfrak{j} \leq s } E_{\mathfrak{j}} (t) \Big)  \\
		& + \sum_{0 \leq \mathfrak{j} \leq k - 1} D_{\mathfrak{j}} (t) + \sum_{0 \leq \mathfrak{j} \leq k} \big( E_{\mathfrak{j}} (t) + E_{\mathfrak{j}}^\frac{1}{2} (t) \big) + \sqrt{\epsilon} \sum_{0 \leq \mathfrak{j} \leq k-1} D_{\mathfrak{j}}^\frac{1}{2} (t) \sum_{ 0 \leq \mathfrak{j} \leq s } E_{\mathfrak{j}} (t) \bigg\}
	\end{aligned}
\end{equation}
for $\ss \in \mathbb{N}^3$ with $1 \leq |\ss| = k \leq s$ ($s \geq 2$), where the constant $d_\star = d - \frac{25 \sqrt[4]{3} \nu_0 }{ c_1 \lambda_0 } > 0$, and
\begin{equation}
	\begin{aligned}
		& \Phi_s (\rho_0, \Lambda_0) = ( 1 + \| \nabla_x \rho_0 \|^{2s}_{H^{s+2}_x} ) ( 1 + \| \nabla_x \Lambda_0 \|^{3s^2 + 4s}_{ H^{s+2}_x } ) ( 1 + \| \partial_t \Lambda_0 \|_{ H^{s+2}_x } ) \,, \\
		& \Psi_s (f_1) = 1 + \| \tfrac{f_1}{ M_{ \Lambda_0 } } \|_{H^s_x L^2_A ( M_{ \Lambda_0 } )} + \| \tfrac{ R ( f_1 ) }{ M_{ \Lambda_0 } } \|_{H^s_x L^2_A ( M_{ \Lambda_0 } )} \,.
	\end{aligned}
\end{equation}
Theorem \ref{LWP} indicates that 
$$ \Phi_s (\rho_0, \Lambda_0) \leq C ( \| \nabla_x \rho_0^{in} \|_{H^{s+2}}, \| \nabla_x \Lambda_0^{in} \|_{H^{s+2}} ) \,, $$
and Lemma \ref{Lmm-f1} reads
\begin{equation*}
	\begin{aligned}
		\Psi_s (f_1) \leq C ( \| \nabla_x \rho_0^{in} \|_{H^{s+3}}, \| \nabla_x \Lambda_0^{in} \|_{H^{s+3}} ) \,.
	\end{aligned}
\end{equation*}
Summing up \eqref{Hder-1} for $|\ss| = k$ and employing the Young's inequality, one then obtains
\begin{equation}\label{Hk-bnd}
	\begin{aligned}
		\tfrac{\d}{\d t} E_k (t) + d_\star D_k (t) \leq C \Big[ \sum_{0 \leq \mathfrak{j} \leq k - 1} D_{\mathfrak{j}} (t) + \sum_{ 0 \leq \mathfrak{j} \leq s } \big( E_{\mathfrak{j}}^\frac{1}{2} (t) + E_{\mathfrak{j}} (t) + \epsilon E_{\mathfrak{j}}^2 (t) \big) \Big]
	\end{aligned}
\end{equation}
for $1 \leq k \leq s$ with $s \geq 2$, where the constant $C > 0$ depends on $\| \nabla_x \rho_0^{in} \|_{H^{s+3}}$ and $ \| \nabla_x \Lambda_0^{in} \|_{H^{s+3}} $.

\subsubsection{Close the estimates: finishing the proof of Lemma \ref{UEE}}

In this subsection, based on the estimates \eqref{L2-bnd} and \eqref{Hk-bnd}, we will finish the proof of the uniform-in-$\epsilon$ estimate \eqref{UEE-ineq} in Lemma \ref{UEE}. Observe that the main difficulty is to absorb the term $ \sum_{0 \leq \mathfrak{j} \leq k - 1} D_{\mathfrak{j}} (t) $ in \eqref{Hk-bnd}. As inspired in Section 4 of \cite{JL-APDE-2022}, we can apply the induction arguments for $k = 0, 1, \cdots, s$ with $s \geq 2$. More precisely, there exist some positive constants $c_k'$, $c_k'' > 0$ such that
\begin{equation*}
	\begin{aligned}
		\tfrac{\d}{\d t} \big( \sum_{k=0}^s c_k' E_k (t) \big) + \sum_{k=0}^s c_k'' D_k (t) \lesssim \sum_{ 0 \leq \mathfrak{j} \leq s } \big( E_{\mathfrak{j}}^\frac{1}{2} (t) + E_{\mathfrak{j}} (t) + \epsilon E_{\mathfrak{j}}^2 (t) \big) \\
		\lesssim \big( \sum_{k=0}^s c_k' E_k (t) \big)^\frac{1}{2} + \sum_{k=0}^s c_k' E_k (t) + \epsilon \big( \sum_{k=0}^s c_k' E_k (t) \big)^2 \,.
	\end{aligned}
\end{equation*}
Hence the bound \eqref{UEE-ineq} holds. Then the proof of Lemma \ref{UEE} is finished.

\subsection{Hydrodynamic limit for SOKB system: Complete of proof of Theorem \ref{Hydrodynamic}}

    In this subsection, we give the proof of Theorem $\ref{Hydrodynamic}$ by virtue of Lemma \ref{UEE}. Given initial data $f_{R}^{\epsilon,in},$ it follows from the local existence result of the remainder equation \eqref{RE} that, there exists a unique local solution $f_{R}^{\epsilon}(t\,,x\,,A)$ on $[0\,,T_{\epsilon}]$ with $T_{\epsilon}>0$ be the maximal lifespan. For simplicity, the details of proof for the local existence are omitted here. The similar arguments can be found in Section 5 of \cite{JLZ-ARMA-2020}.
    
    Based on the uniform-in-$\epsilon$ estimate \eqref{UEE-ineq} in Lemma \ref{UEE}, we mainly show that the lifespan $T_\epsilon \geq T$ uniformly in sufficiently small $\epsilon > 0$, where $T > 0$ is constructed in Theorem \ref{LWP} with $m = s + 4$. Recall that
    \begin{equation*}
    	\begin{aligned}
    		\tfrac{\d}{\d t} \mathcal{E} (t) + \mathcal{D} (t) \leq C ( 1 + \mathcal{E} (t) ) + C \epsilon \mathcal{E}^{2} (t) \,,
    	\end{aligned}
    \end{equation*}
    which means that for all $t \in [0, T_\epsilon]$,
    \begin{equation}
    	\begin{aligned}
    		\tfrac{\d}{\d t} \mathbb{E} (t) + \mathcal{D} (t) \leq C \mathbb{E} (t) + C \epsilon \mathbb{E}^{2} (t) \,,
    	\end{aligned}
    \end{equation}
    where $ \mathbb{E} (t) = 1 + \mathcal{E} (t) $.
    
    Note that the initial data $f_R^{\epsilon, in}$ satisfies \eqref{IC-fR-bnd}, i.e.,
    \begin{equation*}
    	\begin{aligned}
    		K^{in} : = \sup_{\epsilon \in (0,1)} \big( \| \rho_R^{\epsilon , in} \|_{H^s_x} + \| \tfrac{ f_R^{ \epsilon, in } }{ M_{ \Lambda_0^{in} } } - \rho_R^{\epsilon , in} \|_{ H^s_x L^2_A ( M_{ \Lambda_0^{in} } ) } \big) < \infty \,,
    	\end{aligned}
    \end{equation*}
    where $\rho_R^{\epsilon , in} = \int_{\mathrm{SO}(3)} f_R^{\epsilon , in} \d A $. Then it is easy to see 
    \begin{equation*}
    	\begin{aligned}
    		\mathbb{E} (0) = & 1 + \sum_{0 \leq k \leq s} c_k' \sum_{|\ss| = k} \big( \| \partial_x^{\ss} \rho_R^{\epsilon , in} \|_{H^s_x} + \| \partial_x^{\ss} ( \tfrac{ f_R^{ \epsilon, in } }{ M_{ \Lambda_0^{in} } } - \rho_R^{\epsilon , in} ) \|_{ H^s_x L^2_A ( M_{ \Lambda_0^{in} } ) } \big) \\
    		\leq & 1 + C_1 K^{in} : = K_*^{in} \,.
    	\end{aligned}
    \end{equation*}
    Denote by 
    \begin{equation*}
    	\begin{aligned}
    		\Upsilon = 2 e^{C T} ( K_*^{in} + C T ) > \mathbb{E} (0) \,.
    	\end{aligned}
    \end{equation*}
    Let 
    \begin{equation}\label{T-tilde}
    	\begin{aligned}
    		\widetilde{T} = \sup \{ t \in [ 0, T_\epsilon ]; \mathbb{E} (t) \leq \Upsilon \} \leq T_\epsilon \,.
    	\end{aligned}
    \end{equation}
    The continuity of $\mathbb{E} (t)$ in $t$ shows that $ \widetilde{T} > 0 $. Then for any $t \in [ 0, \widetilde{T} ]$,
    \begin{equation*}
    	\begin{aligned}
    		\tfrac{\d}{\d t} \mathbb{E} (t) + \mathcal{D} (t) \leq C \mathbb{E} (t) + C \epsilon \mathbb{E}^{2} (t) \leq C \mathbb{E} (t) + C \epsilon \Upsilon^2 \,.
    	\end{aligned}
    \end{equation*}
    Then the Gr\"onwall inequality implies
    \begin{equation}\label{E-bnd}
    	\begin{aligned}
    		\mathbb{E} (t) \leq e^{C t} ( \mathbb{E} (0) + C \epsilon \Upsilon^2 t ) \leq e^{C t} ( K_*^{in} + C \epsilon \Upsilon^2 t )
    	\end{aligned}
    \end{equation}
    for any $t \in [ 0, \widetilde{T} ]$.
    
    We then claim that $\widetilde{T} \geq T$ when $0 < \epsilon < \epsilon_0 : = \frac{1}{C \Upsilon^2} $. Indeed, if $\widetilde{T} < T$, the bound \eqref{E-bnd} reduces to
    \begin{equation*}
    	\begin{aligned}
    		\sup_{ t \in [0, \widetilde{T}] } \mathbb{E} (t) \leq e^{C T} ( K_*^{in} + C \epsilon \Upsilon^2 T ) \leq e^{C T} ( K_*^{in} + C T ) = \tfrac{1}{2} \Upsilon < \Upsilon \,.
    	\end{aligned}
    \end{equation*}
    Then the continuity of $\mathbb{E} (t)$ shows that there is a $t_* > 0$ such that
    \begin{equation*}
    	\begin{aligned}
    		\sup_{ t \in [0, \widetilde{T} + t_*] } \mathbb{E} (t) \leq \Upsilon \,,
    	\end{aligned}
    \end{equation*}
    which contracts to the definition of $\widetilde{T}$ in \eqref{T-tilde}. As a result, $T_\epsilon \geq \widetilde{T} \geq T$. Consequently, one easily has
    \begin{equation*}
    	\begin{aligned}
    		\mathbb{E} (t) + \int_0^t \mathcal{D} (\tau) \d \tau \leq C \Upsilon T + C T
    	\end{aligned}
    \end{equation*}
    for all $t \in [ 0, T ]$ and $0 < \epsilon < \epsilon_0$. Then the proof of Theorem \ref{Hydrodynamic} is completed.

\section{Bounds for expanded term $f_1$: Proof of Lemma \ref{Lmm-f1}}\label{Sec:f1}

In this section, we will give the proof of Lemma \ref{Lmm-f1}. We rewrite \eqref{f1} as
\begin{equation}\label{f1-eq}
	\begin{aligned}
		\tfrac{1}{ M_{ \Lambda_0 } } \L_{ M_{ \Lambda_0 } } f_1 = h_0 : = \tfrac{1}{ M_{ \Lambda_0 } } ( \partial_t + A e_1 \cdot \nabla_x ) f_0 
	\end{aligned}
\end{equation}
with $\int_{\mathrm{SO}(3)} f_1 \d A = 0$ and $ P_{ T_{ \Lambda_0 } } ( \lambda [ f_1 ] ) = 0 $.

Multiplying \eqref{f1-eq} by $f_1$ and integrating by parts over $(x, A) \in \R^3 \times \mathrm{SO(3)}$, one has
\begin{equation*}
	\begin{aligned}
		\iint_{\R^3 \times \mathrm{SO(3)}} \tfrac{1}{ M_{ \Lambda_0 } } \nabla_A \cdot \big[ M_{ \Lambda_0 } \nabla_A ( \tfrac{f_1}{ M_{ \Lambda_0 } } ) \big] f_1 \d A \d x = \iint_{\R^3 \times \mathrm{SO(3)}} h_0 f_1 \d A \d x \,,
	\end{aligned}
\end{equation*}
which immediately implies
\begin{equation*}
	\begin{aligned}
		\| \nabla_A ( \tfrac{f_1}{ M_{ \Lambda_0 } } ) \|^2_{ L^2_{x,A} ( M_{ \Lambda_0 } ) } \leq \| h_0 \|_{ L^2_{x,A} ( M_{ \Lambda_0 } ) } \| \tfrac{f_1}{ M_{ \Lambda_0 } } \|_{ L^2_{x,A} ( M_{ \Lambda_0 } ) } \,.
	\end{aligned}
\end{equation*}
Note that $\int_{\mathrm{SO}(3)} \tfrac{f_1}{ M_{ \Lambda_0 } } M_{ \Lambda_0 } \d A = 0 $. Then the Poincar\'e inequality in Lemma \ref{CE} indicates that
\begin{equation*}
	\begin{aligned}
		\lambda_0 \| \tfrac{f_1}{ M_{ \Lambda_0 } } \|^2_{ L^2_{x,A} ( M_{ \Lambda_0 } ) } \leq \| \nabla_A ( \tfrac{f_1}{ M_{ \Lambda_0 } } ) \|^2_{ L^2_{x,A} ( M_{ \Lambda_0 } ) } \,.
	\end{aligned}
\end{equation*}
As a result,
\begin{equation}\label{L2-0}
	\begin{aligned}
		\| \tfrac{f_1}{ M_{ \Lambda_0 } } \|^2_{ L^2_{x,A} ( M_{ \Lambda_0 } ) } + \| \nabla_A ( \tfrac{f_1}{ M_{ \Lambda_0 } } ) \|^2_{ L^2_{x,A} ( M_{ \Lambda_0 } ) } \lesssim \| h_0 \|^2_{ L^2_{x,A} ( M_{ \Lambda_0 } ) } \,.
	\end{aligned}
\end{equation}

Next we control the higher order spatial derivatives of $\tfrac{f_1}{ M_{ \Lambda_0 } }$. For any multi-index $\ss \in \mathbb{N}^3$ with $1 \leq |\ss| = k \leq s$ ($s \geq 2$), we apply the derivative operator $\partial_x^{\ss}$ to \eqref{f1-eq}, and then obtain
\begin{equation}\label{Hf1-eq}
	\begin{aligned}
		\partial_x^{\ss} \big[ \tfrac{1}{ M_{ \Lambda_0 } } \L_{ M_{ \Lambda_0 } } f_1 \big] = \partial_x^{\ss} h_0 \,.
	\end{aligned}
\end{equation}
Multiplying \eqref{Hf1-eq} by $ \partial_x^{\ss} ( \tfrac{f_1}{ M_{ \Lambda_0 } } ) M_{ \Lambda_0 } $ and integrating over $(x, A) \in \R^3 \times \mathrm{SO(3)}$, it infers
\begin{equation}\label{f1-1}
	\begin{aligned}
		\iint_{\R^3 \times \mathrm{SO(3)}} \partial_x^{\ss} \big[ \tfrac{1}{ M_{ \Lambda_0 } } \L_{ M_{ \Lambda_0 } } f_1 \big] \cdot \partial_x^{\ss} ( \tfrac{f_1}{ M_{ \Lambda_0 } } ) M_{ \Lambda_0 } \d A \d x = \iint_{\R^3 \times \mathrm{SO(3)}} \partial_x^{\ss} h_0 \cdot \partial_x^{\ss} ( \tfrac{f_1}{ M_{ \Lambda_0 } } ) M_{ \Lambda_0 } \d A \d x \,.
	\end{aligned}
\end{equation}
Following the same procedure of estimate \eqref{HEC-bnd},
	\begin{align}\label{f1-2}
		\no & - \iint_{\R^3 \times \mathrm{SO(3)}} \partial_x^{\ss} \big[ \tfrac{1}{ M_{ \Lambda_0 } } \L_{ M_{ \Lambda_0 } } f_1 \big] \cdot \partial_x^{\ss} ( \tfrac{f_1}{ M_{ \Lambda_0 } } ) M_{ \Lambda_0 } \d A \d x \\
		\no \geq & d \| \nabla_A \partial_x^{\ss} ( \tfrac{f_1}{ M_{ \Lambda_0 } } ) \|^2_{L^2_{x,A} ( M_{ \Lambda_0 } ) } \\
		\no & - C ( 1 + \| \nabla_x \Lambda_0 \|^s_{ H_x^{ s + 1 } } ) \| \nabla_x \Lambda_0 \|^{s^2 + s}_{ H_x^{ s + 1 } } \sum_{0 \neq \ss' \leq \ss } \| \nabla_A \partial_x^{\ss - \ss'} ( \tfrac{f_1}{ M_{ \Lambda_0 } } ) \|^2_{L^2_{x,A} ( M_{ \Lambda_0 } ) } \\
		\no & - C ( 1 + \| \nabla_x \Lambda_0 \|^s_{ H_x^{ s + 1 } } ) \| \nabla_x \Lambda_0 \|^s_{ H_x^{ s + 1 } } \| \nabla_A \partial_x^{\ss} ( \tfrac{f_1}{ M_{ \Lambda_0 } } ) \|_{L^2_{x,A} ( M_{ \Lambda_0 } ) } \\
		\no & \qquad \qquad \qquad \qquad \qquad \qquad \times \sum_{0 \neq \ss' \leq \ss } \| \nabla_A \partial_x^{\ss - \ss'} ( \tfrac{f_1}{ M_{ \Lambda_0 } } ) \|_{L^2_{x,A} ( M_{ \Lambda_0 } ) } \\
		\no \geq & \tfrac{d}{2} \| \nabla_A \partial_x^{\ss} ( \tfrac{f_1}{ M_{ \Lambda_0 } } ) \|^2_{L^2_{x,A} ( M_{ \Lambda_0 } ) } \\
		& - C ( 1 + \| \nabla_x \Lambda_0 \|^{2s}_{ H_x^{ s + 1 } } ) \| \nabla_x \Lambda_0 \|^{s^2 + s}_{ H_x^{ s + 1 } } \sum_{0 \neq \ss' \leq \ss } \| \nabla_A \partial_x^{\ss - \ss'} ( \tfrac{f_1}{ M_{ \Lambda_0 } } ) \|^2_{L^2_{x,A} ( M_{ \Lambda_0 } ) }
	\end{align}
Moreover, considering the fact $\int_{\mathrm{SO}(3)} f_1 \d A = 0$ and following the inequality \eqref{HPoincare}, one has
\begin{equation}\label{f1-3}
	\begin{aligned}
		\lambda_0 \| \partial_x^{\ss} ( \tfrac{f_1}{ M_{ \Lambda_0 } } ) \|_{ L^2_{x,A} ( M_{ \Lambda_0 } ) }^2 \leq & \| \nabla_A \partial_x^{\ss} ( \tfrac{f_1}{ M_{ \Lambda_0 } } ) \|_{ L^2_{x,A} ( M_{ \Lambda_0 } ) }^2 \\
		& + C \| \nabla_x \Lambda_0 \|^{2s^2}_{H^{s+1}_x} \sum_{0 \neq \ss' \leq \ss} \| \nabla_A \partial_x^{\ss - \ss'} ( \tfrac{f_1}{ M_{ \Lambda_0 } } ) \|^2_{ L^2_{x,A} ( M_{ \Lambda_0 } ) } \,.
	\end{aligned}
\end{equation}
Furthermore, the H\"older inequality implies
\begin{equation}\label{f1-4}
	\begin{aligned}
		| \iint_{\R^3 \times \mathrm{SO(3)}} \partial_x^{\ss} h_0 \cdot \partial_x^{\ss} ( \tfrac{f_1}{ M_{ \Lambda_0 } } ) M_{ \Lambda_0 } \d A \d x | \leq \| \partial_x^{\ss} h_0 \|_{ L^2_{x,A} ( M_{ \Lambda_0 } ) } \| \partial_x^{\ss} ( \tfrac{f_1}{ M_{ \Lambda_0 } } ) \|_{ L^2_{x,A} ( M_{ \Lambda_0 } ) } \,.
	\end{aligned}
\end{equation}
Consequently, the relations \eqref{f1-1}-\eqref{f1-2}-\eqref{f1-3}-\eqref{f1-4} tell us that
\begin{equation}\label{Hder-f1}
	\begin{aligned}
		& \| \partial_x^{\ss} ( \tfrac{f_1}{ M_{ \Lambda_0 } } ) \|_{ L^2_{x,A} ( M_{ \Lambda_0 } ) }^2 + \| \nabla_A \partial_x^{\ss} ( \tfrac{f_1}{ M_{ \Lambda_0 } } ) \|^2_{L^2_{x,A} ( M_{ \Lambda_0 } ) } \\
		\lesssim & ( 1 + \| \nabla_x \Lambda_0 \|^{2s^2 + 3s}_{H^{s+2}_x} ) \sum_{0 \neq \ss' \leq \ss} \| \nabla_A \partial_x^{\ss - \ss'} ( \tfrac{f_1}{ M_{ \Lambda_0 } } ) \|^2_{ L^2_{x,A} ( M_{ \Lambda_0 } ) } + \| \partial_x^{\ss} h_0 \|_{ L^2_{x,A} ( M_{ \Lambda_0 } ) }
	\end{aligned}
\end{equation}
for $\ss \in \mathbb{N}^3$ with $1 \leq |\ss| = k \leq s$ ($s \geq 2$). Observe that Theorem \ref{LWP} shows
\begin{equation*}
	\begin{aligned}
		1 + \| \nabla_x \Lambda_0 \|^{2s^2 + 3s}_{H^{s+2}_x} \leq 1 + C_0 (\| \nabla_x \Lambda_0^{in} \|_{H^{s+2}_x}, \| \nabla_x \rho_0^{in} \|_{H^{s+2}_x}) : = \mathbb{C}_{s+2}^{in} > 0 \,.
	\end{aligned}
\end{equation*}
Then the bound \eqref{Hder-f1} reduces to
\begin{equation}\label{L2-k}
	\begin{aligned}
		& \sum_{|\ss| = k} \| \partial_x^{\ss} ( \tfrac{f_1}{ M_{ \Lambda_0 } } ) \|_{ L^2_{x,A} ( M_{ \Lambda_0 } ) }^2 + \sum_{|\ss| = k} \| \nabla_A \partial_x^{\ss} ( \tfrac{f_1}{ M_{ \Lambda_0 } } ) \|^2_{L^2_{x,A} ( M_{ \Lambda_0 } ) } \\
		\lesssim & \mathbb{C}_{s+1}^{in} \sum_{0 \leq \mathfrak{j} \leq k-1 } \sum_{ |\ss'| = \mathfrak{j} } \| \nabla_A \partial_x^{\ss'} ( \tfrac{f_1}{ M_{ \Lambda_0 } } ) \|^2_{ L^2_{x,A} ( M_{ \Lambda_0 } ) } + \| \partial_x^{\ss} h_0 \|^2_{ L^2_{x,A} ( M_{ \Lambda_0 } ) }
	\end{aligned}
\end{equation}

Based on the bounds \eqref{L2-0} and \eqref{L2-k}, the induction arguments for $k = 0, 1, \cdots, s$ shows that
\begin{equation}\label{f1-Hs-1}
	\begin{aligned}
		& \| \tfrac{f_1}{ M_{ \Lambda_0 } } \|_{ H^s_x L^2_A ( M_{ \Lambda_0 } ) }^2 + \| \nabla_A ( \tfrac{f_1}{ M_{ \Lambda_0 } } ) \|^2_{H^s_x L^2_A ( M_{ \Lambda_0 } ) } \\
		= & \sum_{|\ss| \leq s} \| \partial_x^{\ss} ( \tfrac{f_1}{ M_{ \Lambda_0 } } ) \|_{ L^2_{x,A} ( M_{ \Lambda_0 } ) }^2 + \sum_{|\ss| \leq s} \| \nabla_A \partial_x^{\ss} ( \tfrac{f_1}{ M_{ \Lambda_0 } } ) \|^2_{ L^2_{x,A} ( M_{ \Lambda_0 } ) } \\
		\leq & C ( \mathbb{C}_{s+2}^{in} ) \| h_0 \|^2_{ H^s_x L^2_A ( M_{ \Lambda_0 } ) } \,.
	\end{aligned}
\end{equation}
for $s \geq 2$. Recall that 
\begin{equation}\label{h0-express}
	\begin{aligned}
		h_0 = \tfrac{1}{ M_{ \Lambda_0 } } ( \partial_t + A e_1 \cdot \nabla_x ) f_0 = ( \partial_t + A e_1 \cdot \nabla_x ) \rho_0 + \tfrac{\nu_0}{d} \rho_0 \big[ A \cdot ( \partial_t + A e_1 \cdot \nabla_x ) \Lambda_0 \big]
	\end{aligned}
\end{equation}
where the last equality is derived from $f_0 = \rho_0 M_{ \Lambda_0 } $ and the relation \eqref{M-der}. It is thereby easy to see that
\begin{equation}\label{h0-bnd}
	\begin{aligned}
		\| h_0 \|^2_{ H^s_x L^2_A ( M_{ \Lambda_0 } ) } \leq C ( \| \partial_t ( \rho_0, \Lambda_0 ) \|^2_{H^s_x} + \| \nabla_x ( \rho_0, \Lambda_0 ) \|^2_{H^s_x} + \| \nabla_x ( \rho_0, \Lambda_0 ) \|^4_{H^s_x} ) \leq C ( \mathbb{C}_s^{in} )
	\end{aligned}
\end{equation}
by using Theorem \ref{LWP}. Then, by \eqref{f1-Hs-1} and \eqref{h0-bnd}, the first inequality in \eqref{f1-Rf1-bnd} holds.

It remains to control the quantity $ \| \tfrac{ R ( f_1 ) }{ M_{ \Lambda_0 } } \|_{H^s_x L^2_A ( M_{ \Lambda_0 } )} $. Recalling that $R (f_1) = - ( \partial_t + A e_1 \cdot \nabla_x ) f_1 $, one has
\begin{equation}
	\begin{aligned}
		\tfrac{ R ( f_1 ) }{ M_{ \Lambda_0 } } = - ( \partial_t + A e_1 \cdot \nabla_x ) ( \tfrac{ f_1 }{ M_{ \Lambda_0 } } ) - \tfrac{\nu_0}{d} ( \tfrac{ f_1 }{ M_{ \Lambda_0 } } ) [ A \cdot ( \partial_t + A e_1 \cdot \nabla_x ) \Lambda_0 ] \,,
	\end{aligned}
\end{equation}
where the relation \eqref{M-der} has been used. Then
\begin{equation}\label{Rf1-1}
	\begin{aligned}
		\| \tfrac{ R ( f_1 ) }{ M_{ \Lambda_0 } } \|_{H^s_x L^2_A ( M_{ \Lambda_0 } )} \lesssim \| \partial_t ( \tfrac{ f_1 }{ M_{ \Lambda_0 } } ) \|_{H^s_x L^2_A ( M_{ \Lambda_0 } )} + \| \nabla_x ( \tfrac{ f_1 }{ M_{ \Lambda_0 } } ) \|_{H^s_x L^2_A ( M_{ \Lambda_0 } )} \\
		+ \| \tfrac{ f_1 }{ M_{ \Lambda_0 } } \|_{H^s_x L^2_A ( M_{ \Lambda_0 } )} \big( \| \partial_t \Lambda_0 \|_{H^s_x} + \| \nabla_x \Lambda_0 \|_{H^s_x} \big) \,.
	\end{aligned}
\end{equation}
Together with the first inequality in \eqref{f1-Rf1-bnd} and Theorem \ref{LWP}, one easily has
\begin{equation}\label{Rf1-2}
	\begin{aligned}
		\| \nabla_x ( \tfrac{ f_1 }{ M_{ \Lambda_0 } } ) \|_{H^s_x L^2_A ( M_{ \Lambda_0 } )} + \| \tfrac{ f_1 }{ M_{ \Lambda_0 } } \|_{H^s_x L^2_A ( M_{ \Lambda_0 } )} \big( \| \partial_t \Lambda_0 \|_{H^s_x} + \| \nabla_x \Lambda_0 \|_{H^s_x} \big) \\
		\leq C ( \| \nabla_x \rho_0^{in} \|_{H^{s+3}_x}, \| \nabla_x \Lambda_0^{in} \|_{H^{s+3}_x} ) \,.
	\end{aligned}
\end{equation}

We next control the norm $ \| \partial_t ( \tfrac{ f_1 }{ M_{ \Lambda_0 } } ) \|_{H^s_x L^2_A ( M_{ \Lambda_0 } )} $. Applying the time derivative operator $\partial_t$ to \eqref{f1-eq}, we have
\begin{equation}
	\begin{aligned}
		\tfrac{d}{ M_{ \Lambda_0 } } \nabla_A \cdot \big[ M_{ \Lambda_0 } \nabla_A \partial_t ( \tfrac{f_1}{ M_{ \Lambda_0 } } ) \big] = & - \partial_t ( \tfrac{d}{ M_{ \Lambda_0 } } ) \nabla_A \cdot \big[ M_{ \Lambda_0 } \nabla_A ( \tfrac{f_1}{ M_{ \Lambda_0 } } ) \big] \\
		& - \tfrac{d}{ M_{ \Lambda_0 } } \nabla_A \cdot \big[ \partial_t M_{ \Lambda_0 } \nabla_A ( \tfrac{f_1}{ M_{ \Lambda_0 } } ) \big] + \partial_t h_0 \,.
	\end{aligned}
\end{equation}
Following the similar estimate of \eqref{f1-Hs-1}, one can derive that
\begin{equation}
	\begin{aligned}
		& \| \partial_t ( \tfrac{f_1}{ M_{ \Lambda_0 } } ) \|_{ H^s_x L^2_A ( M_{ \Lambda_0 } ) }^2 + \| \nabla_A \partial_t ( \tfrac{f_1}{ M_{ \Lambda_0 } } ) \|^2_{H^s_x L^2_A ( M_{ \Lambda_0 } ) } \\
		\leq & C ( \mathbb{C}_{s+2}^{in} ) \big( \| \tfrac{f_1}{ M_{ \Lambda_0 } } \|_{ H^s_x L^2_A ( M_{ \Lambda_0 } ) }^2 + \| \nabla_A ( \tfrac{f_1}{ M_{ \Lambda_0 } } ) \|^2_{H^s_x L^2_A ( M_{ \Lambda_0 } ) } + \| \partial_t h_0 \|^2_{ H^s_x L^2_A ( M_{ \Lambda_0 } ) } \big) \\
		\leq & C ( \mathbb{C}_{s+2}^{in} ) \big( \| h_0 \|^2_{ H^s_x L^2_A ( M_{ \Lambda_0 } ) } + \| \partial_t h_0 \|^2_{ H^s_x L^2_A ( M_{ \Lambda_0 } ) } \big) \\
		\leq & C ( \mathbb{C}_{s+2}^{in} ) \big( 1 + \| \partial_t h_0 \|^2_{ H^s_x L^2_A ( M_{ \Lambda_0 } ) } \big) \,,
	\end{aligned}
\end{equation}
where the bounds \eqref{f1-Hs-1} and \eqref{h0-bnd} have been utilized. By \eqref{h0-express}, one has
\begin{equation*}
	\begin{aligned}
		\partial_t h_0 = ( \partial_t + A e_1 \cdot \nabla_x ) \partial_t \rho_0 + \tfrac{\nu_0}{d} \partial_t \rho_0 \big[ A \cdot ( \partial_t + A e_1 \cdot \nabla_x ) \Lambda_0 \big] + \tfrac{\nu_0}{d} \rho_0 \big[ A \cdot ( \partial_t + A e_1 \cdot \nabla_x ) \partial_t \Lambda_0 \big] \,,
	\end{aligned}
\end{equation*}
which reduces to
\begin{equation}
	\begin{aligned}
		\| \partial_t h_0 \|^2_{ H^s_x L^2_A ( M_{ \Lambda_0 } ) } \lesssim & \| \partial_t^2 \rho_0 \|^2_{H^s_x} + (1 + \| \nabla_x \rho_0 \|^2_{H^s_x} ) \| \partial_t^2 \Lambda_0 \|^2_{H^s_x} \\
		& + \| \partial_t \rho_0 \|^2_{H^{s+1}_s} + (1 + \| \nabla_x \rho_0 \|^2_{H^s_x} ) \| \partial_t \Lambda_0 \|^2_{H^{s+1}_x} \\
		& + \partial_t \rho_0 \|^2_{H^s_x} ( \| \partial_t \Lambda_0 \|^2_{H^s_x} + \| \nabla_x \Lambda_0 \|^2_{H^s_x} ) \\
		\lesssim & C ( \| \nabla_x \rho_0^{in} \|_{H^{s+1}_x}, \| \nabla_x \Lambda_0^{in} \|_{H^{s+1}_x} ) ( \| \partial_t^2 \rho_0 \|^2_{H^s_x} + \| \partial_t^2 \Lambda_0 \|^2_{H^s_x} ) \,,
	\end{aligned}
\end{equation}
where the last inequality is derived from Theorem \ref{LWP}. Due to $(\rho_0, \Lambda_0)$ obeying the equation \eqref{SOHB} or \eqref{SOHB-1}, it is easy to know that $(\partial_t^2 \rho_0, \partial_t^2 \Lambda_0) \thicksim ( \nabla_x^2 \rho_0, \nabla_x^2 \Lambda_0 ) $. Together with the structure of \eqref{SOHB} or \eqref{SOHB-1}, one can derived that
\begin{equation}
	\begin{aligned}
		\| \partial_t^2 \rho_0 \|^2_{H^s_x} + \| \partial_t^2 \Lambda_0 \|^2_{H^s_x} \leq & C ( 1 + \| \nabla_x \rho_0 \|_{H^{s+2}_x}^4 + \| \nabla_x \Lambda_0 \|^4_{H^{s+2}_x} ) \\
		\leq & C ( \| \nabla_x \rho_0^{in} \|_{H^{s+2}_x}, \| \nabla_x \Lambda_0^{in} \|_{H^{s+2}_x} ) \,.
	\end{aligned}
\end{equation}
Collecting the above all estimates, we have
\begin{equation}\label{Rf1-3}
	\begin{aligned}
		\| \partial_t ( \tfrac{f_1}{ M_{ \Lambda_0 } } ) \|_{ H^s_x L^2_A ( M_{ \Lambda_0 } ) }^2 + \| \nabla_A \partial_t ( \tfrac{f_1}{ M_{ \Lambda_0 } } ) \|^2_{H^s_x L^2_A ( M_{ \Lambda_0 } ) } \leq C ( \| \nabla_x \rho_0^{in} \|_{H^{s+2}_x}, \| \nabla_x \Lambda_0^{in} \|_{H^{s+2}_x} ) \,.
	\end{aligned}
\end{equation}
As a result, the bounds \eqref{Rf1-1}, \eqref{Rf1-2} and \eqref{Rf1-3} complete the second inequality in \eqref{f1-Rf1-bnd}. Then the proof of Lemma \ref{Lmm-f1} is finished.


\section*{Acknowledgments}

This work was supported by National Key R\&D Program of China under the grant 2023YFA 1010300, the National Natural Science Foundation of China under contract No. 12201220, the Guang Dong Basic and Applied Basic Research Foundation under contract No. 2024A15150123 58, and the Fundamental Research Funds for the Central Universities under contract No. 531118011008.

\bigskip

\bibliography{reference}

\end{document}